\documentclass[11pt]{article}
\usepackage{graphicx}
\usepackage{amssymb,amsfonts,amsthm}
\usepackage{amsmath,amsthm,amssymb}
\usepackage{amsfonts}
\newcommand{\bbb}{{\cal B}}

\newcommand{\N}{{\mathbb N}}
\newcommand{\R}{{\mathbb R}}

\newtheorem{theorem}{Theorem}[section]
\newtheorem{lemma}[theorem]{Lemma}

\theoremstyle{definition}
\newtheorem{df}[theorem]{Definition}

\newtheorem{rk}[theorem]{Remark}

\newtheorem{ct}[theorem]{Quasi-Theorem}

\newcounter{ppp}

\newcommand{\la}{\langle}
\newcommand{\ra}{\rangle}

\newcommand{\area}{{\rm area}}
\newcommand{\me}{\medskip}

\newcommand{\Lab}{\phi}

\newcommand{\tool}{\stackrel{\ell}{\too} }

\newcommand{\ttt}{{\cal T}}

\newcommand{\eee}{{\cal D}}
\newcommand{\aaa}{{\cal A}}

\newcommand{\LL}{{\cal L}}

\newcommand{\bb}{{\cal B}}
\newcommand{\topp}{{\bf top}}
\newcommand{\ttopp}{{\bf ttop}}
\newcommand{\tbott}{{\bf tbot}}
\newcommand{\bott}{{\bf bot}}
\newcommand{\vk}{van Kampen }

\newcommand{\ccc}{{\cal C}}

\newcommand{\iv}{^{-1}}

\newcommand{\too}{\to }

\newcommand{\xxx}{{\cal X} }

\newcommand{\gc}{{\cal G} }

\newcommand{\qq}{{\cal Q} }
\newcommand{\QQ}{{\mathbf{Q}}}
\newcommand{\TT}{{\mathbf{T}}}
\newcommand{\sss}{{\cal S} }

\newcommand{\base}{\mathrm{base}}
\begin{document}

\renewcommand{\theequation}{\thesection.\arabic{equation}}

\title{Groups with small Dehn functions and bipartite chord diagrams}
 \author{A.Yu. Ol'shanskii, M.V. Sapir\thanks{Both authors were supported in part by
the NSF grants DMS 0245600 and DMS 0455881. In addition, the
research of the first author was supported in part by the Russian
Fund for Basic Research 05-01-00895,  the research of the second
author was supported in part by the NSF grant DMS 9978802 and the
US-Israeli BSF grant 1999298.}}
\date{}
\maketitle

\begin{abstract}
We introduce a new invariant of bipartite chord diagrams and use it
to construct the first examples of groups with Dehn function
$n^2\log n$. Some of these groups have undecidable conjugacy
problem. Our groups are multiple HNN extensions of free groups. We
show that $n^2\log n$ is the smallest Dehn function of a multiple
HNN extension of a free group with undecidable conjugacy problem.
\end{abstract}

\tableofcontents

\section{Introduction}

Recall that the {\em Dehn function} of a finite presentation $\la
X\mid R\ra$ of a group $G$ is the smallest function $f(n)$ such that
any word of length at most $n$ in $X$ that represents the identity
of $G$ is freely equal to a product of at most $f(n)$ conjugates of
elements of $R$. The Dehn functions $f_1, f_2$ of any two finite
presentations of the same group $G$ are {\em equivalent}, that is
$f_2(n)\le Cf_1(Cn) +Cn+C$, $f_1(n)<Cf_2(Cn)+Cn+C$ for some constant
$C$. As usual, we do not distinguish equivalent functions.

The purpose of this paper is to prove the following statement.

\begin{theorem}\label{thmain} There exist finitely presented
multiple HNN extensions of free groups with Dehn function $n^2\log
n$. Some of these groups have undecidable conjugacy problem.
Conversely, if $d(n)$ is the Dehn function of a multiple HNN
extension of a free group and $\lim_{n\to\infty}^c
\frac{d(n)}{n^2\log n}=0$ then the group has decidable conjugacy
problem. Here $\lim^c$ stands for the {\em constructive limit} (see
Definition \ref{dfcon} below).
\end{theorem}

\begin{rk}{\rm For those unfamiliar with the definition of constructive limit,
it is the same as the ordinary $(\epsilon,N)$-definition, only $N$
must recursively depend on $\epsilon$. Theorem \ref{thmain}
implies that a multiple HNN extension of a free group whose Dehn
function does not exceed, say, $n^2\sqrt{\log n}$ or
$n^2\frac{\log n}{\log\log\log...\log n}$, must have decidable
conjugacy problem.\footnote{The authors are able to show that one
cannot replace constructive limits by the ordinary limits in
Theorem \ref{thmain}. The proof is technically more difficult than
Theorem \ref{thmain} and is left out of this paper.}}
\end{rk}

In particular, this theorem gives the first example of a Dehn
function between $n^2$ and $n^4$ not of the form $n^\alpha$,
$\alpha\in\mathbb{R}$ (see \cite{Bri}, \cite{BrBr}). The set of
Dehn functions $\ge n^4$ is known to be large, and contains
functions of the form $n^\alpha$, $n^\alpha\log n$,
$n^\alpha\log\log n$ for any rational $\alpha\ge 4$ (and in fact
any relatively fast computable $\alpha>4$), and much more
complicated functions. An almost complete description of all Dehn
functions $\ge n^4$ has been found in \cite{SBR}. Unfortunately,
the methods from \cite{SBR} do not give any information about Dehn
functions between $n^2$ and $n^4$, and methods of \cite{Bri} and
\cite{BrBr} give Dehn functions only of the form $n^\alpha$. Thus
the picture of the class of small Dehn functions is incomplete.

Recall also that it is still unknown if there exists a finitely
presented group with Dehn function $n^2$ and undecidable conjugacy
problem. We believe that such groups do not exist and  $n^2\log n$
is the lowest Dehn function of a group with undecidable conjugacy
problem in a large class of groups than the HNN extensions of free
groups (see Section \ref{why}). Groups with subquadratic Dehn
functions are hyperbolic \cite{Gro}, \cite{Ol}, and so they have
solvable conjugacy problem \cite{Gro}.

Our groups belong to the class introduced earlier in \cite{SBR} by
the second author. These are the hub-free realizations of
$S$-machines (in \cite{SBR} these groups were denoted by
$G'_N(\sss)$). It follows from \cite[Lemma 8.1]{SBR} that the Dehn
function of any $G'_N(\sss)$ does not exceed $n^3$. E. Rips and the
second author conjectured that some of these groups have Dehn
functions $n^2\log n$ and showed how such groups could be
constructed. Unfortunately existing methods of finding upper bounds
of Dehn functions did not give precise upper bounds of Dehn function
of $G'_N(\sss)$. The reason is that $G'_N(\sss)$ does not have
distorted cyclic subgroups to apply methods from \cite{Bri},
\cite{BrBr}, and \vk diagrams over these groups do not have
hyperbolic structure which allows one to apply surgeries from
\cite{SBR} or \cite{BORS}.

In this paper, we use a new method of finding the upper bounds of
Dehn functions introduced by the first author. The method is based
on a  new invariant (dispersion) of bipartite chord diagrams
naturally associated with van Kampen diagrams over the
presentations of our groups. These invariants resemble the
invariants of cord diagrams studied in \cite{PV, Polyak} by Polyak
and Viro in relation to problems of Arnold, computing Vassiliev
invariants of knots, etc. As far as we know, this is the first use
of such invariants in geometric group theory, but we are convinced
that similar invariants will be applied to other problems
including those which are far from Dehn functions. The main idea
is that these invariants measure the complexity of a diagram, and
allow us to perform surgeries decreasing the complexity.

The paper is constructed as follows. In Section \ref{why}, we give a
``quasi-proof" of the fact that any finitely presented group with
undecidable conjugacy problem has Dehn function at least $n^2\log
n$, and a complete proof of this statement for multiple HNN
extensions of free groups.

In Section \ref{smach}, we start by introducing general properties
of $S$-machines viewed as groups (multiple HNN-extensions of free
groups). Thus we identify an $S$-machine and the corresponding
group. In particular, we introduce the standard notions of {\em
bands} and {\em trapezia} from \cite{SBR}.

Then we show how to slow down any $S$-machine $\sss$ so that the
space function of the new $S$-machine $\sss\circ Z$ becomes
equivalent to the logarithm of the time function. Later, in Section
\ref{upperb}, these properties will translate into the upper bound
$n^2\log n$ of the Dehn function of $\sss\circ Z$.

Several basic properties of the group $\sss\circ Z$ are proved in
Section \ref{propgr}. In \vk diagrams over an $S$-machine, there are
two types of bands, $Q$-bands and $\theta$-bands, that start and end
on the boundary of the diagram. These bands form a bipartite chord
diagram (BCD) since bands of the same type do not intersect and a
band of one type intersects a band of another type at most once.

In Section \ref{dispersion}, we introduce a new invariant of BCDs,
the {\em dispersion}, and prove that the dispersion of any BCD is
bounded from above by a quadratic polynomial in the number of chords
of one of the types. In Section \ref{upperb}, we prove that the area
of a \vk diagram $\Delta$ over $\sss\circ Z$ with
$|\partial\Delta|\le n$ does not exceed $C(n^2\log n+\eee)$ where
$\eee$ is the dispersion of the corresponding BCD.

This gives an upper bound of $n^2\log n$ for the Dehn function of
$\sss\circ Z$. In Section \ref{lower}, we give the similar lower
bound of the Dehn function and show that $\sss\circ Z$ has
undecidable conjugacy problem provided $\sss$ has undecidable
halting problem.  This completes the proof of Theorem \ref{thmain}.
Finally we show how to generalize our construction to obtain groups
with other unusual Dehn functions between $n^2$ and $n^3$.

{\bf Acknowledgement.} The authors are grateful to the referee for
many helpful remarks.

\section{Why $n^2\log n$?}
\label{why}

In this section, we shall give a ``quasi-proof" of the following
conjecture. Then we show that the conjecture is true for multiple
HNN extensions of free groups.

We shall need the well known constructive version of a limit of a
sequence of numbers. In fact we are going to use that definition
only in the case when the limit is $0$.

\begin{df} {\rm Let $g\colon \N\to \R$ be a function. We say that the
constructive limit of $g(n)$ as $n\to \infty$ is $0$ if for every
integer $A>0$ there exists $N=N(A)$ such that for every $n>N$,
$|f(n)|\le 1/A$, and the function $N(A)$ is recursive. In that case
we shall write $\lim_{n\to\infty}^c g(n)=0.$ It is easy to see that
$\lim_{n\to\infty}^c g(n)=0$ if and only if there exists an
increasing recursive function $f(n)$ such that $g(k)\le \frac{1}{n}$
for every $k\ge f(n)$.} \label{dfcon}
\end{df}

\begin{ct} \label{ct1}Let $d(n)$ be the Dehn
function of a finite group presentation $P$. Suppose that
$\lim_{n\to\infty}^c \frac{d(n)}{n^2\log n}=0$. Then $P$ has
decidable conjugacy problem.
\end{ct}

\proof We shall need the following Lemma. We call a \vk or annular
diagram over a group presentation {\em minimal} if it has minimal
area among all diagrams over that presentation with the same labels
of the contour (contours).

\begin{lemma}\label{lmnlogn} Let $\Delta$ be a minimal
annular diagram with contours $p, p'$ over a finite group
presentation $P$. Let $x$ be a shortest path connecting $p$ and
$p'$. Then the area of $\Delta$ is at least $C|x|\log |x|$ for some
constant $C$ depending on $P$.
\end{lemma}

\proof Consider the following construction. Let $p_0=p$ (considered
as a cyclic path) be the inner contour of the diagram $\Delta$.
Suppose that we have constructed a cyclic path $p_i$ surrounding the
hole of the diagram in $\Delta$ such that $p_i$ does not have common
vertices with $p'$. Let $K_i$ be the annulus bounded by $p_0$ and
$p_i$. Let $M_{i+1}$ be the set of cells of $\Delta$ outside $K_i$
that have common vertices with $p_i$. Then let $K_{i+1}$ be the
minimal annular subdiagram of $\Delta$ with simple contours that
contains $K_i$ and all cells from $M_{i+1}$. Let $p_{i+1}$ be the
outer contour of $K_i$ (the inner contour of $K_i$ is $p=p_0$).

It follows that every edge of the path $p_{i+1}$ belongs to the
contour of one of the cells of $M_{i+1}$. Hence every vertex of
$p_{i+1}$ can be connected with a vertex of $p_i$ by a path such
that

(0) the length of the path is bounded by a constant,

(1) it can be connected with $p_0$ by a path of length at most
\footnote{We use the Computer Science ``big-O" notation assuming
that $f(n)=O(g(n))$ if $\frac1Cg(n)<f(n)<Cg(n)$ for some positive
constant $C$.} $O(i)$ and

(2) the number of cells in $M_{i+1}$ is at least $O(|p_{i+1}|)$.

From (1), it follows that the number of subdiagrams $K_i$ is
$O(|x|)$. Furthermore, more than a half of the paths $p_i$ have
length at most $\log_c |x|$ where $c$ is, say, four times the number
of letters in the alphabet of the presentation $P$. Indeed,
otherwise we would have two paths $p_i$ and $p_j$, $i\ne j$ with the
same labels, and we could remove the annular subdiagram between
$p_i$ and $p_j$ reducing the area of $\Delta$ (that would contradict
the minimality of $\Delta$).

From (2), it follows that at least half of the subsets $M_i$ contain
at least $O(\log |x|)$ cells each. Since these sets do not
intersect, the number of cells in $\Delta$ is at least $O(|x|\log
|x|)$.
\endproof

Now the ``quasi-proof" of the Quasi-Theorem \ref{ct1} proceeds as
follows. Suppose that $P$ is a finite presentation with
undecidable conjugacy problem. Suppose that the constructive limit
of $\frac{d(n)}{n^2\log n}$ is $0.$ Then, in particular, $d(n)$ is
bounded from above by a recursive function, and $P$ has solvable
word problem.

Note that if an annular diagram $\Delta$ with contour labels $u$ and
$v$ has a simple path $x$ with label $t$ connecting the contours,
then we can cut $\Delta$ along $x$ and obtain a disc \vk diagram
with boundary label $t^{\pm 1}ut^{\mp 1}v\iv$. So if $|t|$ is
recursively bounded in terms of $|u|$ and $|v|$ (for every $u$ and
$v$ that are conjugate modulo $P$) then the conjugacy of $u$ and $v$
can be algorithmically verified.

Pick an increasing recursive function $f(n)$ with
$\frac{d(3k)}{Ck^2\log k}< \frac1n$ for every $k>f(n)$ where $C$
is the constant from Lemma \ref{lmnlogn} (as in Definition
\ref{dfcon}). Since the conjugacy problem for $P$ is undecidable,
there exists a minimal annular diagram $\Delta$ with contours $p$,
$p'$ such that any path in $\Delta$ connecting $p$ and $p'$ has
length at least $f(|p|+|p'|)$. Let $n=|p|+|p'|$. Let $x$ be a
shortest path connecting $p$ and $p'$. Thus \begin{equation}|x|\ge
f(n)\label{xn},\end{equation} and so

\begin{equation}\label{eqcon}\frac{d(3|x|)}{C|x|^2\log |x|}<\frac1n.\end{equation}

Since $x$ is a shortest path connecting $p$ and $p'$, $x$ is simple.
Let us cut $\Delta$ along $x$ and obtain a disc diagram $\Gamma$
with boundary label $zuz\iv v\iv$ where $z$ is the label of $x^{\pm
1}$. By Lemma \ref{lmnlogn}, the area of $\Gamma$ is at least
$C|x|\log|x|$.

Now we can take an integer $m$ between $\frac{|x|}{n}-1$ and
$\frac{|x|}{n}$. We attach $m$ copies of $\Gamma$ consecutively to
each other along the sides labeled by $z$ to get a \vk diagram
$\Pi$ with boundary label $zu^mz\iv v^{-m}$. Notice that $\Pi$ is
reduced because it covers $\Delta$ with multiplicity $m$ (after
identification of the two $z$-sides of its boundary): for any
2-cell subdiagram $\Sigma$ of $\Pi$ where the cells share an edge,
the annular diagram $\Delta$ contains a copy of $\Sigma$, so
$\Sigma$ cannot be reducible since $\Delta$ is minimal. The
perimeter $r$ of $\Pi$ is between $2|x|$ and $3|x|$, and the area
is $m$ times the area of $\Delta$. So, by Lemma \ref{lmnlogn}, the
area of $\Pi$ is at least $\frac{C|x|^2\log |x|}{n}$. By
(\ref{eqcon}), we can deduce that the area of $\Pi$ is bigger than
$d(r)$. This contradicts the definition of Dehn function of a
group presentation.
\endproof

\begin{rk}\label{rk91} The only gap in the preceding argument is contained in the last
phrase. Even though $\Pi$ is reduced, we cannot guarantee that
$\Pi$ has minimal area among all diagrams with the same boundary
label, and, in principal, the area of a minimal diagram with this
boundary label may be even quadratic in terms of the perimeter
$r$. Still we do not know any groups for which this proof does not
work. Note that we do not need $\Pi$ to be minimal: only that the
minimal diagram with the same boundary label does not have too few
cells compared to $\Pi$. Also we have freedom of choosing $u, v$,
$\Delta$ and $x$. We do not need $x$ to be a minimal length path
connecting the boundary components of $\Delta$. We only need that
the area of $\Delta$ exceeds $O(|x|\log |x|)$ divided by a
recursive function in $n$ (depending only on the presentation). In
addition, the number $m$ should only be $O(\frac{|x|}{|u|+|v|})$.
Thus Conjecture \ref{ct1} seems true for a very large class of
groups and possibly for all groups.
\end{rk}

Let $P$ be the standard presentation of a multiple HNN extension
of a free group $F_X$ with stable letters $t_1,...,t_k$ and pairs
of finitely generated associated subgroups $A_i=\la
a_{i,1},...,a_{i,j_i}\ra, B_i=\la b_{i,1},...,b_{i,j_i}\ra$ given
by their free generating sets. So the defining relations of the
presentation $P$ are $a_{i,s}^{t_i}=b_{i,s}$, $i=1,...,k$,
$s=1,...,j_i$. Here and below $a^t$ means $tat\iv$.

As usual when one works with HNN extensions, $t$-bands play
significant role (they are also called strips and corridors). We
shall give a more general definition of bands in Section
\ref{bands}.

For every letter $a$, an $a$-edge in a \vk diagram is an edge
labeled by $a^{\pm 1}$. A $t_i$-band in a diagram over $P$ is a
sequence of cells containing $t_i$-edges, such that every two
consecutive cells share a $t_i$-edge. It is well known \cite{MS}
that in a reduced \vk diagram over $P$, there are no $t_i$-annuli,
i.e. the first and the last $t_i$-edges of a $t_i$-band cannot
coincide. So every maximal $t_i$-band in a diagram must connect
two $t_i$-edges on the contour of the diagram. In an annular
diagram over $P$, every maximal $t_i$-band either connects two
edges belonging to the boundary or is an annulus surrounding the
hole of the diagram. The contour of a $t_i$-band has the form
$epf\iv q$ where $e, f$ are $t_i$-edges, and $e$ and $f$ do not
have $t$-edges (these are the {\em sides} of the band).

The following theorem is a part of Theorem \ref{thmain}.

\begin{theorem} \label{th56} Let $d(n)$ be the Dehn function of $P$
and $\lim^c_{n\to \infty}\frac{d(n)}{n^2\log n}=0$. Then $P$ has
decidable conjugacy problem.
\end{theorem}

\proof Suppose that the conjugacy problem is undecidable. We use the
same notation as in the quasi-proof above.

With every reduced \vk diagram $\Psi$ over $P$, one can associate a
chord diagram $C(\Psi)$ where the disc is the diagram and chords are
the $t$-bands (more precisely, their medians).

\begin{lemma}\label{uniquely1} Let $\Psi$ and $\Psi'$ be two reduced diagrams with
the same boundary, such that $C(\Psi)=C(\Psi')$. Then the areas of
$\Psi$ and $\Psi'$ are the same.
\end{lemma}

\proof Let $u$ be the common boundary label of $\Psi$ and $\Psi'$,
and $C=C(\Psi)=C(\Psi')$. Then there exists a one-to-one
correspondence $\ttt\mapsto \ttt'$ between the maximal $t$-bands in
$\Psi$ and $\Psi'$. A side of each maximal $t$-band $\ttt$ connects
two vertices on $\partial\Psi$. The label of the subpath of
$\partial\Psi$ connecting these vertices is a subword of the label
of $\partial\Psi$. Note that this subword is the same for the
corresponding side of $\ttt'$. Hence the labels $\ell, \ell'$ of the
corresponding sides of $\ttt$ and $\ttt'$ are equal modulo $P$.
Since both $\ell$ and $\ell'$ are words from the base group of an
HNN extension, they are equal in the base group as well. Since the
base group is free, these labels are freely equal. Since $a_{i,j}$
and $b_{i,j}$ freely generate subgroups $A_i$ and $B_i$, the number
of cells in $\ttt$ is the same as the number of cells in $\ttt'$
(both numbers are equal to the length of the element $\ell$ in the
corresponding subgroup $A_i$ or $B_i$, and the bands $\ttt$, $\ttt'$
are reduced diagrams by our assumption). Since every cell in $\Psi$
and $\Psi'$ belongs to a $t$-band, the areas of $\Psi$ and $\Psi'$
are the same.
\endproof

Recall that a {\em pinch} is a word of the form $t_iut_i\iv$ or
$t_i\iv vt_i$ where $u\in A_i$ and $v\in B_i$. Note that if the
area of a \vk diagram over $P$ is greater than $0$, then its
boundary label has at least two pinches since it is equal to 1
modulo $P$.

\begin{lemma} \label{uniquely}
Suppose the boundary label $l$ of $\Psi$ has only two pinches as a
cyclic word. Then $C(\Psi)$ is uniquely determined by $l$.
\end{lemma}

\proof The boundary $\partial\Psi$ is a product of two paths $pq\iv$
where the labels of $p$ and $q\iv$ do not have pinches. Hence every
$t$-band connects a $t$-edge of $p$ with a $t$-edge of $q$. Since
$t$-bands do not intersect, $C(\Psi)$ is reconstructed
uniquely.\endproof

We say that a word $W$ is {\em cyclically minimal} if none of the
cyclic shifts of it has pinches. Note that if $W$ is cyclically
minimal then any power of $W$ is cyclically minimal as well.

\begin{lemma} \label{onepinch} Suppose that a word $W$ is cyclically minimal.
Suppose also that a word $U$ has no pinches. Then the word $UWU\iv$
has at most one pinch.
\end{lemma}

\proof Indeed, if this word contains two pinches then $W=W_1W_2W_3$,
$U=U_1U_2$ where $U_2W_1$ is a pinch and $W_3U_2\iv$ is a pinch. But
then $W_3W_1$ is a pinch, and $W$ is not cyclically minimal.
\endproof

Lemma \ref{onepinch} immediately implies

\begin{lemma} \label{twopinches} Suppose that the words $W_1$ and $W_2$ are cyclically
minimal, and $U$ does not have pinches. Then the word $UW_1U\iv
W_2\iv$ has at most two pinches (as a cyclic word).
\end{lemma}

Now let us return to the {\em proof of Theorem \ref{th56}}. Note
that by Remark \ref{rk91}, we can do the following operations with
$\Delta$ and $x$:
\begin{itemize}
\item replace $\Delta$ by a minimal diagram $\Delta_1$ whose boundary
labels are equal to the boundary labels of $\Delta$ modulo $P$ and
have lengths that are recursively bounded in terms of $n=|u|+|v|$,
and
\item replace $x$ by a path $x'$ connecting the boundary components
such that $|x'|/|x|$ is recursively bounded in terms of $|u|+|v|$.
\end{itemize}
In order to be able to replace $(\Delta,x)$ by $(\Delta_1,x_1)$, one
needs to replace the condition (\ref{xn}) by the condition
$|x|>f(f_1(n))f_2(n)$ where $f_1$ and $f_2$ are some fixed
increasing recursive functions.

Our goal is to choose $\Delta$ and the path $x$ so that the diagram
$\Pi$ is minimal.

We may assume that the words $u$ and $v$ are cyclically minimal
because we can replace all pinches in these (cyclic) words by words
without $t$'s (this can be done effectively since the word problem
is solvable, and the lengths of $u$ and $v$ would increase only
recursively). Note that a $t$-band cannot connect two edges on the
same boundary component of $\Delta$ because otherwise $u$ or $v$
would contain pinches (as cyclic words). Thus there are two cases:
(1) $u$ and $v$ do not contain $t$-letters, and (2) $u$ contains
$t$-letters (then $v$ also contains $t$-letters).

In the first  case the maximal $t$-bands in $\Delta$ form annuli
surrounding the hole, the outer side of one annulus is the inner
side of the next one. In the second case, the maximal $t$-bands are
radial, connecting the inner contour with the outer boundary
component of $\Delta$.

Since $\Pi$ is reduced, Lemmas \ref{twopinches}, \ref{uniquely1} and
\ref{uniquely} tell us that we can claim that $\Pi$ is minimal
provided we can ensure that the label $z$ of $x$ does not contain
pinches.

{\bf Case I.} Suppose that $u$ and $v$ do not have $t$-letters. Let
$\tau$ be the number of $t$-annuli in $\Delta$.

Clearly no two sides of these $t$-annuli have the same labels
(otherwise we could remove the subdiagram bounded by these two
sides), which implies as in the proof of Lemma \ref{lmnlogn} that
the area of $\Delta$ is at least $O(\tau\log \tau)$. We can also
assume that $\tau>f(n)$ where $f(n)$ is the recursive function from
the quasi-proof. Indeed, if $\tau$ is bounded by a recursive
function for every $\Delta$ then the area of $\Delta$ is bounded by
a recursive function too.

Clearly there exists a path in $\Delta$ connecting $p$ and $p'$ and
having length $O(\tau)$. So we can assume that $x$ is that path (and
not the shortest path connecting $p$ and $p'$ as in the
quasi-proof).

Suppose that the label $z$ of $x$ contains a pinch and $x'$ is the
corresponding subpath of $x$. Then the first and the last edges of
$x'$ are $t_j$-edges belonging to two consecutive $t_j$-annuli
$\ttt$ and $\ttt'$ in $\Delta$. Since the label $z'$ of $x'$ is a
pinch, there exists a diagram $\Sigma$ consisting of one $t_j$-band
such that the label of $\partial\Sigma$ is $z'z''$, and $z''$ does
not contain $t$-letters. Therefore we can cut $\Delta$ along $x'$,
and patch the resulting hole by gluing in a copy of $\Sigma$ and a
copy of the mirror image $\Sigma'$ of $\Sigma$ glued together along
the part of the boundary labeled by $z''$. The resulting annular
diagram $\Delta'$ is not reduced. But instead of two $t_j$-bands
$\ttt$ and $\ttt'$, $\Delta'$ contains one $t_j$-band $\ttt''$ whose
set of cells is the union of the sets of sells in $\ttt$, $\ttt'$,
$\Sigma$ and $\Sigma'$. The annulus $\ttt''$ does not surround the
hole of $\Delta'$. Hence it bounds a disc subdiagram $\Phi$ of
$\Delta'$. The boundary label of $\Phi$ does not contain
$t$-letters. Hence the boundary label of $\Phi$ must be equal to 1
in the free group. Hence $\Phi$ can be replaced by a diagram without
cells. The new diagram $\Delta''$ has fewer cells, a contradiction
with the minimality of $\Delta$. Thus $z$ does not have pinches, and
we are done.

{\bf Case 2.} Suppose that $u$ and $v$ have $t$-letters.

The number $s$ of maximal $t$-bands in $\Delta$ is bounded by
$\min(|u|, |v|)\le n=|u|+|v|$. Let us cut $\Delta$ along a side
$q$ of a $t_i$-band. Let $\Gamma$ be the resulting diagram. We can
assume that $q$ is the shortest among the sides of the $t$-bands
in $\Delta$. We can also assume that $|q|$ is the smallest for all
annular diagrams with the same boundary labels and the same area.
This implies that there is no diagram with boundary label of the
form $\bar zu\bar z\iv v\iv$ whose area does not exceed the area
of $\Gamma$ and $|\bar z|<|q|$.

The \vk diagram $\Gamma$ has contour $pq(p')\iv (q')\iv$ where the
labels of $p$ and $q$ are (cyclic shifts of) $u$ and $v$
respectively, and $q$ and $q'$ have the same label $z$. We shall
denote by $q_-$ and $q_+$ the initial and terminal vertices of $q$.
Two vertices $V$ and $V'$ in $q$ and $q'$ are called {\em co-phase}
if their distances from $q_-$ and $q'_-$ (along $q$ and $q'$)
respectively are the same. We say that $V$ is {\em higher} ({\em
lower}) than $V'$ if the distance from $V$ to $q_-$ is bigger
(smaller) that the distance from $V'$ to $q'_-$.

Let $\ttt_1$, ..., $\ttt_s$ be all maximal $t$-bands  in $\Gamma$
connecting $p$ and $p'$, ordered from $q$ to $q'$. So $q$ is a side
of $\ttt_1$, $q'$ is a side of $\ttt_s$. It is easy to prove that
the length of each maximal $t$-band in $\Gamma$ is at most
$(|q|+n)\exp(Cn)$ for some constant $C$. Indeed, the label $w_i$ of
a side of $\ttt_i$ is equal modulo $P$ to the label of $q$
multiplied on the left and on the right by two words of length at
most $|u|+|v|$. Therefore the length of $w_i$ in the base group
cannot exceed $(|q|+n)\exp(C'n)$ for some constant $C'$: when we
reduce a pinch in a word, the length increases by a constant factor.
Since the number of cells in $\ttt_i$ is $O(|w_i|)$, we obtain the
desired inequality.

Since $|q|>>|u|$,  the area of $\Gamma$ (= the area of $\Delta_i$)
is at most $O(|q| \exp(C(|u|+|v|)))$ for some constant $C$. For
every vertex $V$ on $q$ consider a shortest path $p(V)$ connecting
$V$ with a vertex in $q'$. Note that $|p(V)|$ does not exceed a
constant times $s$. If the vertex $p(V)_+$ is co-phase with $V$ then
we say that $p(V)$ is {\em parallel} to $p$. We can assume that
$p(V_1)$ does not cross $p(V_2)$ for any $V_1\ne V_2$. Indeed, if
$p(V_1)=p_1p_2, p(V_2)=p_1'p_2'$ and the end points of $p_1$, $p_1'$
are the same, then the lengths of $p_2$ and $p_2'$ are the same
(otherwise the path $p(V_1)$ or $p(V_2)$ would not be shortest), and
we can replace $p(V_2)$ by $p_1'p_2$. Thus we can talk about a
subdiagram $\Phi(V_1, V_2)$ bounded by $p(V_1)$ and $p(V_2)$.

Suppose that there exist two vertices $V_1$ and $V_2$ such that
$p(V_1)$ and $p(V_2)$ are parallel to $p$ and the labels of the
paths $p(V_1)$ and $p(V_2)$ are the same. Then we can remove the
subdiagram $\Phi(V_1,V_2)$ of $\Gamma$. The resulting diagram
$\Gamma'$ will have boundary label of the form $\bar z u\bar z\iv v$
with $\bar z$ shorter than $|x|$, a contradiction. Hence the labels
of all paths $p(V)$ that are parallel to $p$ are different. Hence
the number $\pi$ of such paths is at most an exponent in $C|u|$ for
some constant $C$. Therefore $\pi$ is small comparing to $|q|$.
These paths cut $\Gamma$ into at most $e(n)$  conjugacy subdiagrams
where $e$ is a recursive function. One of these subdiagrams must
have area bigger than $C|q|/e(n)$ for some constant $C$. Thus we can
deal with this ``large" subdiagram instead of $\Gamma$. Hence
without loss of generality, we shall assume that $\Gamma$ does not
have paths $p(V)$ that are parallel to $p$ except possibly for $p$
and $p'$.

Now let us number all vertices of $q$ starting with $q_-$:
$V_1,V_2,...$. For each $j=1,2,...$, let $l(j)$ be the distance from
$p(V_j)_+$ to $q'_-$. Since the paths $p(V_j)$ do not intersect, and
none of the paths $p(V_j)$ are parallel to $p$ except possibly for
$p(q_-)$ and $p(q_+)$, either $l(j)\ge j$ for all $j$  or $j\ge
l(j)$ for all $j$ and the inequalities are strict except, possibly,
for $V_j=q_-$ and $V_j=q_+$. We can assume that the first
possibility holds because otherwise we can turn $\Gamma$ upside down
switching  $u$ and $v$.

Let us define a sequence of vertices on $q$ as follows: $Q_1=V_2$
and for every $j=2,3,...$ let $Q_j=V_{l(Q_{j-1})}$. Note that $Q_j$
is co-phase with $p(Q_{j-1})_+$. Let us define the path $p_j$ as the
composition of the subpath $Q_{j}-Q_{j-1}$ of $q$ and $p(Q_{j-1})$
is parallel to $p$. Then for some small enough constant $c$ the
number of paths $p_j$ of length $\le c\log |q|$ is at most
$\sqrt{|q|}$. Since the length of $p(Q_{j-1})$ is recursively
bounded in terms of $|u|$, we can assume that the number of $j$'s
such that the distance between $Q_{j-1}$ and $Q_j$ along $q$ is
smaller than $c\log |q|$ is at most $\sqrt{|q|}$. Since each $Q_j$
is higher than $Q_{j-1}$, the number of $j$'s for which this
distance is not smaller than $c\log |q|$ is at most $|q|/c\log |q|$.
Hence the total number of points $Q_j$ is at most
$O(\sqrt{|q|}+|q|/\log |q|)=O(|q|/\log |q|)$.

Let $q_0$ be the subpath of $q$ connecting $V_1$ and $V_2$ (this
path is simply an edge). Since $p(Q_{j-1})$ and $Q_j$ are co-phase,
images of paths $q_0$ and all $p(Q_j)$ in $\Delta$ is a path $x$
connecting the boundary components of $\Delta$. The length of $x$ is
at most $O(|q|/\log |q|)$ times a recursive function in $n$ as
required. Thus, by Lemmas \ref{uniquely1}-\ref{twopinches}, it
remains only to show that the label $z$ of $x$ does not have
pinches.

For a contradiction, suppose that $z$ contains a pinch $z'$ and
$x'$ is a subpath of $x$ whose label is $z'$. Then $z'=t_j^{\pm
1}u't_j^{\mp 1}$. The two $t$-bands intersecting $x'$ are
consecutive $t$-bands $\ttt_k$ and $\ttt_{k+1}$ for some $k$ in
the annular diagram $\Delta_i$ (we consider $k$ modulo $s$ so
$\ttt_{s+1}=\ttt_1$). Let us connect $x'_-$ and $x'_+$ with
vertices $R$, $R'$ on the inner contour of $\Delta_i$ along the
sides of $\ttt_k$ and $\ttt_{k+1}$ by paths $q_1$, $q_2$. The
vertices $R$, $R'$ can be connected by a subpath $\tilde p$ of the
boundary component such that $\tilde p$ contains the end edges of
the $t$-bands $\ttt_k$ and $\ttt_{k+1}$. If this quadrangle
surrounds the hole of $\Delta_i$, we can repeat the same
construction using another boundary component of $\Delta_i$. The
resulting quadrangle won't surround the hole in that case. Since
these cases are similar, we can assume that the initial quadrangle
does not surround the hole. Then this quadrangle is a \vk diagram
over $P$. Since the label of $x'$ is a pinch, it is equal to a
word without $t$-letters modulo $P$. Since the labels of $q_1$ and
$q_2$ do not contain $t$-letters, the label of $\tilde p$ is equal
modulo $P$ to a word without $t$-letters. Hence the (cyclic) word
$u_i$ contains a pinch, a contradiction.
\endproof

\section{$S$-machines}
\label{smach}

\subsection{$S$-machines as HNN extensions of free groups}

Probably the easiest way to view an $S$-machine $S$ in the sense of
\cite{SBR} is to consider $S$ as a group that is an HNN extension of
a free group $F(Q,Y)$ generated by two sets of letters: state
letters $Q=\cup_{i=1}^N Q_i$ and tape letters $Y=\cup_{i=1}^{N-1}
Y_i$ where $Q_i$ are disjoint and non-empty. The sets $Q_i$ (resp.
$Y_i$) are called {\em parts} of $Q$ (resp. $Y$).

We shall follow the tradition of calling state letters {\em
$q$-letters} and tape letters {\em $a$-letters}, even though we
shall use $k$ with indexes for state letters and $y$ with indexes
for tape letters.

Instead of the set of stable letters we have a collection $\Theta$
of $N$-tuples of  $\theta$-letters or {\em rules}. The components
of $\theta$ are called {\em brothers} $\theta_1,...,\theta_N$. In
this paper, we always assume that all brothers are different. We
set $\theta_{N+1}=\theta_1$, $Y_0=Y_N=\emptyset$.

To every $\theta\in \Theta$, we associate two sequences of
elements in $F(Q\cup Y)$: $B(\theta)=[U_1,...,U_N]$,
$T(\theta)=[V_1,...,V_N]$, and a subset $Y(\theta)=\cup
Y_i(\theta)$ of $Y$, where $Y_i(\theta)\subseteq Y_i$.

The words $U_i, V_i$ satisfy the following restriction:

\begin{itemize}
\item[(*)] For every $i=1,...,N$, the words $U_i$ and $V_i$ have the form
$$U_i=v_{i-1}k_iu_i, \quad V_i=v_{i-1}'k_{i}'u_{i}'$$ where $k_{i}, k_{i}'\in
Q_{i}$, $u_{i}$ and $u_{i}'$ are words in the alphabet $Y_{i}^{\pm
1}$, $v_{i-1}$ and $v_{i-1}'$ are words in the alphabet
$Y_{i-1}^{\pm 1}$.
\end{itemize}

The generating set $\xxx$ of $S$ consists of all $q$-, $a$- and
$\theta$-letters. The relations are:

$$U_i\theta_{i+1}=\theta_i V_i,\,\,\,\, i=1,...,s, \qquad \theta_j a=a\theta_j$$
for all $a\in Y_j(\theta)$. The first type of relations will be
called $(q,\theta)$-{\em relations}, the second type -
$(a,\theta)$-{\em relations}.

Sometimes we will denote the rule $\theta$ by $[U_1\to
V_1,...,U_N\to V_N]$. This notation contains all the necessary
information about the rule except for the sets $Y_i(\theta)$. In
most cases it will be clear what these sets are. In the $S$-machines
used in this paper, the sets $Y_i(\theta)$ will be mostly equal to
either $Y_i$ or $\emptyset$. By default $Y_i(\theta)=Y_i$.

In order to simplify the notation, we will use the notation
$v_ik_iu_i\tool v_i'k_i'u_i'$ for a part of a rule when the
corresponding $Y_i(\theta)$ is empty (a similar notation has been
used in \cite{OSamen}).

Every $S$-rule $\theta=[U_1\to V_1,...,U_s\to V_s]$ has an inverse
$\theta\iv=[V_1\to U_1,...,V_s\to U_s]$; we set
$Y_i(\theta\iv)=Y_i(\theta)$. We always divide the set of rules
$\Theta$ of an $S$-machine into two disjoint parts, $\Theta^+$ and
$\Theta^-$ such that for every $\theta\in \Theta^+$, $\theta\iv\in
\Theta^-$ and for every $\theta\in\Theta^-$, $\theta\iv\in\Theta^+$.
The rules from $\Theta^+$ (resp. $\Theta^-$) are called {\em
positive} (resp. {\em negative}).

\begin{rk} {\rm 1. Every $S$-machine is indeed an HNN-extension of
the free group $F(Y,Q)$ with finitely generated associated
subgroups. The stable letters are $\theta_1$ for every $\theta\in
\Theta$. We leave it as an exercise to find the associated
subgroups.
\medskip

2. Notice that in \cite{SBR}, a slightly different notation for
rules of $S$-machines was used. Instead of, say, the rule
\begin{equation}\label{rule1}
[v_0k_1u_1\tool v_0'k_1'u_1', v_1k_2u_2\to
v_1'k_2'u_2'],\end{equation} one would use the notation
\begin{equation}\label{rule2} [v_0k_1u_1v_1k_2u_2\to
v_0'k_1'u_1'v_1'k_2'v_2'].\end{equation} But two relations
corresponding to the rule (\ref{rule1}): $$\theta_1\iv
v_0k_1u_1=v_0'k_1'u_1'\theta_2\iv$$ and
$$\theta_2\iv v_1k_2u_2=v_1'k_2'u_2'\theta_3\iv$$ are Tietze equivalent to one
relation corresponding to (\ref{rule2}):
$$\theta_1\iv v_0k_1u_1v_1k_2u_2=v_0'k_1'u_1'v_1'k_2'u_2'\theta_3\iv$$
since
$\theta_2$ is expressible in terms the other generators. Therefore
these notations are equivalent.}
\end{rk}

\subsection{Bands}
\label{bands} From now on, we shall only consider {\em reduced} van
Kampen diagrams i.e. diagrams that do not contain cells that have a
common edge and are mirror images of each other. Hence all \vk
diagrams are assumed to be reduced. To study van Kampen diagrams
over the group $S$ we shall use bands and trapezia as in \cite{SBR},
\cite{BORS}, etc.

Here we repeat necessary definitions from \cite{SBR}.

\begin{df}Let $M$ be a subset of ${\cal X}$. An
$M$-band $\bb$ is a sequence of cells $\pi_1,...,\pi_n$ in a \vk
diagram such that

\begin{itemize}
\item Each two consecutive cells $\pi_i$ and $\pi_{i+1}$
in this sequence have a common edge $e_i$ labeled by a letter from
$M$.
\item Each cell $\pi_i$, $i=1,...,n$ has exactly two $M$-edges,
$e_{i-1}$ and $e_i$ (i.e. edges labeled by a letter from $M$).

\item If $n=0$, then $\bb$ is just an $M$-edge.
\end{itemize}
\end{df}

The counterclockwise boundary of the subdiagram formed by the cells
$\pi_1,...,\pi_n$ of $\bb$ has the form $e\iv q_1f q_2\iv$ where
$e=e_0$ is an $M$-edge of $\pi_1$, $f=e_n$ is an $M$-edge of
$\pi_n$. We call $q_1$ the {\em bottom} of $\bb$ and $q_2$ the {\em
top} of $\bb$, denoted $\bott(\bb)$ and $\topp(\bb)$.

Consider lines $l(\pi_i, e_i)$ and $l(\pi_i, e_{i-1})$ connecting a
point inside the cell $\pi_i$ with midpoints of the $M$-edges of
$\pi_i$. The broken line formed by the lines $l(\pi_1,e)$,
$l(\pi_i,e_i)$, $l(\pi_i,e_{i-1})$, $l(\pi_n,f)$ is called the {\em
median} of the band $\bb$. It connects the midpoints of $e$ and $f$
and lies inside the union of $\pi_i$. The $M$-edges $e$ and $f$ are
called the {\em start} and {\em end} edges of the band. If $n=0$,
then the median is the midpoint of $e=f$.

We say that an $M_1$-band and an $M_2$-band {\em cross} if their
medians cross. We say that a band is an {\em annulus} if its start
and end edges coincide. In this case the median of the band is a
simple closed curve.

The subdiagram bounded by the median of an annulus in a disc
diagram is called the {\em inside diagram of this annulus}.

Let $M_1$ and $M_2$ be two disjoint sets of letters, let ($\pi$,
$\pi_1$, \ldots, $\pi_n$, $\pi'$) be an $M_1$-band and let ($\pi$,
$\gamma_1$, \ldots, $\gamma_m$, $\pi'$) be an $M_2$-band. Suppose
that:
\begin{itemize}
\item the medians of these bands intersect in two points $A,B$ inside $\pi$ and $\pi'$,
and the parts of the medians between $A$ and $B$ form a simple
closed curve,
\item on the boundary  of $\pi$ and on the boundary of $\pi'$ the pairs of
$M_1$-edges separate the pairs of $M_2$-edges,
\item the start and end edges of these bands are not contained in the
region bounded by the medians of the bands.
\end{itemize}
Then we say that these bands form an {\em $(M_1,M_2)$-annulus} and
the simple closed curve formed by the parts of medians of these
bands is the {\em median} of this annulus. For every annulus we
define the {\em inside} subdiagram of the annulus as the subdiagram
bounded by the median of the annulus.

We shall call an $M$-band {\em maximal} if it is not contained in
any other $M$-band.

As in \cite{SBR}, we can consider $q$-bands where $M$ is one of
the sets $Q_i$, $\theta$-bands for every $\theta\in\Theta$, and
$a$-bands where $M=\{a\}\subseteq Y$. Every cel of a $q$-band is a
$q$-cell by definition. A $q$-cell is also a $(q,\theta)$-cell
(and also $(\theta, q)$-cell) since it corresponds to a
$(q,\theta)$-relation. The convention is that $a$-bands do not
contain $q$-cells, and so they consist of $(a,\theta)$-cells ($=$
$(\theta,a)$-cells) only.

The following lemma has been essentially proved in \cite{SBR}.

\begin{lemma}\label{NoAnnul}
A reduced van Kampen diagram $\Delta$ over $S$ has no $q$-annuli,
$\theta$-annuli, $(q,\theta)$-annuli,  $a$-annuli,
$(a,\theta)$-annuli.
\end{lemma}
\proof We assume that $\Delta$ is a counterexample with minimal
area. This means in particular that  the boundary of $\Delta$ is the
boundary component of an annulus $\aaa$, where $\aaa$ has one of the
types from the formulation of the lemma.

(1) Suppose $\aaa$ is a $q$-annulus. Then it consists of
$(q,\theta)$-cells. Hence there is a maximal $\theta$-band $\ttt$
in $\Delta$, whose first cell $\Pi_1$ and the last cell $\Pi_2$
belong to $\aaa$. Being members of the same $q$-band $\aaa$ and
$\theta$-band $\ttt$, the cells $\Pi_1$ and $\Pi_2$ cannot be
neighbors in $\ttt$ (the diagram is reduced). Hence $\ttt$ and a
part of $\aaa$ form a $(\theta,q)$-annulus whose area is smaller
than that of $\Delta$. This contradicts the choice of $\Delta$.

(2) Suppose $\aaa$ is a $\theta$-annulus. If it contains
$q$-cells, then we come to a contradiction as in (1). Otherwise
$\Delta$ has no $q$-cells since there is no counter-example of
smaller area. The inner part of $\aaa$ has no $\theta$-edges for
the same reason. Hence $\Gamma$ has no cells corresponding to the
relations of the group $S$. So, on the one hand, the inner label
of $\aaa$ is a cyclically reduced non-empty word in $Y$ since
$\Delta$ is a reduced diagram, and on the other hand, this word is
freely equal to 1, a contradiction.

(3) Suppose $\aaa$ is a $(q,\theta)$-annulus. Then the maximal
$q$-band $\ttt$ of $\aaa$ cannot have more than two cells because
otherwise $\Delta$ would contain a smaller counterexample as in
(2). Hence the length of $\ttt$ is 2, and its cells are mirror
copies of each other, a contradiction (we assumed that the
diagrams are reduced).

(4) Suppose $\aaa$ is an $a$-annulus. Then its boundary labels are
words in $\theta$-letter. This leads, as in (1), to a smaller
$(a,\theta)$-annulus, a contradiction.

(5) Suppose $\aaa$ is a $(\theta,a)$-annulus and let $\ttt$ be the
maximal $a$-band of it. It cannot have more than two
$(\theta,a)$-cells because otherwise there would be a smaller
$(\theta,a)$-annulus. Hence the length of $\ttt$ is 2, and its
cells are mirror copies of each other, a contradiction. The lemma
is proved.
\endproof

\subsection{Trapezia}

If $W=x_1...x_n$ is a word in an alphabet $X$, $Y$ is another
alphabet, and $\phi\colon X\to Y\cup\{1\}$ (where $1$ is the empty
word) is a map, then $\phi(W)=\phi(x_1)...\phi(x_n)$ is called the
{\em projection} of $W$ onto $Y$. We shall consider the projections
of words from $S$ onto $\Theta$ (all $\theta$-letters map to the
corresponding element of $\Theta$, all other letters map to $1$),
and the projection onto the alphabet $\{Q_1,...,Q_n\}$ (every
$q$-letter maps to the corresponding $Q_i$, other letters map to
$1$).

\begin{df}\label{dfsides}{\rm  The projection of the label
of a side of a $q$-band onto the alphabet $\Theta^{\pm 1}$ is
called the {\em history} of the band. The projection of the label
of a side of a $\theta$-band onto the alphabet $\{Q_1,...,Q_n\}$
is called the {\em base} of the band. Similarly we can define the
history of a word and the base of a word. The base of a word $W$
is denoted by $\base(W)$}. It will be convenient instead of
letters $Q_1, ...,Q_n$, in base words, to use representatives of
these sets. For example, if $k_1\in Q_1$, $k_2\in Q_2$, we shall
say that the word $k_1ak_2$ has base $k_1k_2$ instead of $Q_1Q_2$.
\end{df}

\begin{df}\label{dftrap}
{\rm Let $\Delta$ be a reduced \vk diagram which has the contour of
the form $p_1\iv q_1p_2q_2\iv$ where:
\medskip

$(TR_1)$ $p_1$ and $p_2$ are sides of $q$-bands,

\medskip

$(TR_2)$ $q_1$, $q_2$ are maximal parts of the sides of
$\theta$-bands such that $\Lab(q_1)$, $\Lab(q_2)$ start and end with
$q$-letters,

\medskip

$(TR_3)$ for every $\theta$-band $\ttt$ in $\Delta$, the labels of
$\topp(\ttt)$ and $\bott(\ttt)$ are reduced.

\medskip

Then $\Delta$ is called a {\em trapezium}. The path $q_1$ is called
the {\em bottom}, the path $q_2$ is called the {\em top} of the
trapezium, the paths $p_1$ and $p_2$ are called the {\em left and
right sides} of the trapezium. The history of the $q$-band whose
side is $p_1$ is called the {\em history} of the trapezium; the
length of the history is called the {\em height} of the trapezium.
The base of $p_1$ is called the {\em base} of the trapezium.}
\end{df}

\begin{rk}{\rm Property $(TR_3)$ is easy to achieve: by folding edges
with the same labels having the same initial vertex, one can make
the boundary label of any subdiagram in a \vk diagram reduced, see
\cite{SBR}.}\end{rk}

\begin{rk}\label{rk9} {\rm Notice that the top (bottom) side of a
$\theta$-band $\ttt$ does not necessarily coincide with the top
(bottom) side of the corresponding trapezium of height 1, and is
obtained from $\topp(\ttt)$ (resp. $\bott(\ttt)$) by trimming a few
first and last $a$-letters. We shall denote the trimmed top and
bottom sides of $\ttt$ by $\ttopp(\ttt)$ and
$\tbott(\ttt)$.}\end{rk}

\subsection{Admissible words and computations}

Using trapezia, one can now formally define admissible words and
application of a rule to a word.

\begin{df}
Let $\Delta$ be a trapezium of height 1. Let $\theta$ be the
element of $\Theta^{\pm 1}$ whose representative (one of the
brothers) is written on a side of that trapezium. Then the word
$W$ written on the bottom of of $\Delta$ is called an {\em
admissible word} for $\theta$. The word written on the top of
$\Delta$ is called the {\em result of application} of $\theta$ to
$W$ and is denoted by $\theta\cdot W$. Clearly, $\theta\cdot W$ is
uniquely determined by $\theta$ and $W$, and $\theta\iv \cdot
(\theta\cdot W)=W$ (hence the label of the top of the trapezium is
an admissible word for $\theta\iv$). We call a word {\em
admissible} if it is admissible for some $\theta\in\Theta$. For
every word $f=f_1f_2...f_n$ in $\Theta$ we define $f\cdot W$ as
$f_n\cdot(...(f_2\cdot(f_1\cdot W))...)$. In particular, $1\cdot
W=W$ for every word $W$.
\end{df}

Recall that we assume that all components in every $N$-tuple
$\theta$ are different. In order for a rule $\theta$ to be
applicable to a word $W=q_1u_1q_2...u_{n-1}q_n$ where $q_i$ are from
$Q_{j(i)}^{\pm 1}$, $u_i$ is a word in $Y^{\pm 1}$, the
$(q,\theta)$-relations involving $q_i$ and $q_{i+1}$ must share a
$\theta$-brother, and that $\theta$-brother must commute with all
letters of $u_i$ that are not involved in the relations containing
$q_i, q_{i+1}$. This and the definition of admissible word
immediately imply the following lemma.

\begin{lemma}\label{lmad} Every admissible word of a rule $\theta$ has the form
$q_1u_1q_2...u_{n-1}q_n$ where for every $i$ from $1$ to $n$ there
exists $j(i)$ such that $q_i\in Q_{j(i)}^{\pm 1}$, and
\begin{itemize}
\item If $q_i\in Q_{j(i)}$ then $u_i$ is a group word in
$Y_{j(i)}$ and $q_{i+1}\in Q_{j(i)+1}\cup Q_{j(i)}\iv$;

\item If $q_i\in Q_{j(i)}\iv$ then $u_i$ is a group word in
$Y_{j(i)-1}$ and $q_{i+1}\in Q_{j(i)}\cup Q_{j(i)-1}\iv$.
\end{itemize}
\end{lemma}

By Lemma \ref{NoAnnul}, any trapezium $\Delta$ of height $h\ge 1$
can be decomposed into $\theta$-bands $\ttt_1,...,\ttt_h$ connecting
the left and the right sides of the trapezium. The word written on
the trimmed top side of one of the bands $\ttt_i$ is the same as the
word written on the trimmed bottom side of $\ttt_{i+1}$,
$i=1,...,h$. Therefore with every trapezium $\Delta$ we can
associate a sequence of words $W_1, W_2,...,W_{h+1}$ and a sequence
of rules $\theta_1,\theta_2,...,\theta_h$ such that
$W_2=\theta_1\cdot W_1, W_3=\theta_2\cdot
W_2,...,W_{h+1}=\theta_h\cdot W_h$.

This pair of sequences will be called a {\em computation} of $S$
connecting $W_1$ and $W_{h+1}$.  We shall denote the computation by
$$W_1\to_{\theta_1}\dots\to_{\theta_{h-1}} W_h\to _{\theta_{h}}W_{h+1},$$ or simply
$$W_1\to \dots\to W_{h}\to W_{h+1}.$$

The number $h$ is called the {\em length} of the computation. Since
we consider only reduced diagrams, the history of every trapezium is
a reduced word. That word $t=\theta_1\theta_2...\theta_h$ is called
the  {\em history of computation}. The area of the trapezium is
called the {\em area of computation}. The length of the longest word
$W_i$ in this computation is called its {\em width}. It is also
convenient to consider {\em empty computations} consisting of one
word $W$. The history of an empty computation is the empty word, the
start and end words of this computation are equal to $W$.

Notice that $W_1f'=fW_{h+1}$ for some words $f, f'$ in $\theta$- and
$a$-letters whose projections onto $\Theta$ are equal to the history
word $t$. It is easy to see that $|f|=O(|t|)$.


\begin{rk}\label{rk87} One can easily see that the computation $W\to W_1\to... $ looks
like a computation of a Turing machine with many heads, the
$q$-letters. Heads can move left and right, change their states, and
change $a$-letters (which play the role of tape letters) around
them.
\end{rk}

As for a usual Turing machine, we choose a distinguished {\em stop
word} $\tilde W$ from $F(Q,Y)$.

We say that a word $W\in F(Q,Y)$ is accepted if there exists a
computation connecting this word and $\tilde W$.

The following lemma immediately follows from the definition of a
computation and Lemma \ref{NoAnnul}.

\begin{lemma} \label{computation}
Let $\Delta$ be a trapezium with bottom label $W$ and top label
$W'$. Then there is a unique computation $W\to...\to W'$ whose
history is the history of $\Delta$.
\end{lemma}

We shall also need the following lemma.

\begin{lemma}\label{lmcut}
Let $W_0\to_{\theta_1}W_1\to...\to_{\theta_t}W_t$ be a computation
of an $S$-machine $S$. Suppose that $W_i=W_j$ and
$\theta_i\ne\theta_{j+1}\iv$ for some $i,j$, $1\le i<j<t$. Then
\begin{equation}\label{cut} W_0\to_{\theta_1}W_1...\to_{\theta_i}
W_i\to_{\theta_{j+1}}...\to_{\theta_t} W_t\end{equation} is again a
computation of $S$.
\end{lemma}

\proof This surgery amounts to removing $\theta$-bands number
$i+1,...,j$ (counted from the bottom) in the trapezium $\Delta$
corresponding to the initial computation. Let us show that the new
diagram $\Delta'$ is reduced. Indeed, pairs of cells in the same
$\theta$-band in $\Delta'$ cannot cancel because the same pair of
cells existed in $\Delta$. The cells from different $\theta$-bands
in $\Delta'$ cannot cancel because either the same pair of cells
exists in $\Delta$ or one of these cells is in the $i$-th
$\theta$-band of $\Delta$, and the other one is in the $j+1$-st
$\theta$-band of $\Delta$, and these cells cannot cancel in
$\Delta'$ because $\theta_i\ne\theta_{j+1}\iv$ by our assumption.
The conditions $(TR_1), (TR_2), (TR_3)$ obviously hold for
$\Delta'$. So $\Delta'$ is a trapezium, and (\ref{cut}) is a
computation.
\endproof

\begin{lemma} \label{obratnaja}
Assume that two admissible words $W$ and $W'$ are conjugate in the
group $S$. Then there exists a computation $W\to...\to W''$ of $S$
where $W''$ is a cyclic conjugate of $W'$ that starts with a
$q$-letter from the same $Q_i$ as the first letter of $W$.
\end{lemma}

\proof By the van Kampen - Schupp lemma, there is a reduced
annular
 diagram $\Delta$ whose boundary components $p$ and $p'$ are
 clockwise labeled by $W$ and  $W'$. It follows from Lemma
 \ref{NoAnnul} that, for some $h\ge 0$, $\Delta$ is a union of
concentric $\theta$-annuli $\ttt_1,\dots,\ttt_h$, where $\ttt_1$ is
attached to
 $p$, $\ttt_2$ has a common boundary component with $\ttt_1$,..., $\ttt_h$
 is attached to $p'$. Let $k$ be the first letter of $W$. Then a
 maximal $k$-band connects $p$ and $p'$. We may assume that $h>0$.
 Cutting $\Delta$ along a
 side of this $k$-band, we get a trapezium. By Lemma
 \ref{computation}, there is a computation of $S$ connecting  $W$ with a word $W''$
 that is a cyclic shift
 of $W'$. Since all $q$-edges of a $q$-band have labels from the
 same  set $Q_i$, the first letters of $W$ and $W'$ are from the
 same $Q_i$.
\endproof

\begin{rk} \label{rk99} {\rm Suppose that the base of a trapezium
$\Delta$ starts with a $Q_1$-letter and ends with a $Q_N$-letter.
Then the labels of the sides of the trapezium are the same and do
not contain $a$-letters: it follows from the agreement that
$Q_{N+1}=Q_1$, and that in every part of an $S$-rule of the form
$v_0k_1u_1\to v_0'k_1u_1'$ (resp. $v_{N-1}k_Nu_N\to
v_{N-1}'k_N'u_N'$) the words $v_0,v_0'$ (resp. $u_N$, $u_N'$) are
empty since they are words over empty alphabets $Y_0$ and $Y_N$.}
\end{rk}

\begin{lemma} \label{conj} Suppose that the stop word $\tilde W$ starts with a letter
from $Q_1$ and ends with a letter from $Q_N$ and has only one
$Q_1$-letter. Suppose that the language of accepted words is not
recursive. Then the set of words that are conjugates of $\tilde W$
in $S$ is not recursive. Hence $S$ has undecidable conjugacy
problem.
\end{lemma}

\proof Since the first letter in $\tilde W$ is from $Q_1$ and the
last letter is from $Q_N$, the left and right sides of any trapezium
with the bottom label $\tilde W$ are the same by Remark \ref{rk99}.
Hence if $W$ is accepted, it is a conjugate of $\tilde W$.
Conversely, if $W$ starts with a $Q_1$-letter and is a conjugate of
$\tilde W$, then by Lemma \ref{obratnaja}, there exists a
computation $W\to...\to \tilde W'$ where $\tilde W'$ is a cyclic
conjugate of $\tilde W$ starting with a $Q_1$-letter. Since $\tilde
W$ contains only one $Q_1$-letter, $\tilde W'=\tilde W$, so $W$ is
accepted.
\endproof

\subsection{A slight modification of the $S$-machine from \cite{SBR}}
\label{secsm}

Let $\LL$ be a recursively enumerable language over an alphabet $X$.
Then by \cite[Proposition 4.1]{SBR}, there exists an $S$-machine
$\sss$ {\em recognizing} $\LL$ in the sense of Lemma \ref{lm1}
below.

For that $S$-machine, $N=17$, so the set of $q$-letters is
partitioned into $17$ subsets which will be more convenient to
denote by $K_1,...,K_N$. The elements of $K_i$ will be denoted by
$k_i(j)$. The stop word $\tilde W$ is $k_1(0)k_2(0)...k_N(0)$. The
set $X$ is contained in $Y_1$, and for every positive word $u$ in
$X$ we denote $\sigma(u)=k_1(1)uk_2(1)...k_N(1)$. The following
lemma is proved in \cite{SBR}.

\begin{lemma}\cite[Proposition 4.1]{SBR} (a) For every $u\in \LL$
there exists a computation $\sigma(u)\to_{\theta_1} \cdot \to \tilde
W$ consisting of positive words.

(b) For every word $u\not\in \LL$ over $X$ the word $\sigma(u)$ is
not accepted.
 \label{lm1}\end{lemma}

Let us modify $\sss$ a little. We add two new sets of state letters
$K_{N+1}=\{k_{N+1}\}, K_{N+2}=\{k_{N+2}\}$, and a new state letter
$\bar k_j$ in every $K_j$ and a new set $Y_{N+1}=\{\alpha\}$ of tape
letters. We also add two rules $\eta_0$, $\eta_1$ of the same form
$[k_{N+1}\to k_{N+1}, k_{N+2}\to \alpha k_{N+2}, \bar k_j\tool \bar
k_j, j=1,...,N]$. Note that now the new number of parts of $Q$ is
$N+2$, so we have to count modulo $N+2$ and instead of the
assumption that $\theta_{N+1}=\theta_1$ we have to assume that
$\theta_{N+3}=\theta_1$.

Notice that the new machine admits computations of the form

\begin{equation}\label{tolst}
\bar k_1...\bar k_Nk_{N+1}k_{N+2}\to_{\eta_0}...\to_{\eta_0}\bar
k_1...\bar k_Nk_{N+1}\alpha^mk_{N+2}
\to_{\eta_1\iv}...\to_{\eta_1\iv}\bar k_1...\bar
k_Nk_{N+1}k_{N+2}.\end{equation}

We keep notation $\sss$ for the new  $S$-machine and we keep
notation $N$ for the number of parts of the set of $Q$-letters
(instead of $N+2$).

Our goal is to cross-breed  $\sss$ with another $S$-machine in
order to slow it down whilst preserving the width of computations.

\subsection{The adding $S$-machine}

Let $A$ be a finite set of letters. Let the set $A_1$ be a copy of
$A$. It will be convenient to denote $A$ by $A_0$. For every
letter $a_0\in A_0$ let $a_1$ denote its copy in $A_1$. Consider
the following auxiliary ``adding" $S$-machine $Z(A)$.

Its set of state letters is $P_1\cup P_2\cup P_3$ where $P_1=\{L\},
P_2=\{p(1),p(2),p(3)\}, P_3=\{R\}$. The set of tape letters is
$Y_1\cup Y_2$ where $Y_1=A_0\cup A_1$ and $Y_2=A_0$.

The machine $Z(A)$ has the following positive rules (there $a$ is an
arbitrary letter from $A$). The comments explain the meanings of
these rules.

\begin{itemize}

\item $r_1(a)=[L\to L, p(1)\to a_1\iv p(1)a_0, R\to R]$.
\me

{\em Comment.} The state letter $p(1)$ moves left searching for a
letter from $A_0$ and replacing letters from $A_1$ by their copies
in $A_0$.

\me

\item $r_{12}(a)=[L\to L, p(1)\to a_0\iv a_1p(2), R\to R]$.

\me

{\em Comment.} When the first letter $a_0$ of $A_0$ is found, it is
replaced by $a_1$, and $p$ turns into $p(2)$.

\me

\item $r_2(a)=[L\to L, p(2)\to a_0p(2)a_0\iv, R\to R]$.

\me

{\em Comment.} The state letter $p(2)$ moves toward $R$.

\me

\item $r_{21}=[L\to L, p(2)\tool p(1), R\to R]$, $Y_1(r_{21})=Y_1,
Y_2(r_{21})=\emptyset$.

\me

{\em Comment.} $p(2)$ and $R$ meet, the cycle starts again.

\me

\item $r_{13}=[L\tool L, p(1)\to p(3), R\to R]$, $Y_1(r_{13})=\emptyset,
Y_2(r_{13})=A_0$.

\me

{\em Comment.} If $p(1)$ never finds a letter from $A_0$, the cycle
ends, $p(1)$ turns into $p(3)$; $p$ and $L$ must stay next to each
other in order for this rule to be executable.

\item $r_{3}(a)=[L\to L, p(3)\to a_0p(3)a_0\iv, R\to R]$,
$Y_1(r_3(a))=Y_2(r_3(a))=A_0$

 \me

{\em Comment.} The letter $r_3$ returns to $R$.

\end{itemize}

For every letter $a\in A$ we set $r_i(a\iv)=r_i(a)\iv$ ($i=1,2,3$).

The natural projection of every admissible word of $Z(A)$ onto $A$
takes $a_1$ and $a_0$ to $a$, and all other letters to $1$. The
following lemma is obvious.

\begin{lemma}\label{qqq2} Let $\theta$ be a rule in $Z(A)$,
$W=q_1u_1q_2...u_nq_{n+1}$ with base $q_1q_2...q_{n+1}$ be an
admissible word for $\theta$. Then $\theta\cdot W=q_1'\bar
w_1u_1w_2q_2'\bar w_2u_2...w_{n+1}q'_{n+1}$ where $w_j$ and $\bar
w_j$ are empty if $q_j\in \{L,R\}^{\pm 1}$ and if
$q_j=p\in\{p(1),p(2),p(3)\}$, then $w_j=w(p)$ and $\bar w_j=\bar
w(p)$ are determined by the following formulas:

\begin{equation}\label{rr78}
w(p)=\left\{\begin{array}{ll} a_1\iv & \hbox{ if }\tau_i=r_1(a), p=p(1),\\
                             a_0 & \hbox{ if }\tau_i=r_2(a), p=p(2),\\
                             a_0 & \hbox{ if }\tau_i=r_3(a), p=p(3),\\
                             (a_0\iv a_1)^{\epsilon}&\hbox{ if
                             }\tau_i=r_{12}(a)^{\epsilon}, \epsilon=\pm 1, p=p(1) \hbox{ if } \epsilon=1, p=p(2) \hbox{ if } \epsilon=-1,\\
                             \emptyset &\hbox{ if
                             }\tau_i\in \{r_{13}^{\pm 1}, r_{21}^{\pm 1}\},

\end{array}\right.
\end{equation}
\begin{equation}\label{rr79}
\bar w(p)=\left\{\begin{array}{ll} a_0 & \hbox{ if }\tau_i=r_1(a),p=p(1),\\
                             a_0\iv & \hbox{ if }\tau_i=r_2(a),p=p(2),\\
                             a_0\iv & \hbox{ if }\tau_i=r_3(a),p=p(3),\\
                             \emptyset & \hbox{ if
                             }\tau_i=r_{12}(a)^{\pm 1},\\
                             \emptyset &\hbox{ if }\tau_i\in \{r_{13}^{\pm1}, r_{21}^{\pm1}\}.
\end{array}\right.
\end{equation}
Finally if $q_j=p\in \{p(1)\iv, p(2)\iv, p(3)\iv\}$, then $w_j=\bar
w(p\iv)\iv, \bar w_j=w(p\iv)\iv$.
\end{lemma}

Lemma \ref{qqq2} immediately implies.

\begin{lemma} \label{lm569} Suppose that an admissible word $W$ has the form
$LupvR$ (resp. $p\iv upvR$) where $u,v$ are words in $(A_0\cup
A_1)^{\pm 1}$. Let $\theta\cdot W=Lu'p'v'R$ (resp. $\theta\cdot
W=(p')\iv u'p'v'R$). Then the projections of $uv$ and $u'v'$ (resp.
$v\iv uv$ and $(v')\iv u'v'$) onto $A$ are freely equal.
\end{lemma}

\begin{rk}\label{rk90} {\rm If we replace every letter in $A_i$ by its index $i$, then every
word $u$ in the alphabet $A_0\cup A_1$ turns into a binary number
$b(u)$. If the machine starts with the word $Lup(1)R$ where $u$ is a
positive word in $A_0$, then $b(u)=0$ and each cycle of the machine
amounts to adding $1$ to $b(u)$. After $2^{|u|}$ cycles the machine
stops, the admissible word becomes $Lup(3)R$. Let us compute the
length of this computation.

Notice that if $u=(b_1)_0...(b_n)_0$ for some $b_i\in A$ and at
the beginning of a cycle of the computation the last $k$
$a$-letters in the admissible word are from $A_1$ and the
$a$-letter number $n-k$ is from $A_0$ ($k<n$), then this cycle of
the computation has the following history

\begin{equation}\label{his}
r_1(b_n)...r_1(b_{n-k+1})r_{12}(b_{n-k})r_2(b_{n-k+1})...r_2(b_n)r_{21}
\end{equation}
(the letter $p(1)$ moves left searching for the first letter from
$A_0$ and replacing every $a$-letter from $A_1$ by the
corresponding letter from $A_0$; then $(b_{n-k})_0$ is replaced by
$(b_{n-k})_1$ and $p(1)$ is replaced by $p(2)$; then $p(2)$ moves
right, and after it meets with $R$, it is replaced by $p(1)$
again). The length of the cycle is $2k+2$. The number of cycles of
this length is $2^{n-k-1}$. Therefore the total length of all of
these cycles is $\sum_{k=0}^{n-1} (2k+2)2^{n-k-1}$. The length of
the last cycle is $2n+1$. Hence the total length of the
computation is
$$\sum_{k=0}^{n-1}(2k+2)2^{n-k-1}+2n+1=2^n(\sum_{k=0}^{n-1}
\frac{k+1}{2^k})+2n+1< 6\cdot 2^n$$ for every $n$ since
$\sum_{k=0}^{\infty} \frac{k+1}{2^k}=4$ and $2\cdot 2^n\ge 2n+1$ for
every $n\ge 0$. Hence the length of the computation is between $2^n$
and $6\cdot 2^n$.}
\end{rk}

Informally speaking, the remaining part of the section is devoted to
describing all possible computations of the machine $Z(A)$. In the
next section, we shall use that information to describe computation
of a composition of the $S$-machine from \cite{SBR} and $Z(A)$. We
first consider the case when the base is $LpR$. We show that if the
machine $Z(A)$ works without changing the length of an admissible
word, then the computations are in some sense unique and are
subcomputations of the computations described in Remark \ref{rk90}.
Then we consider computations where the lengths of words can change.
We show (and this is a standard feature of $S$-machines used in
\cite{SBR} and other papers) that as soon as the length of an
admissible word increases during the computation, it cannot decrease
again later, that is there are no trapezia that look like a honey
pot (width in the middle is bigger than the width of the bottom and
the top). If the base is not normal, say, it is $p\iv pR$, we show
that, again, there are no very wide ``honey pots".

\begin{lemma}\label{lm88} Suppose $W$ is a word with $\base(W)=LpR$. Then
there are at most two rules of $Z(A)$ that are applicable to $W$
word without changing its length.
\end{lemma}

\proof Let $W=LupvR$. Assume that $\theta\cdot W=W'=Lu'p'v'R$ for
a rule $\theta$. If $\theta=r_1(a)$ for some $a\in A^{\pm 1}$ then
the only ways equality $|W|=|W'|$ can occur is when either $v$
starts with $a_0\iv$ or $u$ ends with $a_1$ (by Lemma \ref{qqq2}).
Thus there are at most two choices for $a$ in that case. A similar
statement holds for $r_2(a)$, $r_3(a)$. Moreover, since
applicability of $r_i(a)$ determines the value of $p=p(i)$, rules
of only one of these three types can apply to $W$. For
$\theta=r_{12}(a)$, we have $|W|=|W'|$ only if $u$ ends with $a_0$
but not with $a_1\iv a_0$ and $p=p(1)$. Hence if $r_{12}(a)$ does
not change the length then only one other rule does not change the
length (either $r_1(b)$ where $b\iv$ is the first letter of $v$ or
$r_{21}\iv$ if $v$ is empty). Similar arguments hold in all other
cases.
\endproof

\begin{lemma}\label{lm89} For every admissible word $W$ with $\base(W)=LpR$,
every rule $\theta$ applicable to $W$, and every natural number
$t>1$, there is at most one computation $W\to_\theta W_1\to ...\to
W_t$ of length $t$ where the lengths of the words are all the same.
\end{lemma}
\proof Let $W\to_{\theta} W_1\to_{\theta_1}
W_2\to...\to_{\theta_{t-1}}W_{t-1}$ be a computation where all words
have the same length. By Lemma \ref{lm88}, $W_2$ is completely
determined by $W_1$ (because the history of the computation is a
reduced word, and so $\theta_1\ne\theta\iv$), $W_3$ is completely
determined by $W_2$, etc.
\endproof

For every word $f$ and every $i\le |f|$ we denote the prefix of $f$
of length $i$ by $f[i]$.

\begin{lemma}\label{lm92} Let $W=LvpuR$, where $p\in\{p(1),p(2)\}$, $u$ is a
word in $A_0^{\pm 1}$, and $v$ is a word in $(A_0\cup A_1)^{\pm
1}$. Suppose that a non-empty computation
\begin{equation}\label{c} W=W_1\to...\to f\cdot W\end{equation} is
such that all words $W_s$ in the computation have the same length,
$f\cdot W=Lv'p'uR$, $p'\in\{p(1),p(2)\}$, and either

(1) $p=p(2)$ and $f[1]=r_2(a)$ or

(2) $p=p(1)$ and $f[1]=r_1(a)\iv$

where $a$ is the first letter of $u$. Suppose also that $f[j]\cdot
W$ does not have the form $Lv''p''uR$ for every $1<j<|f|$. Then
$u$ is a positive word, $v'=v$ and $p'=p(1)$ in case (1) and
$p'=p(2)$ in case (2). Moreover, under the above conditions, the
computation (\ref{c}) is unique.
\end{lemma}

\proof We shall consider only the case when $p=p(2)$. The other
cases are similar. Let $u=b_1b_2...b_n$, $b_i\in A_0^{\pm 1}$.
Suppose that $u$ is not positive, $u=u[k]b\iv b_{k+2}...b_n$ where
$b_{k+2},...,b_n\in A$ (it could happen that $b_{k+2}...b_n$ is
empty, i.e. $n=k+1$).

Let $f[1]=\theta$. Since all words in the computation have the same
length, we can conclude, by Lemma \ref{lm89}, that for every $j$
between $1$ and $|f|$ there exists exactly one computation of length
$j$ starting with $W\to_\theta \theta\cdot W$. The history of that
computation must be $f[j]$.

Using Remark \ref{rk90} and (\ref{his}) one can find a computation
$W\to_\theta \theta\cdot W\to...$ with history $g$ of the form
$$g=r_2(b_1)...r_2(b_n)r_{21}r_{12}(b_n)r_{21}r_1(b_n)r_{12}(b_{n-1})
r_2(b_{n-1})r_2(b_n) r_{21} ...r_{21}r_1(b_n)...r_1(b_{k+2}),$$ and
$g\cdot W=Lvu_1b\iv p(1)b_{k+2}...b_nR$. Therefore either $g$ is a
prefix of $f$ or $f$ is a prefix of $g$.

Note that $f$ cannot be a prefix of $g$ because $f\cdot W=Lv'p'uR$
and  there is no $j$ such that $g[j]\cdot W$ has that form. Hence
$g$ is a prefix of $f$, i.e. $g=f[j]$ for some $j$. Moreover
$j<|f|$.

Notice that since $b\in A_0$, there is no rule except for
$r_1(b_{k+2})\iv$ (if $n\ge k+2$) or $r_{21}\iv$ (if $n=k+1$) which
can be applied to that word without increasing its length. This
contradicts the assumption that $f$ is reduced.  Thus $u$ is
positive.

Now we can consider the computation with the history of the form

$$h=r_2(b_1)...r_2(b_n)r_{21}r_{12}(b_n)r_{21}r_1(b_n)r_{12}(b_{n-1})
r_2(b_{n-1})r_2(b_n) r_{21} ...r_{21}r_1(b_n)...r_1(b_1)$$ (see
(\ref{his}) again) such that $h\cdot W=Lvp(1)uR$ and words
$h[i]\cdot W$, $i<|h|$, do not have the form $Lv''p''uR$. Since no
words $f[i]\cdot W$ have that form by assumption, we can conclude
that $f=h$. This completes the proof.
\endproof

Now we shall find out what happens when the lengths of the words
change during a computation.

\begin{lemma} \label{lm91} Let $W$,  be an admissible word
with $\base(W)=LpR$. Let
 $\theta\in \{r_1(a)^{\pm 1}$,
$r_2(a)^{\pm 1}$, $r_3(a)^{\pm 1}\}$. Suppose that $|\theta\cdot
W|>|W|$. Then for every computation $W\to_\theta W_1\to_{\theta'}
W_2$ we have $|W_2|>|W_1|$. Moreover, if $\theta'$ is not of the
form $r_{12}(a)^{\pm 1}$, then the $p$-letters in $W_2$ and $W_1$
are the same and $W_2$ does not have a 2-letter subword of the
form $pR$.
\end{lemma}
\proof Suppose that $\theta=r_1(a)$ (the other cases are similar).
If $|\tau\cdot W|>|W|$, then $\tau\cdot W=Lua_1\iv p(1)a_1vR$ where
the right hand side is a reduced word. It is easy to see that the
only rule that can apply to $\tau\cdot W$ without increasing the
length is $\tau\iv$, and the only type of rules that can apply are
$r_1(b)$ and $r_{12}(b)$. This immediately implies both statements
of the lemma.
\endproof

\begin{lemma} \label{lm93} Let $W=LvpuR$ and $\base(W)=LpR$. Suppose that
$|\theta\cdot W|>|W|$. Then for every computation $W\to_\theta
W_1\to W_2\to...\to f\cdot W$, we have $|W_i|>|W|$ for every $i\ge
1$.
\end{lemma}

\proof By contradiction, suppose that there exists a computation
$W\to_\theta W_1\to...\to f\cdot W$ such that $|f\cdot W|\le |W|$.
Consider such a computation with the smallest $|f|=t$ and smallest
$|u|$ for all such computations of length $t$. Then
$|W|<|W_1|=...=|W_{t-1}|>|W_t|$. By Lemma \ref{lm91},
$\theta=r_{12}(a)^{\pm 1}$ for some $a$.

We have $\theta\cdot W=Lv(a_0\iv a_1)^{\epsilon} p(m)uR$ for some
$\epsilon\in\{-1,1\}$, $m=(\epsilon+3)/2$, and some words $u, v$
where all letters in $u$ belong to $A_0^{\pm 1}$. We shall assume
that $\epsilon=1$ (the other case is similar). So $m=2$.

Suppose that the letter $a_1$ inserted by $f[1]$ is not touched
during the computation $W_1\to ...\to W_{t-1}$ (i.e. it is not
cancelled with a letter inserted by one of the rules of this
computation). Since the letters of $u$ do not disappear during the
computation $W_1\to ...\to W_{t-1}$ (they may only change indices
from $0$ to $1$ by Lemma \ref{lm569}), the word $W_{t-1}$ has the
form $Lva_0\iv a_1u_1pu_2R$ where $u_2$ is a suffix of $u$ and
$u_1$ is obtained from the corresponding prefix of $u$ by changing
indices of the letters. Applying Lemma \ref{lm91} to the inverse
computation $W_t\to W_{t-1}\to...\to W$ we conclude that the last
rule $\theta'$ in $f$ is $r_{12}(b)^{\pm 1}$ for some $b$. Hence
$W_t=Lv'p'u_2R$. Since $|W_t|<|W_{t-1}|$ and since $f$ is a
minimal counterexample, we conclude that $u_2$ must be equal to
$u$. So $u_1$ must be empty (otherwise the computation $W_t\to
W_{t-1}\to ...\to W$ would be a smaller counterexample since it
has the same length as $W_1\to...\to W_t$, but starts with
$W_t=Lv'p'u_2R$ with $|u_2|<|u|$), and either $p=p(2)$, and the
last rule of $f$ is $r_{12}(a)\iv$ or $p=p(1)$, the first letter
of $u$ is $a_0\iv$, and the last rule of $f$ is $r_1(a)\iv$. By
Lemma \ref{lm92}, both cases are impossible.

Therefore we can assume that $a_1$ is touched by the computation.
Hence, by Lemma \ref{lm92} applied to the computation $W_1\to...\to
W_{t-1}$, $u$ is a positive word, and if $s$ is the number of the
rule that touches $a_1$, we have $s\le t-1$, $f[s]\cdot W=Lva_0\iv
p(1)a_0uR$ and the rule number $s$ in $f$ is $r_1(a)$. Since $u$ is
positive, the word $a_0u$ is reduced. But then every rule of $Z(A)$
except $r_1(a)\iv$ would increase the length of the word, a
contradiction. This completes the proof of the lemma.
\endproof

\begin{lemma}\label{lm90} Let $\base(W)=LpR$.
Then for every computation $W=W_0\to W_1\to \dots\to W_t=f\cdot W$
of the $S$-machine $Z(A)$:

\begin{enumerate}
\item $|W_i|\le \max(|W_0|, |f\cdot W|)$, $i=0,...,t$,

\item If $W=LupR$ where $p=p(1)$ (resp. $p=p(3)$),
$f\cdot W$ contains $p(3)R$ (resp. $p(1)R$) and all $a$-letters in
$W, f\cdot W$ are from $A_0^{\pm 1}$, then the length of $f$ is
between $2^{|u|}$ and $6\cdot 2^{|u|}$, $u$ is a positive word, and
all words in the computation have the same length.
\end{enumerate}
\end{lemma}

\proof 1. Immediately follows from from Lemma \ref{lm93}.

2. We consider the case when $W=Lup(1)R$, the other case is similar.
If $u$ is empty, the statement is obvious. So assume that $u$ is not
empty.  By Lemmas \ref{lm93} and \ref{lm569}, all words in the
computation have the same length.

The letter $p(3)$ can occur only after rule $r_{13}$ is executed.
The admissible word $f[i]\cdot W$ to which $r_{13}$ is applied
must have the form $Lp(1)uR$. Indeed, the $a$-letters in
$f[i]\cdot W$ must be from $A$ (by the definition of $r_{13}$),
and the projection of $W$ and $f[i]\cdot W$ onto $A$ must coincide
by Lemma \ref{lm569}. The word $f[i-1]\cdot W$ in the computation
must be $La_1p(1)u'R$ where $a_0$ is the first letter of $u$,
$a_0u'=u$. Hence during the computation the first letter of $u$
must change the index from $0$ to $1$. Hence for some $i_2< i_1$,
we have $f[i_2]\cdot W=La_0p(2)u'R$ and $u[1]$ is a positive
letter. Applying now Lemma \ref{lm92} to the smallest initial part
of the subcomputation $f[i_2]\cdot W,...,f[i-1]\cdot W$ satisfying
the conditions of this lemma, we conclude that $u$ is a positive
word.

The first letter in $f$ is either $r_{12}(a)$ or $r_{21}\iv$.
Suppose first that $f[1]=r_{21}$.

By Remark \ref{rk90} there exists a computation connecting
$Lup(1)R$ and $Lup(3)R$. It has the history
$g=r_{12}(b_1)r_{21}r_1(b_1)r_{12}(b_2)r_2(b_1)...r_{13}...r_3(b_2)r_3(b_1)$
where $u=b_s...b_2b_1$, $b_i\in A_0$. Hence, by Lemma \ref{lm89},
either $g$ is a prefix of $f$ or $f$ is a prefix of $g$. But note
that there is only one rule applicable to the word $Lup(3)R$ that
does not increase its length (namely $r_3(b_1)\iv$). Hence $f=g$
and so $|f|$ is between $2^{|u|}$ and $6\cdot 2^{|u|}$.

Now let $f[1]=r_{21}\iv$. Then, as in the previous paragraph, we
deduce that $f$ must start with $f[s+1]=r_{21}\iv
r_2(b_1\iv)r_2(b_2\iv)...r_2(b_s\iv)$. But then $f[s+1]\cdot
W=Lp(2)uR$ and there is no rule except $r_2(b_s)$ that can be
applied to this word without increasing the length. The $s+2$-nd
rule in $f$ cannot be $r_2(b_s)$ since $f$ is reduced. Thus this
case is impossible.\endproof

\begin{lemma}\label{lmpp} Let $W=W_0=(p^{(0)})\iv u_0p^{(0)}$ be an admissible word
with $\base(W)=p\iv p$ ($p^{(0)}\in \{p(1),p(2),p(3)\}$. Let
$W\to_{\theta_1} W_1\to_{\theta_2}...\to_{\theta_t} W_t$ be a
computation. Let $W_j=(p^{(j)})\iv u_jp^{(j)}$, $u_j=w_j\iv u_{j-1}
w_j$ where $w_j$ is defined in Lemma \ref{qqq2}, $j=1,2,...$.
Suppose that none of $p^{(j)}$ is equal to $p(3)$. Then:
\begin{enumerate}
\item The word $\theta_1\theta_2...$ is a subword of the word of the
form $r_1...r_1 x_1 r_2...r_2 y_1 r_1...r_1 x_2...$ where $r_i$
stands for any $r_i(a)$, $x_j\in \{r_{12}(a), r_{21}\iv\mid a\in
A\}$, $y_j\in \{r_{12}(a)\iv, r_{21}\mid a\in A\}$.

\item  If none of the
rules $\theta_1,\theta_2,...$ is $r_{21}^{\pm 1}$, then none of the
words $w_1,w_2,...,w_t$ is empty, and the product $w_1...w_t$ is a
freely reduced word.

\item The word $\theta_1...\theta_t$ is completely determined by the
rule $\theta_1$  and the word $w_1w_2...w_t$.
\end{enumerate}
\end{lemma}

\proof 1. The first statement is obvious.

2. Statement 1 and Lemma \ref{qqq2} imply that the sequence $w_1,
w_2,...,w_t$ is a subsequence  of the sequence
\begin{equation}\label{eq3}\begin{array}{l} a_1(1),...,a_1(t_1),(a_0(t_1+1)\iv
a_1(t_1+1))^{\epsilon_1},a_0(t_1+2),...,a_0(t_2),\\(a_0(t_2+1)\iv
a_1(t_2+1))^{\delta_1},a_1(t_2+2),... \end{array}\end{equation}
where $a_i(j)\in A_i^{\pm 1}$, $\epsilon_j\in \{0,1\}$, $\delta_j\in
\{0,-1\}$, and we set $v^0=1$ for every word $v$. Moreover if
$r_{21}^{\pm1}$ does not appear in the computation, none of the
$\delta_j$ and $\epsilon_j$ are equal to 0. Hence  $w_1w_2...w_t$ is
freely reduced.

3. The third statement follows from the form of the sequence
(\ref{eq3}). Indeed, if $\theta_1=r_1(a)$, then $w_1=a_1\in A_1^{\pm
1}$, and $\theta_2$ is completely determined by the next one or two
letters of the word $w_1w_2...w_t$: if the second letter is $b_1$,
then $\theta_2=r_1(b)$; if it is $b_0\iv$ and the third letter is
$b_1$, then $\theta_2=r_{12}(b)$; if it is $b_0$ and the third
letter either does not exist or has index 0, then $\theta_2=r_{21}$.
Similarly for other choices of $\theta_1$,  the second and the third
letter of $w_1w_2...w_2$ completely determine $\theta_2$. Now we can
complete the proof by induction on $t$.
\endproof

\begin{lemma} \label{lm901} Suppose $\base(W)\in\{LpR, p\iv pR\}$, both $W$ and
$f\cdot W$ contain $p(1)R$ (resp. $p(3)R$) and all $a$-letters in
$W$, $f\cdot W$ are from $A_0$. Then $f$ is empty.
\end{lemma}

\proof The statement can be proved in the same way as part 2 of
Lemma \ref{lm90} provided $\base(W)=LpR$.

Let $\base(W)=p\iv pR$, $W=p(1)\iv up(1)R=f\cdot W$ where all
letters in $u$ are from $A_0^{\pm 1}$ (the case of $p(3)$ is
similar). By Lemma \ref{lmpp}, part 1, $u=w\iv uw$ where
$w=w_1w_2...$ where $w_j$'s are determined by the formulas from
Lemma \ref{qqq2}, the product $w_1w_2...$ considered as a word in
the alphabet $(A_0\cup A_1)^{\pm 1}\cup\{1\}$ (i.e. we do not throw
away the empty factors) does not contain subwords $xx\iv$, $x1x\iv$
or $11$. By Lemma \ref{qqq2}, $...\bar w_2\bar w_1=1$. Then, by
Lemma \ref{qqq2}, since all letters of $u$ are from $A_0$, all rules
in $f$ are from $\{r_2(a), r_3(a), r_{21}^{\pm 1}, r_{13}^{\pm
1}\}$. But in that case $\bar w_j=w_j\iv$ by Lemma \ref{qqq2}.
Therefore the product $...\bar w_2\bar w_1$ is reduced. Since that
product is 1, $w_1, w_2,..., \bar w_1, \bar w_2,...$ are empty
words. Since the factorization $w_1w_2...$ does not contain $11$, we
conclude that either $f$ is empty or $f=r_{21}^{\pm 1}$ (by Lemma
\ref{qqq2}). The second option is clearly impossible.
\endproof

We shall need the following general statement from \cite{OScol}.

\begin{lemma}[\cite{OScol}, Lemma 8.1] \label{powers}
For arbitrary elements $u,v,w$ of $F$ and any integer $t\ge 0$, the
length of an arbitrary product $u^{j}wv^{j}$ in $F$ is not greater
than $2(|u|+|v|+|w|)+|u^{t}wv^{t}|$ provided $0\le j\le t$.
\end{lemma}

\begin{lemma}\label{qqq} Suppose that one of the following conditions for an
admissible word $W$ of $Z(A)$ is satisfied (there $p=
\{p(1),p(2),p(3)\}$):

\begin{enumerate}
\item $W$ does not contain a $p$-letter.

\item $\base(W)=Lpp\iv$;

\item $\base(W)=pp\iv p$;

\item $\base(W)=p\iv pR$;

\item $\base(W)=LpR$.
\end{enumerate}
Then the width of any computation $$W=W_0\to_{\theta_0}
W_1\to_{\theta_1}...\to_{\theta_{t-1}} W_t$$ is at most
$C\max(|W|,|W_t|)$ for some constant $C$.
\end{lemma}

\proof 1 is obvious: the length of the admissible word does not
change during the computation.

\me

2. Let $\base(W)=Lpp\iv$. Let $W=W_0\to_{\theta_0}
W_1\to_{\theta_1}...\to_{\theta_{t-1}} W_t$ be a computation. Then
$W_i=Lu_iq_iv_iq_i\iv$, where $q_i\in \{p(1),p(2),p(3)\}$. It is
easy to see that for each $i=0,...,t-1$, $u_{i+1}=u_iw_i$,
$v_{i+1}=\bar w_i\iv v_i\bar w_i$ (equalities are in the free group)
where $w_i$ and $\bar w_i$ are determined by Lemma \ref{qqq2} (see
(\ref{rr78}), (\ref{rr79})).

Therefore $\bar w_i$ is a letter from $A_0^{\pm 1}$ or $1$, and the
projection of $w_i$ onto $A_0$ is freely equal to $\bar w_i$. In
particular, $|\bar w_i|\le |w_i|$.

Since there is no $R$ between $p$ and $p\iv$, none of $\tau_i$ is
$r_{21}$.

If $\theta_i=r_{13}$, then $u_i=\emptyset$, and
$w_0...w_{i-1}=u_0\iv$. Therefore $|v_i|\le
2|u_0|+|v_0|<3\max(|u_0|,|v_0|)$.

Hence it is enough to assume that rules $r_{13}^{\pm 1}$ do not
occur during the computation and prove that in that case, say,
$|W_i|<7\max(|W_0|, |W_t|)$.

We can assume that all rules $\theta_i$ are in $\{r_1(a), r_2(a),
r_{12}(a)^{\pm 1}\}$ (the case when all $\theta_i$ are of the form
$r_3$ is similar but easier). Then the history of the computation is
a subword of a word of the following form (we write $r_i$, $r_{12}$
instead of $r_i(a)$ and $r_{12}(a)$):

\begin{equation}\label{rrr}
r_1...r_1r_{12}r_2...r_2r_{12}^{\iv}r_1...r_1...\end{equation} Since
the history is reduced, by Lemma \ref{lmpp},  $w_0w_1...w_i$, is a
freely reduced word for every $i$. In addition $|\bar w_0\bar
w_1...\bar w_i|\le |w_0w_1...w_i|$. Hence $3i+|u_0|+|v_0|+3\ge
|W_i|>i-|u_0|$ for every $i$. Therefore $|W_i|\le
3i+|W_0|\le3t+|W_0|<|W_0|+3|W_t|+3|u_0|<4|W_0|+3|W_t|\le
7\max(|W_0|,|W_t|)$.

3. Let $\base(W)=pp\iv p$. Let $W=W_0\to_{\theta_0}
W_1\to_{\theta_1}...\to_{\theta_{t-1}} W_t$ be a computation,
$W_i=q_iu_iq_i\iv v_iq_i$, where $q_i\in \{p(1),p(2),p(3)\}$.

Let $W_j$ for some $j$ be a longest word in the computation. Then we
can assume without loss of generality that $|W_0|$ is shorter than
any of the words $W_1,...,W_j$, and $W_t$ is shorter than any of the
words among $W_j,...,W_t$. So assume that $|W_t|<|W_j|>|W_0|$.

Since the base of $W$ does not contain $L$ between $p\iv$ and $p$ or
$R$ between $p$ and $p\iv$, rules $r_{21}^{\pm 1}, r_{13}^{\pm 1}$
are excluded, so $\tau_i\in \{r_1(a), r_2(a), r_3(a), r_{12}(a)^{\pm
1}\}$. As in part 2, we shall assume that $r_3$ does not appear in
the computation.

By Lemma \ref{qqq2}, $u_{i+1}=\bar w_iu_i\bar w_i\iv$,
$v_{i+1}=w_i\iv v_i w_i$ where $w_i=w(q_i)$ and $\bar w_i=\bar
w(q_i)$ are determined by formulas (\ref{rr78}), (\ref{rr79}).

As in part 2, the history $h=\theta_0...\theta_{t-1}$ of the
computation is a subword of the word of the form (\ref{rrr}), the
product $w_0w_1w_2...w_{t-1}$ is a freely reduced word, and $|\bar
w_j|\le |w_j|$ for every $j$.

Therefore if $|v_{i+1}|>|v_i|$ for some $i$, then $|v_i|<...<|v_t|$
and $W_t$ cannot be shorter than $W_{t-1}$: indeed
$|u_{t-1}|-|u_{t}|\le 2$ by (\ref{rr79}) from Lemma \ref{qqq2} and
$|v_t|-|v_{t-1}|\ge 2$ since $v_t$ is a conjugate of $v_{t-1}$.
Similarly if $|v_{i+1}|<|v_i|$, then $W_0$ cannot be shorter than
$W_1$. Hence $|v_0|=|v_1|=...=|v_t|$, i.e. all $v_i$ are cyclic
shifts of $v$. This implies that the word $w_0...w_{t-1}$ is
periodic with period $d\le |v|$. By Lemma \ref{lmpp}, the history
$h$ of the computation is determined by its first letter $\theta_1$
and the word $w=w_1w_2...$. Every letter in $w$ is contained in one
of the $w_i$, so it corresponds to one of $\theta_i$. Consider the
letters in $w$ number $1, d+1, 2d+1,...$. These letters are the same
since the word $w$ is periodic with period $d$. Let $D$ be the total
number of rules in $Z(A)$. Then among the first $D+1$ rules
corresponding to these letters, there are two equal rules. Since
none of the words $w_i$ contain two same letters, we can deduce that
the word $h$ is periodic with period $d_1\le (D+1)d$.

Therefore, by Lemma \ref{qqq2}, the sequence $\bar w_0,...,\bar
w_{t-1}$ is periodic with the same period $d_1$. Let $z=\bar
w_{d_1-1}...\bar w_0$. Then $u_j=z_jz^su_0z^{-s}z_j\iv$,
$u_t=z_tz^{s'}u_0z^{-s'}z_t\iv$ for some words $z_j, z_t$ of length
at most $d_1$ and some $s, s'$. Now we can apply Lemma \ref{powers}
and deduce that $|u_j|\le 2(2d_1+|u_0|)+|u_t|+4d_1\le
C_1(|v_0|+|u_0|+|u_t|)$ for some constant $C_1$. Therefore $|W_j|\le
C(|W_0|+|W_t|)$ for some constant $C$ as required.

4. We can assume that $W$ is a shortest word in the computation.
Using notation similar to Case 3, $W_i=p\iv v_ipu_iR$ and
$v_{i+1}=w_i\iv v_iw_i$, $u_{i+1}=\bar w_i u_i$. As in Case 3, we
can assume that $|v_0|=|v_1|=...$ which implies, as before, that the
history of the computation is periodic with period $d_1\le
(D+1)|v|$. The proof can be completed as in Case 3, using Lemma
\ref{powers}.

Case 5 follows from Lemma \ref{lm90}.
\endproof

\subsection{A composition of $\sss$ and the adding machine}

\label{compsa}

Now let us define a ``composition" $\sss\circ Z$ of $\sss$ and
$Z(A)$. Essentially we insert a $p$-letter between any two
consecutive $k$-letters in admissible words of $\sss$, and treat any
subword $k_i...p...k_{i+1}$ as an admissible word for $Z(A)$.

First, for every $i=0,...,N-1$, we make two copies of the alphabet
$Y_i$ of $\sss$ ($i=1,...,N-1$): $Y_{i,0}=Y_i$ and $Y_{i,1}$. The
set of state letters of the new machine is $$K_1\cup P_1\cup
K_2\cup P_2\cup...\cup P_{N-1}\cup K_N$$ where $P_i=\{p_i,
p_i(\theta,1), p_i(\theta,2), p_i(\theta,3)\mid
\theta\in\Theta\}$, $i=1,...,N-1$. We shall denote the components
of this union by $Q_1,...,Q_{2N-1}$.

The set of state letters is $$\bar Y=(Y_{1,0}\cup Y_{1,1})\cup
Y_{1,0}\cup (Y_{2,0}\cup Y_{2,1})\cup Y_{2,0}\cup...\cup
(Y_{N-1,0}\cup Y_{N-1,1})\cup Y_{N-1,0};$$ the components of this
union will be denoted by $\bar Y_1,...,\bar Y_{2N-2}$.

The set of positive rules $\bar\Theta$ of $\sss\circ Z$ is a union
of the set of suitably modified positive rules of $\sss$ and
$2(N-1)|\Theta|$ copies $Z_i(\theta)^+$ ($\theta\in\Theta,
i=1,...,N$) of positive rules of the machine $Z(Y_i)$ (also suitably
modified).

More precisely, every positive rule $\theta\in \Theta^+$ of the form
$$[k_1u_1\to k_1'u_1', v_1k_2u_2\to
v_1'k_2'u_2',...,v_{N-1}k_N\to v_{N-1}'k_{N}']$$ where $k_i, k_i'\in
K_i$, $u_i$ and $v_i$ are words in $Y$, is replaced by

$$\bar\theta=\begin{array}{l}[k_1u_1\to k_1'u_1', v_1p_1\tool
v_1'p_1(\theta,1), k_2u_2\to k_2'u_2', ...,\\
 v_{N-1}p_{N-1}\tool v_{N-1}'p_{N-1}(\theta,1), k_{N}\to
 k_{N}']\end{array}$$ with
$\bar Y_{2i-1}(\bar\theta)=Y_{i,0}(\theta)$ and $Y_{2i}=\emptyset$
for every $i$.

Thus each modified rule from $\Theta$ turns on $N-1$ copies of the
machine $Z(A)$ (for different $A$'s).

Each machine $Z_i(\theta)$ is a copy of the machine $Z(Y_{i})$ where
every rule $\tau=[U_1\to V_1, U_2\to V_2, U_3\to V_3]$ is replaced
by the rule of the form
$$\bar\tau_i(\theta)=\left[
\begin{array}{l}\bar U_1\to \bar V_1, \bar U_2\to \bar V_2, \bar U_3\to \bar V_3,\\
k_j'\to k_j', p_j(\theta,3)\tool  p_j(\theta, 3),
j=1,...,i-1,\\
p_s(\theta,1)\tool p_s(\theta,1), k_{s+1}'\to k_{s+1}',
s=i+1,...,N-1\end{array}\right]$$ where $\bar U_1, \bar U_2, \bar
U_3, \bar V_1, \bar V_2, \bar V_3$ are obtained from $U_1, U_2, U_3,
V_1, V_2, V_3$, respectively, by replacing $p(j)$ with
$p_i(\theta,j)$, $L$ with $k_{i}'$ and $R$ with $k_{i+1}'$, and for
$s\ne i$, $\bar Y_{2s-1}(\bar\tau_i(\theta))=Y_{i,0}$.


In addition, we need the following {\em transition} rule
$\zeta(\theta)$ that returns all $p$-letters to their original form.

$$[k_i'\to k_i', p_j(\theta,3)\tool
p_j, i=1,...,N, j=1,...,N-1].$$

\me

Thus while the machine $Z_i(\theta)$ works all other machines
$Z_j(.)$, $j\ne i$, must stay idle (their state letters do not
change and do not move away from the corresponding $k$-letters).
After the machine $Z_i(\theta)$ finishes (i.e. the state letter
$p_{i}(\theta,3)$ appears next to $k_{i+1}$), the next machine
$Z_{i+1}(\theta)$ starts working. After all $p$-letters have the
form $p_j(\theta,3)$ we can apply $\zeta(\theta)$ and turn all
$p_j(\tau,3)$ into $p_j$. Thus in order to simulate a computation of
$\sss$ consisting of a sequence of applications of rules $\theta_1,
\theta_2,...,\theta_s$, we first apply all rules corresponding to
$\theta_1$, then all rules corresponding to $\theta_2$, then all
rules corresponding to $\theta_3$, etc.

The modified rules $\bar\theta$ of $\sss$ will be called {\em basic
rules}. We shall call basic rules, transition rules $\zeta(\theta)$,
and their inverses {\em main rules}.

\subsection{Properties of the $S$-machine $\sss\circ Z$}

{\em Notation.} For every word $W$, we shall denote the number of
$a$-letters in $W$ by $|W|_a$.

Every word in the alphabet $\{Q_1,...,Q_{2N-1}\}$ is called a {\em
base word}. Let $\bbb$ be a finite set of base words.

\begin{df}\label{narrow}{\rm
We call a base word {\em $\bbb$-covered} if
\begin{itemize}

\item it is covered by bases from $\bbb$ (i.e. every letter
belongs to a subword from $\bbb$),

\item it starts and ends with the same $q$-letter $x$.
\end{itemize}
We call a base word $w$ {\em $\bbb$-tight} if it has the form $uxvx$
where $xvx$ is a $\bbb$-covered word,  $w$ does not contain any
other $\bbb$-covered subwords. A base word is called {\em
$\bbb$-narrow} if it does not contain $\bbb$-covered subwords.}
\end{df}

\begin{lemma}\label{width}
There exists a finite set of base words $\bbb$ such that
\begin{itemize}
\item[(*)] the length of every $\bbb$-narrow base is smaller than a constant $K_0$,

\item[(**)] for every admissible word $W$ with base from $\bbb$
the width of every computation $W\to W_1\to...\to W_t$ does not
exceed $C(|W|+|W_t|+\log_2 t)$ for some constant $C$.
\end{itemize}
\end{lemma}

\proof Let $\bbb$ be the set of bases of the form
$(k_iq_ik_{i+1})^{\pm 1}$, $(q_i\iv q_ik_{i+1})^{\pm 1}$,
$(k_iq_iq_i\iv)^{\pm 1}$, $k_ik_i\iv$, $k_i\iv k_i$, $(q_iq_i\iv
q_i)^{\pm 1}$ ($q_i=\{p_i(\theta,1),p_i(\theta,2),p_i(\theta,3),
p_i\}$). Lemma \ref{lmad} implies that condition (*) holds for
$\bbb$. Indeed, base words of the form $(q_iq_i\iv q_i)^{\pm 1}$ and
all base words that start and end on the same $k$-letters and do not
contain $(q_iq_i\iv q_i)^{\pm 1}$ as a subword are $\bbb$-covered,
and every word of length at least $6N$ contains one of these
$\bbb$-covered words.

Let $W$ be an admissible word with one of the bases from $\bbb$. By
passing to $W\iv$ if necessary, we can assume that $\base(W)$ has
one of the forms $k_iq_ik_{i+1}$, $q_i\iv q_ik_{i+1}$,
$k_iq_iq_i\iv$, $k_ik_i\iv$, $k_i\iv k_i$, $q_iq_i\iv q_i$.

Let $f=f_1f_2f_3$ where $f_1,f_3$ contain no basic and transition
rules or their inverses, $f_2$ starts and ends with a basic or
transition rule or its inverse (one or more of these subwords may be
empty).

{\bf Case 1.} Suppose that $f=f_1$, so $f$ does not contain main
rules.

Notice that we can represent $f$ as a product $z_1z_2...z_s$ where
each $z_i$ is a history of computation of one of the machines
$Z_{j(i)}(\theta)$ for some $\theta$ and $j(i)$ depending on $i$,
$j(i)\ne j(i+1)$. Since $f^{\pm 1}$ does not contain basic and
transition rules, all $\theta$'s are the same. By Lemma \ref{lm90},
\ref{qqq}, if $\base(W)\not\in \{k_ipk_{i+1},p\iv pk_{i+1}\}$, then
$s=1$,  the width of the computation does not exceed $C\max(W,
f\cdot W)$ for some constant $C$ and no basic or transition rule can
apply to $f\cdot W$.

If $W=k_iupvk_{i+1}$ or $W=p\iv upvk_{i+1}$, then each
subcomputation corresponding to $z_j$ either does not affect the
admissible word, or it is essentially a computation of $Z(Y_i)$
(with $L$ replaced by $k_i$, $R$ replaced by $k_{i+1}$, etc.). In
the first case the width is equal to the length of $W$. In the
second case, the width of the subcomputation corresponding to $z_j$
does not exceed $C$ times the sum of lengths of the beginning and
ending words of the computation (by Lemma \ref{qqq}). Moreover in
the second case the end word of the subcomputation, except possibly
for $f\cdot W$, has the form $k_iwp'k_{i+1}$ (resp. $(p')\iv
wp'k_{i+1}$) for some word $w$ in $Y_{i,0}^{\pm 1}$.  By Lemma
\ref{lm569}, the words $w$ and $uv$ (or $w$ and $v\iv uv$) must have
the same projections onto $Y_i$. Therefore all $w$'s have the same
lengths which does not exceed $\max(|W|,|f\cdot W|)$ (resp.
$2\max(|W|,|f\cdot W|)$). Therefore the width of the whole
computation does not exceed $C(|W|+|f\cdot W|)$. Moreover we can
conclude that if $f\cdot W$ contains a 2-letter subword $pR$, then
$|f\cdot W|\le |W|$.

{\bf Case 2.} Suppose that $f_2$ is of length $1$, i.e. it contains
exactly one main rule. Then the base of $W$ is $k_iq_ik_{i+1}$ or
$q_i\iv q_ik_{i+1}$.  We can repeat the argument from the last
paragraph of Case 1, applied to the subcomputations corresponding to
$f_1$ and $f_3\iv$. Note that a main rule $f_2^{\pm 1}$ must be
applicable to $f_1\cdot W$ and to $f_3\iv \cdot (f\cdot W)$), hence
these words must contain subwords of the form $pk_{i+1}$ and all
$a$-letters must belong to $Y_{i,0}$. Hence the lengths of the end
words of these subcomputations cannot be bigger than $2|W|$ or
$2|f\cdot W|$ by Lemma \ref{lm569}.

{\bf Case 3.} Suppose that $f_2$ contains at least two main rules.

Again $\base(W)$ must belong to $\{k_iq_ik_{i+1}, q_i\iv
q_ik_{i+1}\}$ for some $i$. By the argument in Case 1, $|f_1\cdot
W|\le 2|W|$, $|f_1f_2\cdot W|\le 2|f\cdot W|$, so we can assume that
$f_1$ and $f_3$ are empty, thus $f$ starts and ends with a main
rule. We can assume that this computation has minimal length among
all computations connecting $W$ and $f\cdot W$ and having the same
width.

Let $\tau g\tau'$ be a subword of $f$ where $\tau, \tau'$ are main
rules and $g$ does not contain main rules. Then either $g$ is empty
or $g=g_1...g_m$ and each $g_j$ is a non-empty product of rules of
one of the machines $Z_{s_{j}}(\theta_{j})$, $s_{j}\ne s_{j+1}$.

Suppose that $g$ is not empty and $s_{j}=i$ for some $j$. In this
case we shall call the subword $\tau g\tau'$ {\em active}. Let
$f=f'g_jf''$ for some $f', f''$. Let $W'=f'\cdot W$, $W''=g_j\cdot
W'$. Then both $W'$ and $W''$ must contain subwords $pk_{i+1}$ where
$p\in \{p_i(.,1), p_i(.,3)\}$ and all $a$-letters of $W', W''$ are
in $Y_{i,0}$. Since $g_j$ is not empty, by Lemma \ref{lm901}, either
$W'$ contains $p_i(.,1)$ and $W''$ contains $p_i(.,3)$, or $W'$
contains $p_i(.,3)$ and $W''$ contains $p_i(.,1)$. Since rule
$r_{13}^{\pm1}$ cannot be applicable to admissible words of $Z(A)$
of the form $p\iv pR$, we conclude that $\base(W)=k_iq_ik_{i+1}$. By
Lemma \ref{lm90}, the length of $g_{j}$ is at least $2^{|W'|-3}$
(here $|W'|-3$ is the number of $a$-letters in $W'$), and the
lengths of the words in the subcomputation $W'\to ...\to W''$ are
the same.

If none of $s_j$ are equal to $i$, then we call the subword $\tau
g\tau'$ {\em passive}. In that case all rules of $g$ fix admissible
words with base $\base(W)$. Therefore all words in the
subcomputation corresponding to $g$ are the same. From the
definition of main rules, then either $\tau'=\tau\iv$, or $\tau$ is
the inverse of a basic rule and $\tau'$ is another basic rule, or
$\tau$ is a transition rule and $\tau'$ is the inverse of another
transition rule. If $\tau'\ne\tau\iv$, $g$ is non-empty, and the
subcomputation corresponding to $g$ does not contain the longest
word in the computation, we can use Lemma \ref{lmcut} and remove the
subcomputation corresponding to $g$ and obtain a shorter computation
of the same width connecting $W$ and $f\cdot W$. That would
contradict the assumption that $f$ is the shortest possible. Hence
we can assume that if $\tau g\tau'$ is passive and $g$ is not empty,
and the subcomputation corresponding to $g$ does not contain the
longest word in the computation,  then $\tau'=\tau\iv$.

Notice also that we have proved, that in our computation,
$|W_j|=|W_{j+1}|$ unless the $j+1$-st rule of the computation is a
main rule.

Let $J=|f[t]\cdot W|$ be the width of the computation $W\to...\to
f\cdot W$. We can assume (by the remark made in the previous
paragraph) that the $t$-th rule in $f$ is a main rule. For every
admissible word $U$ with $\base(U)=\base(W)$ and every basic rule
$\tau$ the difference $|U|-|\tau\cdot U|$ is bounded from above by a
constant $c$. We can assume that $J-|W|>4c$. Let $t_0$ be such that
$J-|f[t_0]\cdot W|\ge 4c$ and $t-t_0$ is minimal possible. Then the
suffix of $f[t]$ obtained by removing $f[t_0]$ has the form $\tilde
f=\tau_1g_1\tau_2g_2...\tau_jg_j\tau_{j+1}$ such that
\begin{itemize}
\item $j\ge 3$,
\item all $\tau_l$ are main rules, all $g_l$ do not contain main rules,

\item all admissible words in the subcomputation corresponding to
$\tilde f$ are of length between $J-4c$ and $J$,
\item one of the admissible words in that subcomputation has length
at most $J-4c$.
\end{itemize}

If one of the subwords $\tau_lg_l\tau_{l+1}$ in $\tilde f$ is
active then the length of $f$ is at least $2^{J-4c-3}$. So $J\le
\log_2 |f|+4c+3$ and (**) holds.

Thus we can assume that all subwords $\tau_lg_l\tau_{l+1}$ are
passive. But then, as we proved before, for every $l$ either $g_l$
is empty or $\tau_{l+1}=\tau_l\iv$. Moreover if $g_l$ is empty
then $\tau_{l+1}$ cannot be the inverse of a basic rule or a
transition rule. Hence $g_{l+1}$ cannot be empty. Therefore the
lengths of admissible words in the subcomputation corresponding to
$\tilde f$ can differ by at most $2c$, a contradiction with the
fact that one of these words has length $J$ and the length of
another one is at most $J-4c$.
\endproof

From now on we shall fix the set $\bbb$ of base words from Lemma
\ref{width}, and a constant $K$ such that
\begin{equation}\label{kk}K>2K_0.\end{equation}

\section{Properties of the group $\sss\circ Z$}
\label{propgr}

Let $\Theta_+$ be the set of positive rules of $\sss\circ Z$.

We redefine the map $\sigma$ and the stop word $\tilde W$ (see
Section \ref{secsm}) as follows:
$$\sigma'(v)=k_1(1)vp(1)k_2(1)p(2)k_3(1)...p(N-1)k_N(1)$$ for every
word $v$ over $X$,
$$\tilde W'=k_1(0)p(1)k_2(0)p(2)k_3(0)...p(N-1)k_N(0).$$ Recall that $\LL$
is the language recognized by $\sss$.

\begin{lemma}\label{lmlll1} Let $W_1=\sigma'(v)$ and $W_1\to_{\theta_1}
...\to_{\theta_{t-1}} W_t$ be any computation of $\sss\circ Z$. Let
$\theta_{i_1},...,\theta_{i_s}$ be all basic rules or their inverses
in that computation. Let $\bar\theta_{i_j}$ be the rule of $\sss$
corresponding to $\theta_{i_j}$. For each $j$ let $\bar W_j$ be the
(natural) projection of $W_j$ on the alphabet of $q$- and
$a$-letters of $\sss$. Then  there exists a computation $\bar
W_1\to...\to \bar W_t$ of the machine $\sss$, whose history is the
reduced form of the word $\bar\theta_1...\bar\theta_s$.
\end{lemma}

\proof This immediately follows from Lemmas \ref{lm569}, \ref{lmcut}
and the definition of $\sss\circ Z$.\endproof

\begin{lemma} \label{lmlll}
A word $v$ over $X$ belongs to $\LL$ if and only if $\sigma'(v)$ and
$\tilde W'$ are conjugate in $\sss\circ Z$.
\end{lemma}

\proof Suppose that $\sigma'(v)$ is a conjugate of $\tilde W'$.
Then, by Lemma \ref{obratnaja}, there exists a computation
$\sigma'(v)\to ...\to \tilde W'$ of $\sss\circ Z$. Then, by Lemma
\ref{lmlll1}, there exists a computation
$\overline{\sigma'(v)}\to...\to \bar W'$. Note that
$\overline{\sigma'(v)}=\sigma(v)$, $\bar W'=\tilde W$. Hence $v$ is
accepted by $\sss$. Hence $v\in \LL$.

Conversely suppose that $v\in \LL$. Then, by Lemma \ref{lm1}, there
exists a computation of $\sss$ connecting $\sigma(v)$ with $\tilde
W$ and consisting of positive words. Then, by Lemma \ref{lm92}
(applied several times), there exists a computation of $\sss\circ Z$
connecting $\sigma'(v)$ and $\tilde W'$, so by Remark \ref{rk99},
$\sigma'(v)$ and $\tilde W'$ are conjugate in $\sss\circ Z$.
\endproof

Let $L$ be the maximal length of a defining relation of $\sss\circ
Z$.

\begin{lemma}\label{band} Let $\ttt$ be a $\theta$-band with base
of length $l_b$.  Let $l_a$ be the number of $a$-edges in the top
path $\topp(\ttt)$. Then the length of $\ttt$ is between
$l_a-(L-1)l_b$ and $l_a+(L+1)l_b$.
\end{lemma}

\proof Every $a$-letter in the label of $\topp(\ttt)$ labels an edge
of a $(a,\theta)$-cell or on a $(q,\theta)$-cell. Every
$(q,\theta)$-cell in $\ttt$ has at most $L$ $a$-edges lying on
$\topp(\ttt)$,  and every $(a,\theta)$-cell has at most one $a$-edge
lying on $\topp(\ttt)$.  At most $L$ of these $a$-edges can belong
to the same $(q,\theta)$-cell. Hence $|\ttt|\ge l_a-Ll_b+l_b$.

The number of $(a,\theta)$-cells having no common $a$-edges with
$\topp(\ttt)$, does not exceed $Ll_b$ because at least one of the
two $a$-edges of this cell is glued to an $a$-edge of a
$(q,\theta)$-cell.  This proves the inequality $|\ttt|\le
l_a+Ll_b+l_b$.
\endproof

\begin{lemma}\label{area} Let $\Delta$ be a trapezium of height $h\ge 1$ whose base is
$\bbb$-tight and $\bbb$-covered. Then the area of $\Delta$ does not
exceed $Ch(|W|_a+|W'|_a+ \log h+1)$, where $W, W'$ are the labels of
its top and bottom, respectively for some constant $C$.
\end{lemma}

\proof Notice that the number of letters in the base of $\Delta$
does not exceed $K_0+1$. Since $\base(W)$ is covered by bases from
$\bbb$, $\Delta$ is covered by subtrapezia whose width, by Lemmas
\ref{width} and \ref{band}, does not exceed a constant times
$|W|_a+|W'|_a+\log h+1$, height does not exceed $h$,  and the
number of these subtrapezia does not exceed a constant (since
$\bbb$ is finite and from every cover of $\base(W)$ by words from
$\bbb$, we can find a subcover which covers every letter of
$\base(W)$ at most a constant number of times). Hence the width of
$\Delta$ does not exceed $C(|W|_a+|W'|_a+\log h+1)$ for some
constant $C$. Since, by Lemma \ref{band}, the area of $\Delta$
does not exceed a constant times the hight times the width, the
statement of the lemma follows.
\endproof

\subsection{A modified length function on $\sss\circ Z$}
\label{length}

Let us modify the length function on the group $\sss\circ Z$. The
standard length of a word (path) will be called its {\em
combinatorial length}. From now on we use the word {\em length} for
the modified length. As before let $L$ be the maximal length of a
defining relation of $\sss\circ Z$. We set the length of every
$q$-letter equal 1, and the length of every $a$-letter equal a small
enough number $\delta$ so that

\begin{equation}
2-4L\delta-LK\delta>\delta. \label{param}
\end{equation}

We also set to 1 the length of every word of length $\le L$ which
contains exactly one $\theta$-letter and no $q$-letters (such
words are called $(\theta,a)$-{\em syllables}). The length of a
decomposition of an arbitrary word in a product of letters and
$(\theta,a)$-syllables is the sum of the lengths of the factors.
{\em The length of a word $w$} is the smallest length of such
decompositions. {\em The length of a path} in a diagram is the
length of its label. The {\em perimeter} $|\partial\Delta|$ of a
van Kampen diagram is similarly defined by cyclic decompositions
of the boundary $\partial\Delta$. It follows from this definition
that for any product $s=s_1s_2$ of two paths in a \vk diagram, we
have

\begin{equation}
|s_1|+|s_2|\ge |s|> |s_1|+|s_2|-L\delta \label{delta}
\end{equation}

A maximal $\theta$-band of a van Kampen diagram $\Delta$ is called a
{\em rim band} if its top or its bottom side lies on the contour
$\partial\Delta$.
\bigskip

\begin{lemma} \label{rim} Let $\Delta$ be a van Kampen diagram
whose rim band $\ttt$ has base with at most $K$ letters. Denote by
$\Delta'$ the subdiagram $\Delta\backslash \ttt$. Then
$|\partial\Delta|-|\partial\Delta'|\ge\delta$.
\end{lemma}

\proof Let $s$ be the top side of $\ttt$ and
$s\subset\partial\Delta$. Note that by our assumptions the
difference between the number of $a$-edges in the bottom $s'$ of
$\ttt$ and the number of $a$-edges for $s$ cannot be greater than
$LK$. Hence $|s'|-|s|\le LK\delta$. However, $\Delta'$ is obtained
by cutting off $\ttt$ along $s'$, and its boundary contains two
$\theta$-edges fewer than $\Delta$. Hence we have $|s_0|-|s'_0|\ge
2-2L\delta$ by (\ref{delta}), for the complements $s_0$ and $s'_0$
of $s$ and $s'$, respectively, in the boundaries $\partial\Delta$
and $\partial\Delta'$. Finally,
$$|\partial\Delta|-|\partial\Delta'|\ge
2-2L\delta-LK\delta-2L\delta>\delta$$ by (\ref{param}) and
(\ref{delta}).
\endproof
The definition of length has also the following obvious consequence.

\begin{lemma} \label{ochev}
Let $s$ be a path in a diagram $\Delta$ having $c$ $\theta$-edges
and $d$ $a$-edges. Then

(a) $ |s|\ge \max(c, c+(d-Lc)\delta)$;

(b) $|s|=c$ if $s$ is a top or a bottom of a $q$-band.
\end{lemma}

\subsection{Combs}

\begin{df} {\rm We say that a reduced diagram $\Gamma$ is a {\em comb} if it has a
maximal $q$-band $\qq$ (the {\em handle} of the comb), such that

\begin{enumerate}
\item[$(C_1)$]$\bott(\qq)$ is a part of $\partial \Gamma$, and every maximal
$\theta$-band of $\Gamma$ ends at a cell in $\qq$.
\end{enumerate}

If in addition the following properties hold:
\begin{enumerate}
\item[$(C_2)$] one of the maximal $\theta$-bands $\ttt$ in
$\Gamma$ has a $\bbb$-tight base and

\item[$(C_3)$] other maximal
$\theta$-bands in $\Gamma$ have $\bbb$-tight or $\bbb$-narrow bases
\end{enumerate}
then the comb is called {\em tight}.

The number of cells in the handle $\qq$ is the {\em length} of the
comb, and the maximal length of the bases of the $\theta$-bands of a
comb is called the {\em basic width} of the comb. }
\end{df}

\begin{figure}[ht]
\centering
\unitlength 1mm 
\linethickness{0.4pt}
\ifx\plotpoint\undefined\newsavebox{\plotpoint}\fi 
\begin{picture}(42.75,44.5)(0,0)
\put(38.25,44){\rule{.88\unitlength}{.5\unitlength}}
\put(37.38,44.25){\line(-1,0){16.875}}
\put(20.5,44.25){\line(0,-1){2.875}}
\put(20.5,41.38){\line(1,0){16.875}}
\put(37.38,41.38){\line(-1,0){23.5}}
\put(13.88,41.38){\line(0,-1){2.875}}
\put(13.88,38.5){\line(1,0){23.625}}
\qbezier(37.63,38.63)(10.94,38.63)(9.5,33.63)
\put(9.5,33.63){\line(0,1){.25}}
\put(9.5,33.88){\line(0,-1){3}}
\qbezier(9.5,30.88)(10.13,36)(37.25,36.13)
\qbezier(37.25,35.88)(10.81,35.94)(12.13,30.75)
\put(12.13,31.13){\line(0,-1){2.625}}
\qbezier(12.13,28.63)(14.06,33.31)(37.25,32.75)
\put(37.25,9.5){\rule{3.5\unitlength}{34.88\unitlength}}
\put(37.63,9.5){\line(-1,0){22.125}}
\put(15.5,9.5){\line(0,1){3.375}}
\put(15.5,12.88){\line(1,0){22}}
\put(37.5,12.88){\line(-1,0){19.625}}
\multiput(17.88,12.88)(-.06435644,.03341584){101}{\line(-1,0){.06435644}}
\put(11.38,16.25){\line(0,1){3.5}}
\multiput(11.38,19.75)(.0576087,-.03369565){115}{\line(1,0){.0576087}}
\put(18,15.88){\line(1,0){19.375}}
\put(14.63,17.75){\line(0,1){3.25}}
\multiput(14.63,21)(.0583333,-.0333333){60}{\line(1,0){.0583333}}
\put(18.13,19){\line(1,0){19.625}}
\put(23.38,25.13){\makebox(0,0)[cc]{$\cdots$}}
\put(42.75,28.13){\makebox(0,0)[cc]{$\qq$}}
\put(29.5,24.63){\makebox(0,0)[cc]{$\Gamma$}}
\put(9.75,10.63){\makebox(0,0)[cc]{$\ttt_1$}}
\put(7,17.88){\makebox(0,0)[cc]{$\ttt_2$}}
\put(15.88,43.13){\makebox(0,0)[cc]{$\ttt_l$}}
\end{picture}
\caption{Comb.} \label{fig1}
\end{figure}

Notice that every trapezium is a comb.

\begin{df}\label{topbotom}
{\rm Let $\qq$ be a maximal $q$-band of $\Delta$. Then $\bott(\qq)$
divides $\Delta$ into two parts. The part containing $\qq$ is called
the {\em top part of $\Delta$ with respect to $\qq$}. The other part
is called {\em the bottom part of $\Delta$ with respect to $\qq$}.}
\end{df}

\begin{lemma} \label{est'takaya} Let $\Delta$ be a reduced diagram
with non-zero area. Assume that every rim band of $\Delta$ has base
of length at least $K$. Then there exists a maximal $q$-band $\qq$
in $\Delta$ such that the top part $\Gamma$ of $\Delta$ with respect
to $\qq$ is a tight comb.
\end{lemma}

\proof Let $\ttt_0$ be a rim band of $\Delta$. Its base $w$ is of
length at least $K$, and therefore $w$ has disjoint prefix and
suffix of lengths $K_0$ since $K>2K_0$ by (\ref{kk}). The prefix of
this base word must have its own $\bbb$-tight subprefix $w_1$, by
Lemma \ref{width}, part 1, and the definition of $\bbb$-tight words.
A $q$-edge of $\ttt_0$ corresponding to the last $q$-letter of the
$w_1$ is the start edge of a maximal $q$-band $\qq'$ which bounds a
subdiagram $\Gamma'$ containing a band $\ttt$ (a subband of
$\ttt_0$) satisfying property ($C_2$). It is useful to note that a
minimal suffix $w_2$ of $w$, such that $w_2\iv$ is tight, allows us
to construct another band $\qq''$ and a subdiagram $\Gamma''$ which
satisfies ($C_2$) and has no cells in common with $\Gamma'$.

\begin{figure}[ht]
\centering
\unitlength 1mm 
\linethickness{0.4pt}
\ifx\plotpoint\undefined\newsavebox{\plotpoint}\fi 
\begin{picture}(102.44,48.77)(20,0)
\put(37.63,43.75){\framebox(52,2.75)[cc]{}}
\multiput(30.5,32.88)(-.03368263,-.08757485){167}{\line(0,-1){.08757485}}
\multiput(24.88,18.25)(-.081522,.032609){23}{\line(-1,0){.081522}}
\multiput(23,19)(.03357143,.08357143){175}{\line(0,1){.08357143}}
\multiput(28.88,33.63)(.067308,-.033654){26}{\line(1,0){.067308}}
\qbezier(37.63,43.63)(32.69,40.31)(30.5,32.75)
\qbezier(24.75,18.13)(22.63,13.19)(24.5,12)
\qbezier(24.5,12)(28.31,9.81)(26.88,8.38)
\qbezier(26.88,8.38)(25.38,6.94)(28.88,6.75)
\qbezier(28.88,6.75)(30.44,6.94)(93.25,7.88)
\qbezier(93.25,7.88)(102.44,7.94)(93.88,29.25)
\qbezier(93.88,29.25)(88.31,42.75)(89.5,43.75)
\put(83.25,24){\makebox(0,0)[cc]{$\Delta$}}
\put(62.75,7.38){\rule{2.25\unitlength}{39.13\unitlength}}
\multiput(34.75,40.88)(.033689591,-.062267658){538}{\line(0,-1){.062267658}}
\multiput(33.13,38.63)(.03371063,-.061761811){508}{\line(0,-1){.061761811}}
\put(35.88,36){\circle*{2.37}} \put(37,35.25){\circle*{1.8}}
\put(34.5,38.5){\circle*{1.52}} \put(34.5,39.38){\circle*{1.58}}
\put(34.63,37){\circle*{1.35}} \put(35.5,37.5){\circle*{2}}
\put(34.63,38.13){\circle*{2.12}} \put(37.42,32.74){\circle*{1.71}}
\put(36.84,33.85){\circle*{1.79}} \put(36.58,33.06){\circle*{.53}}
\put(36.53,32.74){\circle*{.24}} \put(36.05,34.69){\circle*{1.27}}
\put(35.9,33.9){\circle*{.32}} \put(35.48,34.84){\circle*{.42}}
\put(34.79,36.26){\circle*{.53}} \put(33.53,38.73){\circle*{.66}}
\put(33.53,38.42){\circle*{.33}} \put(33.74,39.1){\circle*{.54}}
\put(35.37,38.79){\circle*{.77}} \put(35.48,38.21){\circle*{1.12}}
\put(35.05,39.15){\circle*{.9}} \put(37.47,33.48){\circle*{1.69}}
\put(37.74,34.27){\circle*{1.1}} \put(38.05,33.32){\circle*{1.28}}
\put(38.37,32.53){\circle*{1.32}} \put(38.42,32.95){\circle*{.76}}
\put(38.37,31.27){\circle*{2.3}} \put(39.31,29.48){\circle*{2.13}}
\put(40.57,27.7){\circle*{2.33}} \put(41.73,25.02){\circle*{1.89}}
\put(41.2,25.96){\circle*{2.04}} \put(40.52,26.49){\circle*{1.27}}
\put(39.42,28.22){\circle*{1.18}} \put(40.47,29.17){\circle*{.9}}
\put(39.63,30.48){\circle*{1.16}} \put(39.15,31.43){\circle*{1.14}}
\put(38.84,32.06){\circle*{1.05}} \put(37.31,31.85){\circle*{.84}}
\put(38.26,30.17){\circle*{.9}} \put(40.05,29.75){\circle*{1.37}}
\put(41.73,26.59){\circle*{1.27}} \put(41.57,27.28){\circle*{.9}}
\put(42.04,25.75){\circle*{1.49}} \put(42.57,24.49){\circle*{1.81}}
\put(42.52,25.23){\circle*{1.27}} \put(41.89,23.23){\circle*{.71}}
\put(41.83,23.7){\circle*{.97}} \put(43.57,21.92){\circle*{2.27}}
\put(44.83,19.24){\circle*{1.81}} \put(46.04,17.4){\circle*{2.33}}
\put(47.14,15.24){\circle*{2.35}} \put(48.61,13.24){\circle*{1.94}}
\put(47.67,14.14){\circle*{2.16}} \put(49.14,11.77){\circle*{2.32}}
\put(50.14,9.72){\circle*{2.22}} \put(51.19,8.36){\circle*{2.1}}
\put(52.19,7.67){\circle*{.85}} \put(50.45,7.57){\circle*{.64}}
\put(50.03,8.36){\circle*{.74}} \put(51.19,9.46){\circle*{1}}
\put(50.51,10.83){\circle*{.87}} \put(50.14,11.35){\circle*{1.04}}
\put(48.98,10.25){\circle*{.66}} \put(47.83,12.25){\circle*{.57}}
\put(47.62,12.77){\circle*{.54}} \put(48.56,14.4){\circle*{.87}}
\put(48.19,15.08){\circle*{.9}} \put(47.09,16.61){\circle*{1.18}}
\put(46.04,16.03){\circle*{1}} \put(45.99,15.66){\circle*{.53}}
\put(45.62,16.4){\circle*{.63}} \put(45.99,18.82){\circle*{1.18}}
\put(44.83,17.97){\circle*{.82}} \put(43.94,19.66){\circle*{.87}}
\put(44.78,20.23){\circle*{1.92}} \put(45.62,19.81){\circle*{.85}}
\put(45.93,19.45){\circle*{.63}} \put(43.2,23.39){\circle*{1.79}}
\put(43.25,24.23){\circle*{.74}} \put(43.88,23.02){\circle*{.9}}
\put(44.51,21.39){\circle*{1.37}} \put(44.46,22.34){\circle*{.54}}
\put(44.25,22.97){\circle*{.21}} \put(42.73,22.49){\circle*{1.4}}
\put(41.99,22.81){\circle*{.53}} \put(43.83,20.65){\circle*{1.59}}
\put(34.63,40.06){\circle*{.91}}
\qbezier(26.94,28.69)(39.16,23.84)(42.5,7.13)
\qbezier(26.13,26.75)(36.66,22.34)(40.31,7.06)
\qbezier(24.81,23.44)(35.78,16.5)(35.88,7.06)
\qbezier(24.06,21.5)(33.41,15.31)(33.88,6.88)
\put(27.5,27.13){\circle*{1.88}} \put(27.13,27.94){\circle*{1.03}}
\put(28.06,27.75){\circle*{.76}} \put(27.71,27.97){\circle*{.62}}
\put(26.56,27){\circle*{.54}} \put(29.04,26.3){\circle*{1.74}}
\put(28.51,27.05){\circle*{1.51}} \put(29.17,26.91){\circle*{1.15}}
\put(29.57,26.25){\circle*{2.13}} \put(28.11,26.16){\circle*{.45}}
\put(31.38,25.01){\circle*{2.15}} \put(30.63,25.94){\circle*{1.35}}
\put(31.16,26.03){\circle*{.73}} \put(30.54,26.47){\circle*{.62}}
\put(30.23,25.01){\circle*{.9}} \put(29.74,25.15){\circle*{.8}}
\put(30.58,24.57){\circle*{.62}} \put(32.92,23.56){\circle*{2.05}}
\put(32.53,24.62){\circle*{1.24}} \put(33.23,24.53){\circle*{.62}}
\put(33.59,24.17){\circle*{.62}} \put(33.81,23.78){\circle*{.71}}
\put(34.07,23.38){\circle*{.99}} \put(31.91,23.82){\circle*{1.16}}
\put(31.38,23.91){\circle*{.64}} \put(34.78,21.66){\circle*{2.27}}
\put(36.02,19.89){\circle*{2.04}} \put(37.48,18.12){\circle*{2.12}}
\put(38.36,15.95){\circle*{2.05}} \put(39.51,13.48){\circle*{2.21}}
\put(40.57,10.65){\circle*{2.03}} \put(41.01,8.35){\circle*{1.98}}
\put(41.14,9.5){\circle*{1.59}} \put(40.44,9.5){\circle*{1.42}}
\put(40.48,11.49){\circle*{1.79}} \put(39.6,12.37){\circle*{1.51}}
\put(39.64,12.07){\circle*{.97}} \put(39.6,11.49){\circle*{.9}}
\put(39.55,11.09){\circle*{.71}} \put(40.35,12.37){\circle*{1.43}}
\put(40.44,13.3){\circle*{.45}} \put(40.53,13.08){\circle*{.64}}
\put(39.2,14.67){\circle*{1.95}} \put(39.24,15.69){\circle*{1.16}}
\put(39.07,16.62){\circle*{.71}} \put(38.8,17.1){\circle*{.76}}
\put(38.45,17.59){\circle*{.89}} \put(37.96,18.56){\circle*{.99}}
\put(37.12,19.05){\circle*{1.61}} \put(36.99,20.02){\circle*{.73}}
\put(36.64,20.55){\circle*{.95}} \put(36.15,21.08){\circle*{.95}}
\put(35.71,21.7){\circle*{.81}} \put(34.29,22.54){\circle*{1.51}}
\put(34.74,22.85){\circle*{.95}} \put(33.68,22.01){\circle*{.96}}
\put(32.97,22.36){\circle*{.53}} \put(33.45,22.54){\circle*{.71}}
\put(33.19,22.14){\circle*{.36}} \put(35.09,20.42){\circle*{1.29}}
\put(36.33,18.83){\circle*{1.43}} \put(36.77,17.59){\circle*{.76}}
\put(36.37,18.08){\circle*{.62}} \put(37.57,16.88){\circle*{1.63}}
\put(38.23,14.85){\circle*{1.03}} \put(38.49,14.36){\circle*{1.01}}
\put(38.4,13.88){\circle*{.53}} \put(41.94,7.82){\circle*{.92}}
\put(41.85,8.57){\circle*{.81}} \put(41.76,9.1){\circle*{.62}}
\put(41.5,10.21){\circle*{.53}} \put(40.57,7.47){\circle*{.63}}
\put(41.01,7.42){\circle*{.73}} \put(41.5,7.47){\circle*{.62}}
\put(41.98,7.34){\circle*{.73}} \put(25.68,21.48){\circle*{1.82}}
\put(25.01,21.96){\circle*{1.42}} \put(24.97,22.76){\circle*{.88}}
\put(25.59,22.45){\circle*{.8}} \put(26.12,22.14){\circle*{.71}}
\put(24.48,21.52){\circle*{.71}} \put(24.84,21.17){\circle*{.4}}
\put(27.09,20.77){\circle*{1.75}} \put(28.51,18.92){\circle*{1.33}}
\put(28.11,19.84){\circle*{1.9}} \put(29.79,17.81){\circle*{1.63}}
\put(31.25,16.22){\circle*{1.99}} \put(32.57,14.81){\circle*{1.98}}
\put(33.41,12.51){\circle*{1.68}} \put(34.29,10.74){\circle*{2.13}}
\put(34.83,8.44){\circle*{1.86}} \put(34.65,9.59){\circle*{1.99}}
\put(34.83,7.95){\circle*{1.72}} \put(34.07,7.16){\circle*{.45}}
\put(35.53,7.29){\circle*{.62}} \put(34.87,7.2){\circle*{.57}}
\put(34.29,7.16){\circle*{.54}} \put(35.44,9.15){\circle*{.44}}
\put(33.99,8.84){\circle*{.64}} \put(33.94,8.22){\circle*{.36}}
\put(33.94,7.73){\circle*{.27}} \put(34.16,7.47){\circle*{.28}}
\put(33.68,9.99){\circle*{.64}} \put(33.23,11.62){\circle*{.62}}
\put(33.32,11.31){\circle*{.59}} \put(33.94,12.02){\circle*{1.91}}
\put(34.74,11.62){\circle*{.73}} \put(34.12,12.9){\circle*{1.06}}
\put(33.59,13.48){\circle*{1.71}} \put(32.57,13.26){\circle*{.8}}
\put(32.44,13.88){\circle*{1.06}} \put(31.69,14.67){\circle*{.62}}
\put(31.73,14.32){\circle*{.53}} \put(32.35,15.73){\circle*{1.35}}
\put(33.41,14.32){\circle*{1.18}} \put(31.64,15.29){\circle*{1.34}}
\put(30.85,17.19){\circle*{1.82}} \put(31.95,16.93){\circle*{.8}}
\put(32.26,16.35){\circle*{.76}} \put(31.69,17.15){\circle*{.64}}
\put(30.32,16.53){\circle*{.71}} \put(29.92,16.84){\circle*{.62}}
\put(30.41,18.3){\circle*{1.06}} \put(30.94,17.99){\circle*{.8}}
\put(29.79,18.65){\circle*{1.44}} \put(28.9,19.09){\circle*{1.88}}
\put(29.08,18.12){\circle*{1.03}} \put(28.59,18.21){\circle*{.53}}
\put(29.04,19.98){\circle*{.48}} \put(28.06,20.68){\circle*{.76}}
\put(27.71,21.08){\circle*{.45}} \put(26.83,21.52){\circle*{1.03}}
\put(26.34,22.01){\circle*{.54}} \put(26.16,20.64){\circle*{.87}}
\put(26.47,20.29){\circle*{.63}} \put(27.14,19.93){\circle*{.96}}
\put(27.18,19.4){\circle*{.28}} \put(27.44,19.22){\circle*{.69}}
\put(61.28,48.77){\makebox(0,0)[cc]{$\ttt_0$}}
\put(67.27,25.23){\makebox(0,0)[cc]{$\qq'$}}
\put(54.97,23.65){\makebox(0,0)[cc]{$\Gamma'$}}
\put(48.46,41.41){\makebox(0,0)[cc]{$\ttt$}}
\put(34.37,27.96){\makebox(0,0)[cc]{$\Gamma$}}
\put(23.86,28.06){\makebox(0,0)[cc]{$\ttt_1$}}
\put(52.24,4.41){\makebox(0,0)[cc]{$\qq$}}
\put(41.1,4.62){\makebox(0,0)[cc]{$\qq_2$}}
\put(34.06,4.41){\makebox(0,0)[cc]{$\qq_1$}}
\put(27.64,12.51){\makebox(0,0)[cc]{$\Gamma_1$}}
\put(46.27,18.21){\circle*{1.31}} \put(45.25,17.85){\circle*{1.41}}
\put(44.5,18.38){\circle*{.52}} \put(41.68,24.04){\circle*{1}}
\put(40.75,25.06){\circle*{.45}} \put(39.95,26.74){\circle*{.59}}
\put(40.75,28.86){\circle*{.75}} \put(48.79,10.69){\circle*{.73}}
\put(49.36,12.9){\circle*{.87}} \put(47.69,16.04){\circle*{.88}}
\put(47.16,17.1){\circle*{.71}} \put(50.25,8){\circle*{.71}}
\put(51.49,7.56){\circle*{.69}} \put(39.55,30.98){\circle*{.9}}
\put(37.21,19.71){\circle*{.81}} \put(38.94,12.46){\circle*{.59}}
\put(47.33,12.99){\circle*{.45}} \put(39.47,27.62){\circle*{.53}}
\put(39.6,27.31){\circle*{.53}} \put(44.28,22.63){\circle*{.73}}
\put(38.8,32.75){\circle*{.53}} \put(38.71,33.1){\circle*{.44}}
\put(50.91,10.34){\circle*{.63}} \put(51.13,9.99){\circle*{.45}}
\put(39.2,11.71){\circle*{.54}} \put(39.16,12.11){\circle*{.57}}
\put(35.13,22.67){\circle*{.44}} \put(39.16,16.22){\circle*{.76}}
\put(47.38,16.44){\circle*{.83}} \put(49.72,8.71){\circle*{.62}}
\put(52.15,8.04){\circle*{.56}}
\end{picture}

\caption{Lemma \ref{est'takaya}.} \label{fig2}
\end{figure}

Thus, there are $\qq$ and $\Gamma$ satisfying ($C_2$). Let us choose
such a pair with minimal $\area(\Gamma)$. Assume that there is a
$\theta$-band in $\Gamma$ which does not cross $\qq$. Then there
must exist a rim band $\ttt_1$ which does not cross $\qq$ in
$\Gamma$. Hence one can apply the construction from the previous
paragraph to $\ttt_1$ and construct two bands $\qq_1$ and $\qq_2$
and two disjoint subdiagrams $\Gamma_1$ and $\Gamma_2$ satisfying
the requirement ($C_2$) for $\Gamma$. Since $\Gamma_1$ and
$\Gamma_2$ are disjoint, one of them, say $\Gamma_1$, is inside
$\Gamma$. But the area of $\Gamma_1$ is smaller than the area of
$\Gamma$, and we come to a contradiction. Hence $\Gamma$ is a comb
and condition ($C_1$) is satisfied.

Assume that the base of a maximal $\theta$-band $\ttt$ of $\Gamma$
is not narrow. Then it has a $\bbb$-tight proper prefix (we may
assume that $\ttt$ terminates on $\qq$), and again one obtain a
$q$-band $\qq'$ in $\Gamma$, which provides us with a smaller
subdiagram $\Gamma'$ of $\Delta$, satisfying ($C_2$), a
contradiction. Hence $\Gamma$ satisfies property ($C_3$) as well.
\endproof

\begin{lemma} \label{comb} Let $l$ and $b$ be the length and
the basic width of a comb $\Gamma$ and let $\ttt_1,\dots \ttt_l$ be
consecutive $\theta$-bands of $\Gamma$ (as in Figure \ref{fig1}). We
can assume that $\bott(\ttt_1)$ and $\topp(\ttt_l)$ are contained in
$\partial\Gamma$. Denote by $\alpha=|\partial\Gamma|_a$ the number
of $a$-edges in the boundary of $\Gamma$, and by $\alpha_1$ the
number of $a$-edges on $\bott(\ttt_1)$. Then $\alpha + 4Clb\ge
2\alpha_1$, and the area of $\Gamma$ does not exceed $2Cbl^2+2\alpha
l$ for some constant $C$.
\end{lemma}

\proof For every $i$,  $\topp(\ttt_i)$ and $\bott(\ttt_{i+1})$ can
have an initial segment in common. Let $\alpha_i$ and $\alpha_i'$ be
the numbers of the $a$-edges in $\bott(\ttt_i)$ and $\topp(\ttt_i)$
respectively that belong to the boundary of $\Gamma$.

Let $n_i$ be the  length of $\ttt_i$. It follows from Lemma
\ref{band} that $|n_{i+1}-n_i|\le 2Cb+ \alpha'_i+\alpha_{i+1}$.
Since $n_1\le Cb+\alpha_1$ by the same lemma, we have for any $i$,
$$n_i\le (Cb+\alpha)
+2C(i-1)b+\alpha_1'+\alpha_2+\dots + \alpha_{i-1}'+\alpha_i\le 2Clb
+2\alpha$$ This inequality provides us with a required upper bound
for the area $\sum_{i=1}^l n_i$. Finally, $$\begin{array}{l}\alpha
-\alpha_1 \ge \sum_{i=1}^{l-1} \alpha_i' +\sum_{i=2}^l \alpha_i
+\alpha'_l\ge \sum_{i=1}^{l-1} |n_{i+1}-n_i| - 2Cb(l-1)+(n_l-Cb)\\
\ge \sum_{i=1}^{l-1} (n_i-n_{i+1}) - 2Cb(l-1)+(n_l-Cb) \ge
n_1-3Cbl\ge \alpha_1 - 4Cbl.\end{array}$$
\endproof

\section{Dispersion of bipartite chord diagrams}

\label{dispersion}

\begin{df}\label{grating} {\rm Recall that a chord diagram is a system of
chords in a disc such that the intersecting point of any two chords
is in the interior of the disc. We can consider the intersection
graph of a chord diagram whose vertices are chords and two chords
are connected when they intersect. If that graph is bipartite, the
chord diagram is called {\em bipartite}. We shall always consider
bipartite chord diagrams (BCD for short) with fixed subdivision into
two parts $\TT$ and $\QQ$, so chords from one of the parts do not
intersect. The intersection points of $\TT$- and $\QQ$-chords are
the {\em nodes of the BCD}.}
\end{df}

For an integer $K\ge 1$, we fix a $K$-tuple $\alpha$ of numbers
$0<\alpha_1\le \alpha_2\le...\le\alpha_K=1$. For every node $o$ on
a $\TT$-chord $\ttt$ and choice of the left-right direction on
$\TT$, we assign the weight $\alpha_l$ where $l$ is the minimum of
$K$ and the number of nodes on $\ttt$ to the left of $o$
(including $o$). Thus the weight of a node $o$ is smaller if it is
closer to the boundary (according to the chosen direction on the
$\TT$-chord containing $o$). Note that we assign two weights to
each node $o$ according to the two possible directions on a chord
containing $o$.

Let $\ccc$ be a $\QQ$-chord. Let $(o_1,o_2)$ be a pair of nodes
from $\ccc$. Fixing an orientation from $o_1$ to $o_2$ on $\ccc$,
we determine a left-to-right orientation on every $\TT$-chord
crossing $\ccc$. The node $o_i$ ($i=1,2$) lies on some $\TT$-chord
$\ttt_i$. Let $o_i'$ be the next node to the left of $o_i$ on
$\ttt_i$ or the intersection point of $\ttt_i$ and $\partial D$ if
there are no nodes on $\ttt_i$ between $o_i$ and $\partial D$ to
the left of $o_i$ (according to the chosen direction of $\ttt_i$).
If either both $o_i'$ lie on the same $\QQ$-chord or both lie on
$\partial D$, then we call the pair $(o_1,o_2)$ {\em good}.
Otherwise we call the pair $(o_1,o_2)$ {\em bad}. By definition,
the {\em weight} of the pair $(o_1,o_2)$ is the product of their
weights (corresponding to the direction from $o_1$ to $o_2$ on the
$\QQ$-chord). We shall set the weight of every good pair to 0.
Note that a pair $(o_1,o_2)$ may be good while the pair
$(o_2,o_1)$ is bad.

\begin{figure}[ht]
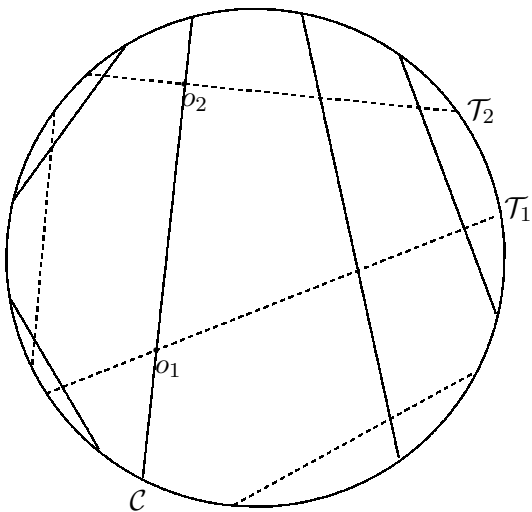

\centering
\unitlength .6 mm 
\linethickness{0.4pt}
\ifx\plotpoint\undefined\newsavebox{\plotpoint}\fi 

\caption{$(o_1,o_2)$ is a bad pair, $(o_2,o_1)$ is a good pair.}
\label{fig3}
\end{figure}

The $\alpha$-{\em dispersion} $\eee_\alpha (\cal C)$ of the system
of $\TT$-chords on a $\QQ$-chord $\cal C$ is the sum of weights of
all bad pairs of nodes $(o_1,o_2)$ such that $o_1$ and $o_2$  lie on
$\cal C$. For example, let ${\bf 1}$ be the $1$-tuple of numbers
$(1)$. Then $\eee_{\bf 1}(\cal C)$$=1$ for Figure \ref{fig3}.

\begin{df}\label{disorder}{\rm
The sum of dispersions $\eee_\alpha (\cal C)$ over all
$\QQ$-chords $\cal C$ is called the $\alpha$-dispersion
$\eee_\alpha(\gc)$ of the BCD $\gc$.}
\end{df}

Clearly $\eee_{\bf 1}(\gc)$ is the number of bad pairs of $\gc$,
and $\eee_\alpha(\gc)\le \eee_{\bf 1}(\gc)$ for every $\alpha$.

We need a quadratic upper bound for the $\alpha$-dispersion $\eee$
of a BCD $\gc$ in terms of the number of $\TT$-chords in $\gc$.

\begin{figure}[ht]
\centering
\unitlength .6mm 
\linethickness{0.4pt}
\ifx\plotpoint\undefined\newsavebox{\plotpoint}\fi 



\caption{Lemma \ref{entropy}. The chords $\ttt_1$, $\ttt_3$ are in
$X_k$, $\ttt_2$, $\ttt_4$ are in $X_{k+1}$. } \label{fig6}
\end{figure}

\begin{lemma} \label{entropy} Let $r$ be the number of $\TT$-chords of
a BCD $\gc$ on a disc $D$. Then the $\bf 1$-dispersion $\eee_{\bf
1}(\gc)$ of the BCD $\gc$ does not exceed $r^2-r$.
\end{lemma}

\proof We will induct on the number of $\TT$-chords in ${\cal G}$.
If the set $\TT$ is empty, then $\eee_{\bf 1}({\cal G})=0$, and the
statement is obviously true.

Thus, we may assume that there is a chord ${\cal T}$ in $\TT$. We
can assume that $\ttt$ is close to the boundary that is there are no
$\TT$-chords in one of the half-discs obtained by cutting $D$ along
$\ttt$. We denote by ${\cal C}_1,\dots, {\cal C}_l$ all the
$\QQ$-chords intersected by ${\cal T}$, where $l\ge 0$. We enumerate
and orient the $\QQ$-chords so that $\ccc_2$ is to the right of
$\ccc_1$, $\ccc_3$ is to the right of $\ccc_2$, etc., and there are
no nodes on $\ccc_k$ above $o_k=\ccc_k\cap \ttt$, for each
$k=1,...,l$.

Let $X_k$ be the set of $\TT$-chords $\ttt'$ of $\gc$ such that the
pair $(\ttt\cap \ccc_k, o_k)$ is bad. Denote by $L_k$ the number of
chords in $X_k$. Thus $L_k$ is the number of bad pairs of the form
$(o,o_k)$.

Note that a chord $\ttt'$ from $X_{k+1}$ cannot intersect $\ccc_k$
because otherwise the pair $(\ttt'\cap \ccc_k, o_k)$ would be good.
Hence the sets $X_k$ are pairwise disjoint, and $\sum L_k\le r-1$
because $\ttt$ does not belong to any $X_k$.

Similarly let $R_k$ be the number of bad pairs of the form
$(o_k,o)$. Then as above, the sum of all $R_k$ does not exceed
$r-1$. Hence

\begin{equation}\label{LR}
\sum L_k+\sum R_k\le 2(r-1)
\end{equation}

Consider the BCD $\gc_0$ obtained from $\gc$ by deleting the
$\TT$-chord $\ttt$. Notice that every bad pair $(o',o'')$ in $\gc_0$
is also bad in $\gc$ and the difference $\eee_{\bf 1}(\gc)-\eee_{\bf
1}(\gc_0)$ between the number of bad pairs for $\gc$ and the number
of bad pairs for $\gc_0$ is the number of bad pairs of the forms
$(o,o_k)$ and $(o_k,o)$ which does not exceed $2(r-1)$ by
(\ref{LR}). By the inductive assumption, $\eee_{\bf 1}(\gc_0)$ does
not exceed $(r-1)^2-(r-1)$. Hence $\eee_{\bf 1}(\gc)$ does not
exceed
$$(r-1)^2-(r-1)+2(r-1)= r^2-r.$$
\endproof

\begin{rk} {\rm The estimate $r^2-r$ in Lemma
\ref{entropy} is optimal: for every $r$ one can easily construct
(using the proof of Lemma \ref{entropy}) a BCD with $r$ $\TT$-chords
whose $\bf 1$-dispersion is exactly $r^2-r$.}
\end{rk}

\begin{lemma} \label{deletion} Let a BCD ${\cal G}_0$ be
obtained by a deleting (a) a $\TT$-chord or (b) a $\QQ$-chord $\ccc$
of a BCD ${\cal G}$. Then the $\alpha$-dispersion of $\gc_0$ does
not exceed the $\alpha$-dispersion of $\gc$.
\end{lemma}

\proof (a) If we remove a $\TT$-chord, then every bad pair in
$\gc_0$ is a bad pair in $\gc$ having the same weight. This implies
part (a).

\begin{figure}
\centering
\unitlength .6mm 
\linethickness{0.4pt}
\ifx\plotpoint\undefined\newsavebox{\plotpoint}\fi 
\begin{picture}(109.18,98.93)(30,25)
\put(109.18,66.5){\line(0,1){1.474}}
\put(109.15,67.97){\line(0,1){1.472}}
\multiput(109.07,69.45)(-.032,.3671){4}{\line(0,1){.3671}}
\multiput(108.95,70.91)(-.0298,.24384){6}{\line(0,1){.24384}}
\multiput(108.77,72.38)(-.03278,.20799){7}{\line(0,1){.20799}}
\multiput(108.54,73.83)(-.0311,.16079){9}{\line(0,1){.16079}}
\multiput(108.26,75.28)(-.033,.14365){10}{\line(0,1){.14365}}
\multiput(107.93,76.72)(-.03164,.11868){12}{\line(0,1){.11868}}
\multiput(107.55,78.14)(-.032997,.108469){13}{\line(0,1){.108469}}
\multiput(107.12,79.55)(-.031845,.092957){15}{\line(0,1){.092957}}
\multiput(106.64,80.95)(-.032864,.086057){16}{\line(0,1){.086057}}
\multiput(106.12,82.32)(-.033725,.079872){17}{\line(0,1){.079872}}
\multiput(105.54,83.68)(-.032639,.070373){19}{\line(0,1){.070373}}
\multiput(104.92,85.02)(-.033311,.065737){20}{\line(0,1){.065737}}
\multiput(104.26,86.33)(-.03234,.058673){22}{\line(0,1){.058673}}
\multiput(103.54,87.62)(-.032865,.055014){23}{\line(0,1){.055014}}
\multiput(102.79,88.89)(-.033308,.051596){24}{\line(0,1){.051596}}
\multiput(101.99,90.13)(-.033677,.048391){25}{\line(0,1){.048391}}
\multiput(101.15,91.34)(-.03272,.043697){27}{\line(0,1){.043697}}
\multiput(100.26,92.52)(-.032996,.041015){28}{\line(0,1){.041015}}
\multiput(99.34,93.66)(-.033214,.03847){29}{\line(0,1){.03847}}
\multiput(98.38,94.78)(-.033379,.03605){30}{\line(0,1){.03605}}
\multiput(97.38,95.86)(-.033495,.033744){31}{\line(0,1){.033744}}
\multiput(96.34,96.91)(-.035802,.033645){30}{\line(-1,0){.035802}}
\multiput(95.26,97.92)(-.038223,.033498){29}{\line(-1,0){.038223}}
\multiput(94.15,98.89)(-.040769,.033298){28}{\line(-1,0){.040769}}
\multiput(93.01,99.82)(-.043453,.033042){27}{\line(-1,0){.043453}}
\multiput(91.84,100.71)(-.046289,.032725){26}{\line(-1,0){.046289}}
\multiput(90.64,101.56)(-.051347,.033689){24}{\line(-1,0){.051347}}
\multiput(89.4,102.37)(-.054769,.033271){23}{\line(-1,0){.054769}}
\multiput(88.14,103.14)(-.058432,.032774){22}{\line(-1,0){.058432}}
\multiput(86.86,103.86)(-.06237,.032188){21}{\line(-1,0){.06237}}
\multiput(85.55,104.53)(-.070129,.03316){19}{\line(-1,0){.070129}}
\multiput(84.22,105.16)(-.075197,.03241){18}{\line(-1,0){.075197}}
\multiput(82.86,105.75)(-.085812,.0335){16}{\line(-1,0){.085812}}
\multiput(81.49,106.28)(-.092718,.032533){15}{\line(-1,0){.092718}}
\multiput(80.1,106.77)(-.100492,.031385){14}{\line(-1,0){.100492}}
\multiput(78.69,107.21)(-.11844,.03252){12}{\line(-1,0){.11844}}
\multiput(77.27,107.6)(-.13036,.03097){11}{\line(-1,0){.13036}}
\multiput(75.84,107.94)(-.16055,.03229){9}{\line(-1,0){.16055}}
\multiput(74.39,108.23)(-.18177,.03003){8}{\line(-1,0){.18177}}
\multiput(72.94,108.47)(-.24361,.0316){6}{\line(-1,0){.24361}}
\multiput(71.48,108.66)(-.29347,.02775){5}{\line(-1,0){.29347}}
\put(70.01,108.8){\line(-1,0){1.471}}
\put(68.54,108.89){\line(-1,0){1.473}}
\put(67.06,108.93){\line(-1,0){1.474}}
\put(65.59,108.91){\line(-1,0){1.472}}
\put(64.12,108.84){\line(-1,0){1.469}}
\multiput(62.65,108.73)(-.29286,-.03359){5}{\line(-1,0){.29286}}
\multiput(61.18,108.56)(-.20823,-.03124){7}{\line(-1,0){.20823}}
\multiput(59.73,108.34)(-.18114,-.03365){8}{\line(-1,0){.18114}}
\multiput(58.28,108.07)(-.14389,-.03194){10}{\line(-1,0){.14389}}
\multiput(56.84,107.75)(-.12972,-.03356){11}{\line(-1,0){.12972}}
\multiput(55.41,107.38)(-.108711,-.032192){13}{\line(-1,0){.108711}}
\multiput(54,106.96)(-.099847,-.033381){14}{\line(-1,0){.099847}}
\multiput(52.6,106.5)(-.086298,-.032225){16}{\line(-1,0){.086298}}
\multiput(51.22,105.98)(-.08012,-.033133){17}{\line(-1,0){.08012}}
\multiput(49.86,105.42)(-.070613,-.032117){19}{\line(-1,0){.070613}}
\multiput(48.52,104.81)(-.065982,-.032823){20}{\line(-1,0){.065982}}
\multiput(47.2,104.15)(-.061716,-.033424){21}{\line(-1,0){.061716}}
\multiput(45.9,103.45)(-.055256,-.032456){23}{\line(-1,0){.055256}}
\multiput(44.63,102.7)(-.051841,-.032925){24}{\line(-1,0){.051841}}
\multiput(43.39,101.91)(-.048639,-.033317){25}{\line(-1,0){.048639}}
\multiput(42.17,101.08)(-.045628,-.033641){26}{\line(-1,0){.045628}}
\multiput(40.98,100.21)(-.041258,-.032691){28}{\line(-1,0){.041258}}
\multiput(39.83,99.29)(-.038715,-.032928){29}{\line(-1,0){.038715}}
\multiput(38.71,98.34)(-.036296,-.033111){30}{\line(-1,0){.036296}}
\multiput(37.62,97.34)(-.033991,-.033244){31}{\line(-1,0){.033991}}
\multiput(36.56,96.31)(-.032816,-.034405){31}{\line(0,-1){.034405}}
\multiput(35.55,95.25)(-.032654,-.036708){30}{\line(0,-1){.036708}}
\multiput(34.57,94.14)(-.0336,-.040521){28}{\line(0,-1){.040521}}
\multiput(33.63,93.01)(-.033363,-.043207){27}{\line(0,-1){.043207}}
\multiput(32.72,91.84)(-.033067,-.046045){26}{\line(0,-1){.046045}}
\multiput(31.86,90.65)(-.032706,-.049053){25}{\line(0,-1){.049053}}
\multiput(31.05,89.42)(-.033676,-.054521){23}{\line(0,-1){.054521}}
\multiput(30.27,88.17)(-.033206,-.058188){22}{\line(0,-1){.058188}}
\multiput(29.54,86.89)(-.032649,-.06213){21}{\line(0,-1){.06213}}
\multiput(28.86,85.58)(-.033678,-.069882){19}{\line(0,-1){.069882}}
\multiput(28.22,84.25)(-.032966,-.074954){18}{\line(0,-1){.074954}}
\multiput(27.62,82.9)(-.032127,-.080528){17}{\line(0,-1){.080528}}
\multiput(27.08,81.53)(-.033219,-.092475){15}{\line(0,-1){.092475}}
\multiput(26.58,80.15)(-.032128,-.100257){14}{\line(0,-1){.100257}}
\multiput(26.13,78.74)(-.0334,-.1182){12}{\line(0,-1){.1182}}
\multiput(25.73,77.33)(-.03193,-.13013){11}{\line(0,-1){.13013}}
\multiput(25.38,75.89)(-.03348,-.16031){9}{\line(0,-1){.16031}}
\multiput(25.08,74.45)(-.03138,-.18155){8}{\line(0,-1){.18155}}
\multiput(24.82,73)(-.03341,-.24337){6}{\line(0,-1){.24337}}
\multiput(24.62,71.54)(-.02992,-.29326){5}{\line(0,-1){.29326}}
\put(24.47,70.07){\line(0,-1){1.471}}
\put(24.38,68.6){\line(0,-1){2.947}}
\put(24.33,65.65){\line(0,-1){1.473}}
\put(24.39,64.18){\line(0,-1){1.47}}
\multiput(24.49,62.71)(.03142,-.2931){5}{\line(0,-1){.2931}}
\multiput(24.65,61.25)(.0297,-.20845){7}{\line(0,-1){.20845}}
\multiput(24.86,59.79)(.03231,-.18138){8}{\line(0,-1){.18138}}
\multiput(25.12,58.34)(.03087,-.14412){10}{\line(0,-1){.14412}}
\multiput(25.43,56.89)(.0326,-.12997){11}{\line(0,-1){.12997}}
\multiput(25.78,55.47)(.031386,-.108946){13}{\line(0,-1){.108946}}
\multiput(26.19,54.05)(.032641,-.100091){14}{\line(0,-1){.100091}}
\multiput(26.65,52.65)(.033691,-.092304){15}{\line(0,-1){.092304}}
\multiput(27.15,51.26)(.032538,-.080363){17}{\line(0,-1){.080363}}
\multiput(27.71,49.9)(.033349,-.074785){18}{\line(0,-1){.074785}}
\multiput(28.31,48.55)(.032333,-.066223){20}{\line(0,-1){.066223}}
\multiput(28.95,47.23)(.032966,-.061962){21}{\line(0,-1){.061962}}
\multiput(29.65,45.92)(.033503,-.058017){22}{\line(0,-1){.058017}}
\multiput(30.38,44.65)(.03254,-.052083){24}{\line(0,-1){.052083}}
\multiput(31.16,43.4)(.032956,-.048885){25}{\line(0,-1){.048885}}
\multiput(31.99,42.18)(.033302,-.045876){26}{\line(0,-1){.045876}}
\multiput(32.85,40.98)(.033584,-.043036){27}{\line(0,-1){.043036}}
\multiput(33.76,39.82)(.032641,-.038958){29}{\line(0,-1){.038958}}
\multiput(34.71,38.69)(.032842,-.03654){30}{\line(0,-1){.03654}}
\multiput(35.69,37.6)(.032991,-.034236){31}{\line(0,-1){.034236}}
\multiput(36.72,36.53)(.034161,-.03307){31}{\line(1,0){.034161}}
\multiput(37.77,35.51)(.036465,-.032925){30}{\line(1,0){.036465}}
\multiput(38.87,34.52)(.038883,-.03273){29}{\line(1,0){.038883}}
\multiput(40,33.57)(.042959,-.033682){27}{\line(1,0){.042959}}
\multiput(41.16,32.66)(.045799,-.033407){26}{\line(1,0){.045799}}
\multiput(42.35,31.79)(.048809,-.033068){25}{\line(1,0){.048809}}
\multiput(43.57,30.97)(.052009,-.032659){24}{\line(1,0){.052009}}
\multiput(44.82,30.18)(.05794,-.033636){22}{\line(1,0){.05794}}
\multiput(46.09,29.44)(.061886,-.033108){21}{\line(1,0){.061886}}
\multiput(47.39,28.75)(.066149,-.032485){20}{\line(1,0){.066149}}
\multiput(48.71,28.1)(.074708,-.03352){18}{\line(1,0){.074708}}
\multiput(50.06,27.5)(.080288,-.032723){17}{\line(1,0){.080288}}
\multiput(51.42,26.94)(.086462,-.031784){16}{\line(1,0){.086462}}
\multiput(52.81,26.43)(.100016,-.03287){14}{\line(1,0){.100016}}
\multiput(54.21,25.97)(.108874,-.031636){13}{\line(1,0){.108874}}
\multiput(55.62,25.56)(.12989,-.0329){11}{\line(1,0){.12989}}
\multiput(57.05,25.2)(.14405,-.0312){10}{\line(1,0){.14405}}
\multiput(58.49,24.89)(.18131,-.03272){8}{\line(1,0){.18131}}
\multiput(59.94,24.62)(.20838,-.03018){7}{\line(1,0){.20838}}
\multiput(61.4,24.41)(.29303,-.03209){5}{\line(1,0){.29303}}
\put(62.86,24.25){\line(1,0){1.47}}
\put(64.33,24.14){\line(1,0){1.473}}
\put(65.81,24.08){\line(1,0){1.474}}
\put(67.28,24.08){\line(1,0){1.473}}
\put(68.75,24.12){\line(1,0){1.471}}
\multiput(70.23,24.22)(.29333,.02925){5}{\line(1,0){.29333}}
\multiput(71.69,24.36)(.24344,.03285){6}{\line(1,0){.24344}}
\multiput(73.15,24.56)(.18162,.03096){8}{\line(1,0){.18162}}
\multiput(74.61,24.81)(.16038,.03311){9}{\line(1,0){.16038}}
\multiput(76.05,25.11)(.1302,.03163){11}{\line(1,0){.1302}}
\multiput(77.48,25.45)(.11827,.03313){12}{\line(1,0){.11827}}
\multiput(78.9,25.85)(.10033,.031898){14}{\line(1,0){.10033}}
\multiput(80.31,26.3)(.092551,.033006){15}{\line(1,0){.092551}}
\multiput(81.69,26.79)(.080602,.031942){17}{\line(1,0){.080602}}
\multiput(83.06,27.34)(.07503,.032794){18}{\line(1,0){.07503}}
\multiput(84.41,27.93)(.069959,.033518){19}{\line(1,0){.069959}}
\multiput(85.74,28.56)(.062205,.032506){21}{\line(1,0){.062205}}
\multiput(87.05,29.25)(.058264,.033072){22}{\line(1,0){.058264}}
\multiput(88.33,29.97)(.054598,.033551){23}{\line(1,0){.054598}}
\multiput(89.59,30.74)(.049128,.032593){25}{\line(1,0){.049128}}
\multiput(90.82,31.56)(.046121,.032961){26}{\line(1,0){.046121}}
\multiput(92.01,32.42)(.043283,.033264){27}{\line(1,0){.043283}}
\multiput(93.18,33.31)(.040598,.033507){28}{\line(1,0){.040598}}
\multiput(94.32,34.25)(.038051,.033693){29}{\line(1,0){.038051}}
\multiput(95.42,35.23)(.03448,.032737){31}{\line(1,0){.03448}}
\multiput(96.49,36.24)(.033322,.033915){31}{\line(0,1){.033915}}
\multiput(97.53,37.3)(.033195,.03622){30}{\line(0,1){.03622}}
\multiput(98.52,38.38)(.033017,.038639){29}{\line(0,1){.038639}}
\multiput(99.48,39.5)(.032785,.041183){28}{\line(0,1){.041183}}
\multiput(100.4,40.66)(.032496,.043863){27}{\line(0,1){.043863}}
\multiput(101.27,41.84)(.033429,.048563){25}{\line(0,1){.048563}}
\multiput(102.11,43.05)(.033043,.051765){24}{\line(0,1){.051765}}
\multiput(102.9,44.3)(.032583,.055181){23}{\line(0,1){.055181}}
\multiput(103.65,45.57)(.033565,.061639){21}{\line(0,1){.061639}}
\multiput(104.36,46.86)(.032974,.065907){20}{\line(0,1){.065907}}
\multiput(105.02,48.18)(.032279,.070539){19}{\line(0,1){.070539}}
\multiput(105.63,49.52)(.033316,.080043){17}{\line(0,1){.080043}}
\multiput(106.2,50.88)(.032423,.086224){16}{\line(0,1){.086224}}
\multiput(106.71,52.26)(.03361,.09977){14}{\line(0,1){.09977}}
\multiput(107.19,53.66)(.032441,.108637){13}{\line(0,1){.108637}}
\multiput(107.61,55.07)(.03104,.11884){12}{\line(0,1){.11884}}
\multiput(107.98,56.49)(.03227,.14382){10}{\line(0,1){.14382}}
\multiput(108.3,57.93)(.03028,.16094){9}{\line(0,1){.16094}}
\multiput(108.57,59.38)(.03172,.20815){7}{\line(0,1){.20815}}
\multiput(108.8,60.84)(.02855,.24399){6}{\line(0,1){.24399}}
\multiput(108.97,62.3)(.0301,.3672){4}{\line(0,1){.3672}}
\put(109.09,63.77){\line(0,1){2.729}}
\multiput(45.5,103)(.0337389381,-.0854535398){904}{\line(0,-1){.0854535398}}
\multiput(74.75,108)(.033712121,-.106439394){660}{\line(0,-1){.106439394}}
\multiput(32.25,91.5)(-.03370787,-.21067416){178}{\line(0,-1){.21067416}}
\multiput(30.75,44.5)(.051903114,-.033737024){578}{\line(1,0){.051903114}}
\multiput(26.18,79.18)(.13755,.03185){7}{\line(1,0){.13755}}
\multiput(28.11,79.63)(.13755,.03185){7}{\line(1,0){.13755}}
\multiput(30.03,80.07)(.13755,.03185){7}{\line(1,0){.13755}}
\multiput(31.96,80.52)(.13755,.03185){7}{\line(1,0){.13755}}
\multiput(33.88,80.96)(.13755,.03185){7}{\line(1,0){.13755}}
\multiput(35.81,81.41)(.13755,.03185){7}{\line(1,0){.13755}}
\multiput(37.73,81.86)(.13755,.03185){7}{\line(1,0){.13755}}
\multiput(39.66,82.3)(.13755,.03185){7}{\line(1,0){.13755}}
\multiput(41.59,82.75)(.13755,.03185){7}{\line(1,0){.13755}}
\multiput(43.51,83.19)(.13755,.03185){7}{\line(1,0){.13755}}
\multiput(45.44,83.64)(.13755,.03185){7}{\line(1,0){.13755}}
\multiput(47.36,84.09)(.13755,.03185){7}{\line(1,0){.13755}}
\multiput(49.29,84.53)(.13755,.03185){7}{\line(1,0){.13755}}
\multiput(51.21,84.98)(.13755,.03185){7}{\line(1,0){.13755}}
\multiput(53.14,85.42)(.13755,.03185){7}{\line(1,0){.13755}}
\multiput(55.06,85.87)(.13755,.03185){7}{\line(1,0){.13755}}
\multiput(56.99,86.31)(.13755,.03185){7}{\line(1,0){.13755}}
\multiput(58.92,86.76)(.13755,.03185){7}{\line(1,0){.13755}}
\multiput(60.84,87.21)(.13755,.03185){7}{\line(1,0){.13755}}
\multiput(62.77,87.65)(.13755,.03185){7}{\line(1,0){.13755}}
\multiput(64.69,88.1)(.13755,.03185){7}{\line(1,0){.13755}}
\multiput(66.62,88.54)(.13755,.03185){7}{\line(1,0){.13755}}
\multiput(68.54,88.99)(.13755,.03185){7}{\line(1,0){.13755}}
\multiput(70.47,89.44)(.13755,.03185){7}{\line(1,0){.13755}}
\multiput(72.4,89.88)(.13755,.03185){7}{\line(1,0){.13755}}
\multiput(74.32,90.33)(.13755,.03185){7}{\line(1,0){.13755}}
\multiput(76.25,90.77)(.13755,.03185){7}{\line(1,0){.13755}}
\multiput(78.17,91.22)(.13755,.03185){7}{\line(1,0){.13755}}
\multiput(80.1,91.67)(.13755,.03185){7}{\line(1,0){.13755}}
\multiput(82.02,92.11)(.13755,.03185){7}{\line(1,0){.13755}}
\multiput(83.95,92.56)(.13755,.03185){7}{\line(1,0){.13755}}
\multiput(85.88,93)(.13755,.03185){7}{\line(1,0){.13755}}
\multiput(87.8,93.45)(.13755,.03185){7}{\line(1,0){.13755}}
\multiput(89.73,93.9)(.13755,.03185){7}{\line(1,0){.13755}}
\multiput(91.65,94.34)(.13755,.03185){7}{\line(1,0){.13755}}
\multiput(93.58,94.79)(.13755,.03185){7}{\line(1,0){.13755}}
\multiput(95.5,95.23)(.13755,.03185){7}{\line(1,0){.13755}}
\multiput(44.18,30.68)(.063112,.032779){14}{\line(1,0){.063112}}
\multiput(45.95,31.6)(.063112,.032779){14}{\line(1,0){.063112}}
\multiput(47.71,32.52)(.063112,.032779){14}{\line(1,0){.063112}}
\multiput(49.48,33.43)(.063112,.032779){14}{\line(1,0){.063112}}
\multiput(51.25,34.35)(.063112,.032779){14}{\line(1,0){.063112}}
\multiput(53.02,35.27)(.063112,.032779){14}{\line(1,0){.063112}}
\multiput(54.78,36.19)(.063112,.032779){14}{\line(1,0){.063112}}
\multiput(56.55,37.1)(.063112,.032779){14}{\line(1,0){.063112}}
\multiput(58.32,38.02)(.063112,.032779){14}{\line(1,0){.063112}}
\multiput(60.08,38.94)(.063112,.032779){14}{\line(1,0){.063112}}
\multiput(61.85,39.86)(.063112,.032779){14}{\line(1,0){.063112}}
\multiput(63.62,40.78)(.063112,.032779){14}{\line(1,0){.063112}}
\multiput(65.39,41.69)(.063112,.032779){14}{\line(1,0){.063112}}
\multiput(67.15,42.61)(.063112,.032779){14}{\line(1,0){.063112}}
\multiput(68.92,43.53)(.063112,.032779){14}{\line(1,0){.063112}}
\multiput(70.69,44.45)(.063112,.032779){14}{\line(1,0){.063112}}
\multiput(72.45,45.36)(.063112,.032779){14}{\line(1,0){.063112}}
\multiput(74.22,46.28)(.063112,.032779){14}{\line(1,0){.063112}}
\multiput(75.99,47.2)(.063112,.032779){14}{\line(1,0){.063112}}
\multiput(77.76,48.12)(.063112,.032779){14}{\line(1,0){.063112}}
\multiput(79.52,49.04)(.063112,.032779){14}{\line(1,0){.063112}}
\multiput(81.29,49.95)(.063112,.032779){14}{\line(1,0){.063112}}
\multiput(83.06,50.87)(.063112,.032779){14}{\line(1,0){.063112}}
\multiput(84.82,51.79)(.063112,.032779){14}{\line(1,0){.063112}}
\multiput(86.59,52.71)(.063112,.032779){14}{\line(1,0){.063112}}
\multiput(88.36,53.62)(.063112,.032779){14}{\line(1,0){.063112}}
\multiput(90.12,54.54)(.063112,.032779){14}{\line(1,0){.063112}}
\multiput(91.89,55.46)(.063112,.032779){14}{\line(1,0){.063112}}
\multiput(93.66,56.38)(.063112,.032779){14}{\line(1,0){.063112}}
\multiput(95.43,57.3)(.063112,.032779){14}{\line(1,0){.063112}}
\multiput(97.19,58.21)(.063112,.032779){14}{\line(1,0){.063112}}
\multiput(98.96,59.13)(.063112,.032779){14}{\line(1,0){.063112}}
\multiput(100.73,60.05)(.063112,.032779){14}{\line(1,0){.063112}}
\multiput(102.49,60.97)(.063112,.032779){14}{\line(1,0){.063112}}
\multiput(104.26,61.89)(.063112,.032779){14}{\line(1,0){.063112}}
\multiput(106.03,62.8)(.063112,.032779){14}{\line(1,0){.063112}}
\multiput(107.8,63.72)(.063112,.032779){14}{\line(1,0){.063112}}
\put(80,91.75){\circle*{1}} \put(91.5,55.25){\circle*{1}}
\put(52.25,85.25){\circle*{1.41}} \put(68.75,43.75){\circle*{1.12}}
\put(94.25,53.75){\makebox(0,0)[cc]{$o_1$}}
\put(83.5,88.75){\makebox(0,0)[cc]{$o_2$}}
\put(48.75,80){\makebox(0,0)[cc]{$o_2'$}}
\put(61.75,45.5){\makebox(0,0)[cc]{$o_1'$}}
\put(82.75,24){\makebox(0,0)[cc]{$\ccc$}}
\end{picture}
\caption{Lemma \ref{deletion}} \label{fig4}
\end{figure}

(b) Suppose $\ccc$ is a $\QQ$-chord. Then the weight of any pair
$(o_1,o_2)$ in $\gc_0$ cannot exceed the weight of that pair in
$\gc$. Indeed, if the neighbor $o_i'$ of $o_i$ in $\gc$ did not
belong to $\ccc$, then the weight of that neighbor is the same or
smaller in $\gc_0$ (since that neighbor can only become closer to
the boundary). If $o_i'$ belongs to $\ccc$ and its weight in $\gc$
is $\alpha_k$, then the weight of the neighbor of $o_i$ in $\gc_0$
cannot be bigger than $\alpha_k$ (because the neighbor of $o_i$ in
$\gc_0$ is closer to the boundary than the neighbor of $o_i$ in
$\gc$).

Let $X_0$ be the set of bad pairs in $\gc_0$. Then $X_0$ is the
union of the set $Y_1$ of all bad pairs of $\gc$ which do not belong
to $\ccc$ (the weight of that pair in $\gc_0$ does not exceed its
weight in $\gc$, we noted that in the previous paragraph) and the
set $Y_2$ of good pairs $(o_1,o_2)$ of $\gc$ whose neighbors
$o_1',o_2'$ belong to $\ccc$ but the pair $(o_1',o_2')$ is bad in
$\gc$ (the weight of $(o_1,o_2)$ in $\gc_0$ equals the weight of
$(o_1',o_2')$ in $\gc$). Since $(o_1',o_2')$ is uniquely determined
by $(o_1,o_2)$, the total weight of the pairs from $X_0=Y_1\cup Y_2$
in $\gc_0$ does not exceed the total weight of bad pairs from $\ccc$
in $\gc$.  Hence $\eee_\alpha(\gc_0)\le \eee_\alpha(\gc)$.
\endproof

\begin{df}\label{dfclose}{\rm
We say that a $\QQ$-chord  $\ccc'$ is {\em close } to a $\QQ$-chord
$\ccc$ of a BCD $\gc$ if every chord crossing $\ccc'$ also crosses
$\ccc$.}
\end{df}

From now on let $$\alpha=(\frac1K,\frac2K,...,1).$$ We shall write
$\eee(\gc)$ instead of $\eee_{\alpha}(\gc)$ and call it simply {\em
dispersion} of the grading $\gc$.

\begin{lemma}\label{close}
Suppose that a $\QQ$-chord $\ccc'$ is close to a $\QQ$-chord
$\ccc$ in a BCD $\gc$ of a disc $D$. Let $D_0$ be one of the two
subdiscs of the disc $D$ divided by $\ccc$ that contains $\ccc'$.
Suppose that the number of nodes from $D_0$ on any $\TT$-chord in
$\gc$ is at most $K$. Denote the numbers of nodes on $\ccc$ and
$\ccc'$ by $l$ and $l'$, respectively. Let ${\gc}_0$ be the BCD
obtained from $\gc$ by removing $\ccc'$. Then
\begin{equation}\label{alph}
\eee(\gc_0)\le \eee(\gc)-\frac1{K^2} l'(l-l').\end{equation}
\end{lemma}

\proof Let $V_1$ be the set of nodes on $\ccc$ that belong to
$\TT$-chords intersecting $\ccc'$, $V_2$ be the set of other nodes
on $\ccc$. Then $|V_1|=l'$, $|V_2|=l-l'$, and the number of pairs
$(o,o')$ of nodes on $\ccc$ such that $\ccc'$ is to the left of
$\ccc$ according to the direction $o\to o'$ and either ($o\in V_1$
and $o'\in V_2$) or ($o\in V_2$ and $o'\in V_1$),  is $l'(l-l')$.

Consider one of such pairs, $(o,o')$. Let $o=o_0,o_1,o_2,...$ be a
sequence of nodes from $D_0$ on the $\TT$-chord $\ttt$ passing
through $o$. Similarly consider a sequence of nodes $o'=o_0',o_1',
o_2'...$ from $D_0$ on a $\TT$-chord $\ttt'$. Since one of the nodes
$o$ or $o'$ is in $V_1$ and another one in $V_2$, there exists
$i=0,1,...$ such that $(o_i,o_i')$ belongs to some $\QQ$-chord
$\ccc_i\ne \ccc'$ and is bad.

The weight in $\gc_0$ of the node from the pair $(o_i,o_i')$ that
belongs to a $\TT$-chord that does not cross $\ccc'$, is the same as
its weight in $\gc$ and is at least $\frac{1}{K}$.  The weight of
the other node from that pair decreases by at least $\frac{1}{K}$
when we pass from $\gc$ to $\gc_0$ since by the assumption the
number of nodes on every $\TT$-chord in $D_0$ is at most $K$ and we
removed $\ccc'$ that was to the left of $\ccc_i$. Hence the weight
of the pair decreases by at least $\frac1{K^2}$. Since the number of
such pairs $(o_i,o_i')$ is $l'(l-l')$, the total weight of such
pairs decreases by at least $\frac{1}{K^2}l'(l-l')$.

Now as in Lemma \ref{deletion} let $X_0$ be the set of bad pairs
in $\gc_0$. Then $X_0$ is the union of the set $Y_1$ of all bad
pairs of $\gc$ which do not belong to $\ccc'$ and the set $Y_2$ of
good pairs $(o_1,o_2)$ of $\gc$ whose neighbors $o_1',o_2'$ belong
to $\ccc'$ but the pair $(o_1',o_2')$ is bad in $\gc$ (the weight
of $(o_1,o_2)$ in $\gc_0$ equals the weight of $(o_1',o_2')$ in
$\gc$). As in the proof of Lemma \ref{deletion} the weight of any
pair of $Y_1$ in $\gc_0$ does not exceed its weight in $\gc$.
Since all pairs $(o_i,o_i')$ that we considered in the previous
paragraph belong to $Y_1$, we get the inequality (\ref{alph}).
\endproof

\section{The upper bound of the Dehn function}
\label{upperb}

Lemma \ref{NoAnnul} implies that we can associate a bipartite
chord diagram $\gc$ to any reduced diagram $\Delta$, where
$\TT$-chords are the medians of the maximal $\theta$-bands, and
$\QQ$-chords are the medians of $q$-bands.

\begin{rk} \label{rk2}
{\rm If $\Delta$ is not a topological disc then $\partial\Delta$
can be transformed into a circle by an arbitrary small
deformation, and so the topological structure of $\gc$ is
well-defined.}
\end{rk}

The dispersion of this BCD is called the {\em dispersion}
$\eee=\eee(\Delta)$ of the diagram $\Delta$.

In the following lemma, we estimate the area of a \vk diagram over
$\sss\circ Z$ in terms of its perimeter and dispersion. Namely we
show that for some constant $M$ the area of any reduced diagram
$\Delta$ of perimeter $n$ does not exceed $Mn^2\log n+M\eee(\Delta)$
for some constant $M$. (Then using the quadratic upper bound for
$\eee(\Delta)$ we will deduce that the area is bounded by $M'n^2\log
n$ for some constant $M'$.) Roughly speaking, we are doing the
following. We use induction on the perimeter of the diagram. First
we remove rim $\theta$-bands (those with one side on the boundary of
the diagram) with short bases. This operation decreases the
perimeter and preserves the sign of $Mn^2\log
n+M\eee(\Delta)-\area(\Delta)$, so we can assume that the diagram
does not have such bands. Then we use Lemma \ref{est'takaya} and
find a tight comb inside the diagram with a handle $\ccc$. We also
find a long enough $q$-band $\ccc'$ that is close to $\ccc$. We use
a surgery which amounts to removing a part of the diagram between
$\ccc'$ and $\ccc$ and then gluing the two remaining parts of
$\Delta$ together. The main difficulty is to show that, as a result
of this surgery, the perimeter decreases and the area and the
dispersion change in such a way that the expression $Mn^2\log
n+M\eee(\Delta)-\area(\Delta)$ does not change its sign. In the
proof, we need to consider several cases depending on the shape of
the subdiagram between $\ccc'$ and $\ccc$. Note that neither
$Mn^2\log n$ nor $M\eee(\Delta)$ nor $\area(\Delta)$ alone behave in
the appropriate way as a result of the surgery, but the expression
$Mn^2\log n+M\eee(\Delta)-\area(\Delta)$ behaves as needed.

Let us take a big enough constant $M$. Here ``big enough" means that
$M$ satisfies the inequalities used in the proof of Lemma \ref{main}
(i.e. (\ref{param3}),(\ref{param4}), (\ref{param5}), (\ref{param6}),
(\ref{param7}), (\ref{param8}), (\ref{param9}), (\ref{param10})).
Each of them has the form $M>C$ for some constant $C$ that does not
depend on $M$ (but depends on the constants introduced earlier), and
the number of inequalities is finite, so the choice of $M$ is
possible.

\begin{lemma}\label{main} The area of a reduced diagram $\Delta$
does not exceed $Mn^2\log' n +M\eee(\Delta)$ where
$n=|\partial\Delta|$, and $\log' n =\max (\log_2 n, 1)$.
\end{lemma}

\proof Arguing by contradiction, we consider a counter-example
$\Delta$ with minimal perimeter $n$. Of course, its area is
positive, and, by Lemma \ref{NoAnnul}, we have at least 2
$\theta$-edges on the boundary $\partial\Delta$, and so $n\ge 2$.

{\bf Step 1.} Assume that there are two $\QQ$-chords $\ccc$ and
$\ccc'$, where $\ccc'$ is close to $\ccc$ in the BCD $\gc(\Delta)$,
and suppose the number $l'$ of the nodes on $\ccc'$ does not exceed
a half of the number $l$ of nodes lying on $\ccc$. Let $\xxx$ and
$\xxx'$ be the $q$-bands corresponding in $\Delta$ to $\ccc$ and
$\ccc'$, respectively. Without loss of generality assume that the
top subdiagram $\Gamma$ with respect to $\xxx$ contains both $\xxx$
and $\xxx'$ and the top subdiagram $\Gamma'<\Gamma$ with respect to
$\xxx'$ does not contain $\xxx$. We suppose in addition that
$\Gamma'$ is a comb with handle $\xxx'$. Let us prove that then the
basic width $b$ of the comb $\Gamma'$ is not smaller than $K$, and
in particular the comb cannot be tight.

By contradiction suppose that $b<K$.

We note that the bands $\xxx$ and $\xxx'$ contain $l$ and $l'$
cells, respectively. It follows from Lemma \ref{comb} and our
assumptions, that the area of $\Gamma'$ does not exceed $
C_1(l')^2+2\alpha l'$ for $\alpha =|\Gamma'|_a$ and some constant
$C_1$.

Let $\Delta'$ be the diagram obtained by deleting the subdiagram
$\Gamma'$ from $\Delta$. Since the boundary of $\Delta'$ has at
least two $q$-edges fewer than $\Delta$, we have $|\partial
\Delta'|\le |\partial\Delta|-2$. Moreover, we have from Lemma
\ref{ochev} and Lemma \ref{NoAnnul} that
\begin{equation}\label{eq679}
|\partial \Delta|-|\partial\Delta'|\ge\gamma =\max(2,
\delta(\alpha-Ll'))
\end{equation}
because the top or the bottom of $\xxx'$ has at most $Ll'$
$a$-edges.

The difference of the dispersions $\eee(\Delta)-\eee(\Delta')$ is at
least $\frac1{K^2} l'(l-l')$ by Lemmas \ref{close} and
\ref{deletion}. Hence $\eee(\Delta)-\eee(\Delta')\ge \frac1{K^2}
(l')^2$ as $l'\le l-l'$. This inequality and the inductive
assumption related to the area of $\Delta'$, imply that the area of
$\Delta'$ is not greater than
$$M(n-\gamma)^2\log'(n-\gamma)+ M\eee(\Delta)-\frac{M}{K^2} (l')^2.$$
Adding the area of $\Gamma'$, we see that the area of $\Delta$ does
not exceed $$Mn^2\log'n +M\eee(\Delta)- M\gamma n \log'n
-\frac{M}{K^2} (l')^2+C_1(l')^2+2\alpha l'.$$ (Keep in mind that
$\gamma\le n$.) This will contradict the choice of the
counter-example $\Delta$ when we prove that
\begin{equation}\label{nado}
- M\gamma n \log'n -\frac{M}{K^2}(l')^2+C_1(l')^2+2\alpha l'<0
\end{equation}
 Consider two cases.

(a) Let $\alpha\le 2Ll'$. Then inequality (\ref{nado}) follows from
the inequality
\begin{equation}\label{param3}
M\ge 2K^2(C_1+ 4L).
\end{equation}

(b) Assume that $\alpha> 2Ll'$. Then by (\ref{eq679}) we have
$\gamma\ge \frac12 \delta\alpha$ and $M\gamma n\log'n >2\alpha l'$
since $n\ge 2l'$ by Lemma \ref{NoAnnul}, and

\begin{equation}\label{param4}
M> 2\delta^{-1}.
\end{equation}
Since $\frac{M}{K^2}(l')^2>C_1(l')^2$ by (\ref{param3}), the
inequality (\ref{nado}) follows.

{\bf Step 2.} Assume that $\Delta$ has a rim $\theta$-band $\ttt$
whose base has  $s\le K$ letters and $\topp(\ttt)$ is in
$\partial(\Delta)$. By deleting $\ttt$, we obtain, by Lemma
\ref{rim}, a diagram $\Delta'$ with $ |\partial\Delta'|\le
n-\delta$. Since $\topp(\ttt)$ lies on $\partial\Delta$, we have
from the definition of the length (Section \ref{length}), that the
number of $a$-edges in $\topp(\ttt)$ is less than
$\delta^{-1}(n-Ls)$. By Lemma \ref{band}, the length of $\ttt$ is at
most $(L+1)s+\delta^{-1}(n-Ls)< \delta^{-1}n$ since $\delta\iv >
\frac{L+1}{L}$ by (\ref{param}). Thus, by applying the inductive
hypothesis to $\Delta'$, we have that area of $\Delta$ is not
greater than $M(n-\delta)^2 \log' (n-\delta)+M\eee(\Delta)
+\delta^{-1}n$ because $\eee(\Delta')\le\eee(\Delta)$ by Lemma
\ref{deletion} (a). But this sum does not exceed $Mn^2\log' n
+M\eee(\Delta)$ provided

\begin{equation}\label{param5}
M\ge\delta^{-2}.
\end{equation}

This contradicts the choice of $\Delta$. Hence the base of every rim
$\theta$-band of $\Delta$ has more than $K$ letters.

{\bf Step 3.} Now we can apply Lemma \ref{est'takaya}. By that
lemma, there exists a tight comb $\Gamma<\Delta$. Let $\ttt$ be a
$\theta$-band of $\Gamma$ with $\bbb$-tight base.

In particular, the basic width of $\Gamma$ is smaller than $K$.
Since the base of $\ttt$ is $\bbb$-tight, it is equal to $uxvx$ for
some $x=Q_i$ (we read the base starting at the boundary of $\Delta$)
where $\base(xvx)$ is $\bbb$-covered.

The second occurrence of $x$ in $uxvx$ corresponds to the last cell
of $\ttt$ belonging to the $x$-band $\qq$. Let $\qq'$ be the maximal
$x$-band of $\Gamma$ crossing $\ttt$ at the cell corresponding to
the first occurrence of $x$ in $uxvx$.

We consider the smallest subdiagram $\Gamma'$ of $\Delta$ containing
all the $\theta$-bands of $\Gamma$ crossing the $x$-band $\qq'$. It
is a comb with handle $\qq_2\subset\qq$. The comb $\Gamma'$ is
covered by a trapezium $\Gamma_2$ placed between $\qq'$ and $\qq$,
and a comb $\Gamma_1$ with handle $\qq'$. The $Q_i$-band $\qq'$
belongs to both $\Gamma_1$ and $\Gamma_2$. The remaining part of
$\Gamma$ is a disjoint union of two combs $\Gamma_3$ and $\Gamma_4$
whose handles $\qq_3$ and $\qq_4$ contain the cells of $\qq$ that do
not belong to the trapezium $\Gamma_2$. The handle of $\Gamma$ is
the composition of handles $\qq_3$, $\qq_2$, $\qq_4$ of $\Gamma_3$,
$\Gamma'$ and $\Gamma_4$ in that order.

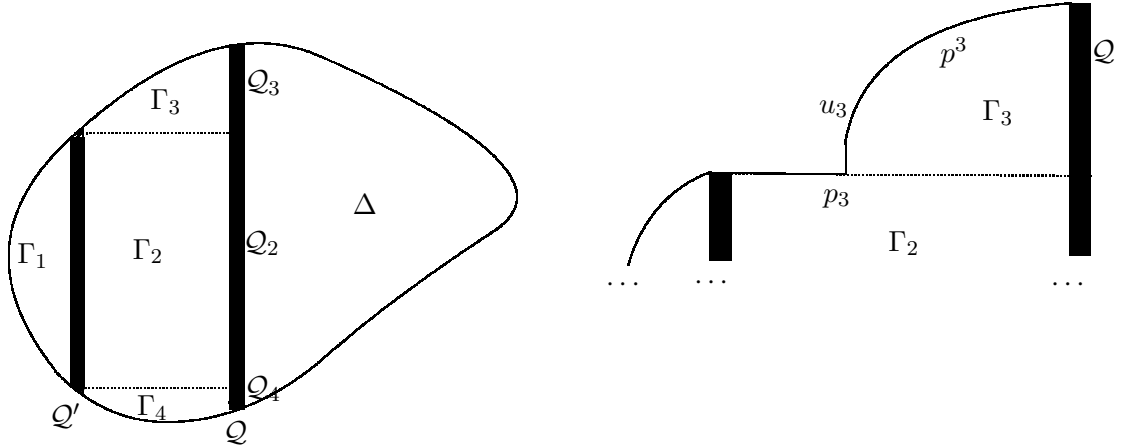
\begin{figure}[ht]
\centering

\unitlength .5mm 
\linethickness{0.4pt}
\ifx\plotpoint\undefined\newsavebox{\plotpoint}\fi 
\begin{picture}(311.87,123.97)(20,0)
\qbezier(33,26.46)(6.84,59.61)(36.57,87.41)
\qbezier(36.57,87.41)(74.03,122.49)(103.16,110)
\qbezier(103.16,110)(173.77,79.23)(149.84,63.33)
\qbezier(149.84,63.33)(121.89,45.04)(103.46,28.54)
\qbezier(103.46,28.54)(83.69,12.19)(62.14,12.49)
\qbezier(62.14,12.49)(45.78,12.78)(33,26.16)
\put(37.76,22){\rule{2.97\unitlength}{.3\unitlength}}
\put(37.46,21.7){\rule{3.57\unitlength}{66.6\unitlength}}
\put(40.73,85.92){\line(0,1){4.757}}
\multiput(40.73,90.68)(-.0335663,-.0479519){62}{\line(0,-1){.0479519}}
\multiput(38.65,87.7)(.0330335,.0330335){45}{\line(0,1){.0330335}}
\put(40.14,89.19){\line(-1,0){1.189}}
\put(38.95,89.19){\line(1,0){1.189}}
\put(40.73,19.32){\line(0,1){4.162}}
\multiput(40.73,23.49)(-.03303,.5946){9}{\line(0,1){.5946}}
\put(39.84,21.7){\circle*{3.03}}
\put(79.38,112.38){\rule{2.97\unitlength}{0\unitlength}}
\put(79.38,17.54){\rule{3.57\unitlength}{.3\unitlength}}
\put(79.68,15.76){\rule{3.86\unitlength}{97.22\unitlength}}
\multiput(38.58,89.12)(.9844,.00661){46}{{\rule{.4pt}{.4pt}}}
\multiput(38.28,21.34)(.9844,0){46}{{\rule{.4pt}{.4pt}}}
\put(27.35,56.78){\makebox(0,0)[cc]{$\Gamma_1$}}
\put(62.73,97.81){\makebox(0,0)[cc]{$\Gamma_3$}}
\put(57.97,57.97){\makebox(0,0)[cc]{$\Gamma_2$}}
\put(59.16,16.35){\makebox(0,0)[cc]{$\Gamma_4$}}
\put(88.6,102.87){\makebox(0,0)[cc]{$\qq_3$}}
\put(88.22,60.05){\makebox(0,0)[cc]{$\qq_2$}}
\put(88.51,21.41){\makebox(0,0)[cc]{$\qq_4$}}
\put(115.65,70.46){\makebox(0,0)[cc]{$\Delta$}}
\put(81.16,9.81){\makebox(0,0)[cc]{$\qq$}}
\put(35.97,14.27){\makebox(0,0)[cc]{$\qq'$}}
\qbezier(185.52,54.41)(191.17,72.99)(207.52,79.08)
\multiput(207.15,78.42)(.99392,-.00583){103}{{\rule{.4pt}{.4pt}}}
\put(243.49,78.49){\line(0,1){8.027}}
\put(302.95,56.78){\rule{5.65\unitlength}{67.19\unitlength}}
\multiput(207.52,78.78)(3.96402,-.03303){9}{\line(1,0){3.96402}}
\qbezier(303.25,123.97)(248.25,119.52)(243.19,86.51)
\put(184.92,49.05){\makebox(0,0)[cc]{$\cdots$}}
\put(208.41,49.65){\makebox(0,0)[cc]{$\cdots$}}
\put(303.25,49.05){\makebox(0,0)[cc]{$\cdots$}}
\put(283.92,94.84){\makebox(0,0)[cc]{$\Gamma_3$}}
\put(240.33,95.89){\makebox(0,0)[cc]{$u_3$}}
\put(272.33,110.89){\makebox(0,0)[cc]{$p^3$}}
\put(241.11,71.95){\makebox(0,0)[cc]{$p_3$}}
\put(258.95,60.05){\makebox(0,0)[cc]{$\Gamma_2$}}
\put(311.87,111.19){\makebox(0,0)[cc]{$\qq$}}
\put(207.22,55.6){\rule{5.95\unitlength}{22.89\unitlength}}
\end{picture}

\caption{Step 3 in Lemma \ref{main}.} \label{fig9}
\end{figure}

Let the lengths of $\qq_3$ and $\qq_4$ be $l_3$ and $l_4$,
respectively. Let $l'$ be the length of the handle of $\Gamma'$.

Then $$l=l'+l_3+l_4,$$ and, as we proved in Step 1, $l'>l/2$.

For $i\in \{3,4\}$ and $\alpha_i=|\partial\Gamma_i|_a$, Lemma
\ref{comb} gives inequalities
\begin{equation}\label{gamma34}
A_i\le C_1l_i^2+2\alpha_i l_i
\end{equation}
for some constant $C_1$ where $A_i$ is the area of $\Gamma_i$. Let
$p_3, p_4$ be the top and the bottom of the trapezium $\Gamma_2$.
Here $p_3$ (resp. $p_4$) share some initial edges with $\partial
\Gamma_3$ (with $\partial \Gamma_4$), the rest of these paths belong
to the boundary of $\Delta$. We denote by $d_3$ the number of
$a$-edges of $p_3$ and by $d'_3$ the number of its edges which do
not belong to $\Gamma_3$. Similarly, we introduce $d_4$ and $d'_4$.
Let $A_2$ be the area of $\Gamma_2$. Then, by Lemma \ref{area},
\begin{equation}\label{gamma2}
A_2\le C_2l'(d_3+d_4+\log l'+1)
\end{equation}
for some constant $C_2$.

Now we note that the handle $\qq_2$ of $\Gamma'$ is a copy of $\qq'$
because both maximal $q$-bands of the trapezium $\Gamma_2$
correspond to the same basic letter $x$. This makes the following
surgery possible. The diagram $\Delta$ is covered by two
subdiagrams: $\Gamma$ and another subdiagram $\Delta_1$, having only
the band $\qq$ in common. We construct a new auxiliary diagram by
attaching $\Gamma_1$ to $\Delta_1$ with identification of the band
$\qq'$ of $\Gamma_1$ and the band $\qq_2$. We denote the constructed
diagram by $\Delta_0$. It is a reduced diagram because every pair of
its cells having a common edge, has a copy either in $\Gamma_1$ or
in $\Delta_1$. It follows from our constructions that the area of
$\Delta$ does not exceed $A_2+A_3+A_4+A_0$, where $A_0$ is the area
of $\Delta_0$.

Let $p^3$ be the segment of the boundary $\partial\Gamma_3$ that
joins $\qq$ and $\Gamma_2$ along the boundary of $\Delta$. It
follows from the definition of $d_3$, $d'_3$, $l_3$ and $\alpha_3$,
that the number of $a$-edges lying on $p^3$ is at least $\alpha_3
-(d_3-d'_3)-Ll_3$.

Let $u_3$ be the part of $\partial\Delta$ that contains $p^3$ and
connects $\qq$ with $\qq'$. It has $l_3$ $\theta$-edges. Hence we
have, by Lemma \ref{ochev}, that the length $|u_3|$ of $u_3$ is at
least
$$\max(l_3, l_3+\delta(|p^3|_a-Ll_3))\ge
\max(l_3, l_3+\delta(\alpha_3 -(d_3-d'_3)-2Ll_3)).$$ Since $u_3$
includes a subpath of length $d'_3$ having no $\theta$-edges, we
also have by inequality (\ref{delta}) that $|u_3|\ge
l_3+\delta(d'_3-L)$.

One can similarly define $p^4$ and $u_4$ for $\Gamma_4$. When
passing from $\partial\Delta$ to $\partial\Delta_0$ we replace the
end edges of $\qq'$, $u_3$ and $u_4$ by two subpaths of
$\partial\qq$ having lengths $l_3$ and $l_4$. Let
$n_0=|\partial\Delta|$. Then it follows from the previous paragraph
that

\begin{equation}\label{nbezn0}
n-n_0\ge2+\delta(\max(0,d'_3-L, \alpha_3-(d_3-d'_3)-2Ll_3)+\max(0,
d'_4-L,\alpha_4-(d_4-d'_4)-2Ll_4))
\end{equation}

In particular, $n_0\le n-2$. By the inductive hypothesis,
\begin{equation}\label{A0}
A_0\le Mn_0^2 \log' n_0+M\eee(\Delta_0)
\end{equation}

We note that the dispersion $\eee(\Delta_0)$ of $\Delta_0$ is not
greater than $\eee(\Delta)-\frac1{K^2} l'(l-l')$ by Lemmas
\ref{close} and \ref{deletion}(b).

Therefore, by inequality (\ref{A0}), the area of $\Delta$ is not
greater than
\begin{equation}\label{ADelta}
Mn^2\log' n +M\eee(\Delta)- M n(n-n_0)\log'n -\frac{M}{K^2}
l'(l-l')+A_2+A_3+A_4
\end{equation}

In view of inequalities (\ref{gamma2}) and (\ref{gamma34}), to
obtain the desired contradiction, we should prove that
\begin{equation}\label{tsel'}
Mn(n-n_0)\log'n+\frac{M}{K^2} l'(l-l')\ge  C_3l'(d_3+d_4+\log'
l'+1)+C_3(l_3^2+l_4^2)+2\alpha_3 l_3 + 2\alpha_4l_4
\end{equation}
where $C_3=\max(C_1,C_2)$ is a constant that does not depend on $M$.
Note that we can assume that
\begin{equation}\label{c3}
C_3>>L.
\end{equation}

First we can choose $M$ big enough so that
$\frac{M}{3K^2}l'(l-l')\ge C_3(l_3+l_4)^2\ge C_3(l_3^2+l_4^2)$.
Indeed

\begin{equation}\label{lll}
l-l'=l_3+l_4<l'
\end{equation}
since $l'>l/2$, and $\frac{M}{3K^2} l'(l-l')\ge \frac{M}{3K^2}
(l_3+l_4)(l_3+l_4)$, so it is enough to assume that
\begin{equation}
\label{param6} M>3K^2C_3.
\end{equation}

We also have that
\begin{equation}\label{a4}
\frac{M}2 n(n-n_0)\log' n\ge C_3l' (\log'l'+1)
\end{equation}
because $n-n_0\ge 2$, $n\ge 2l'$ and $M\ge 2C_3$ by (\ref{param6}).

It remains to prove that \begin{equation}\label{t1} \frac{M}2
n(n-n_0)\log' n+\frac{2M}{3K^2}l'(l-l')>
C_3l'(d_3+d_4)+2\alpha_3l_3+2\alpha_4l_4.
\end{equation}

We assume without loss of generality that $\alpha_3\ge\alpha_4$, and
consider two cases.

\medskip

(a) Suppose $\alpha_3\le 2C_3(l-l')$.

Since $d_i\le \alpha_i+d'_i$ for $i=3,4$, we also, by inequality
(\ref{nbezn0}), have $$d_3+d_4\le
\alpha_3+\alpha_4+d_3'+d_4'<4C_3(l-l')+\delta^{-1}(n-n_0)+2L-2\delta\iv<
4C_3(l-l')+\delta\iv(n-n_0).
$$
since $\delta\iv>L$ by (\ref{param}).

Therefore

\begin{equation}\label{a1}
\frac{M}{3K^2} l'(l-l')+\frac{M}2 n(n-n_0)\log' n\ge
4C_3^2l'(l-l')+C_3\delta\iv(n-n_0)l'>
  C_3l'(d_3+d_4)
\end{equation}
since we can assume that

\begin{equation}\label{param7}
M> 12K^2C_3^2,\qquad M/2> C_3\delta^{-1}.
\end{equation}

We have also by (\ref{lll}):
\begin{equation}\label{a3}
\frac{M}{3K^2} l'(l-l')\ge \frac{M}{3K^2}(l_3+l_4)(l_3+l_4)\ge
\frac{M}{3}\frac{\alpha_3+\alpha_4}{4C_3}(l_3+l_4)>
2\alpha_3l_3+2\alpha_4l_4
\end{equation}
since we can assume that
\begin{equation}\label{param8}
M>24K^2C_3.
\end{equation}

 The sum of inequalities (\ref{a1}) and (\ref{a3}) gives us the
desired inequality (\ref{t1}).

\medskip

(b) Assume now that $\alpha_3>2C_3(l-l')$. Then, applying Lemma
\ref{comb} to the comb $\Gamma_3$, we obtain
$d_3-d'_3<\frac12\alpha_3+2CKl_3\le\frac56\alpha_3$ since $l_3\le
l-l'<\frac{\alpha_3}{2C_3}$ and
\begin{equation}\label{param9} 2C_3>12CK.\end{equation} We also have
$d_4-d'_4<\frac12\alpha_4+2CKl_4\le\frac56 \alpha_3$. These two
inequalities and inequality \ref{nbezn0} lead to
\begin{equation}\label{d34}
d_3+d_4\le \frac53\alpha_3+\delta^{-1}(n-n_0)
\end{equation}

In addition, $$\alpha_3-(d_3-d'_3)-2Ll_3\ge \frac16 \alpha_3 -
\frac{2L}{2C_3}\alpha_3\ge\frac17\alpha_3,$$ since $l_3\le
l-l'<\frac{\alpha_3}{2C_3}$ and $C_3>42L$ by (\ref{c3}). Therefore,
by \ref{nbezn0},
\begin{equation}\label{raznitsa}
n-n_0\ge \frac17 \delta\alpha_3.
\end{equation}
Thus by (\ref{d34})

\begin{equation}
\label{d10} d_3+d_4<13 \delta\iv(n-n_0).
\end{equation}

Since $2l'<n$ and $n-n_0\ge 2$, inequality (\ref{d10}) implies

\begin{equation}\label{b1}
\frac{M}3 n(n-n_0)\log'n > C_3l'(d_3+d_4)
\end{equation}
because we can assume that \begin{equation}\label{param10}
M>>C_3\delta\iv
\end{equation}
($M>21C_3\delta\iv$ is enough).

Inequalities (\ref{raznitsa}), (\ref{param10}),
$\alpha_3\ge\alpha_4$, and $4(l_3+l_4)\le n$ give us
\begin{equation}\label{b2}
\frac{M}6 n(n-n_0)\log n'\ge \frac72 C_3 \delta^{-1}(n-n_0)n\ge
2\alpha_3(l_3+l_4)\ge 2\alpha_3l_3+2\alpha_4 l_4
\end{equation}

The inequality (\ref{t1}) follows now from inequalities (\ref{b1}),
and (\ref{b2}).
\endproof

\begin{lemma}\label{combinat}
Let $n$ be the combinatorial length of a reduced diagram $\Delta$.
Then the area of $\Delta$ is $O (n^2 \log n) $.
\end{lemma}

\proof By Lemmas \ref{NoAnnul} and \ref{entropy}, we have
\begin{equation}\label{quadr}
\eee(\Delta)\le (n/2)^2.
\end{equation}

Let $n'=|\Delta|$. It follows from the definition of the length that
$n'\le n$. By this inequality, inequality (\ref{quadr}) and Lemma
\ref{main}, we have the inequality $\area(\Delta)\le Mn^2\log' n +
\frac{M}4 n^2$. The lemma is proved.
\endproof

\section{The end of the proof}
\label{lower}

\begin{lemma}\label{lb}
The Dehn function of $\sss\circ Z$ is, up to equivalence, at least
$n^2\log n$.
\end{lemma}

\proof Let us use the fact that for some fixed word $W$, $\sss$ had
the computation (\ref{tolst}) for every $n$. The length of that
computation is $2n$ and the width is $n+N$.

Consider  the corresponding computation of $\sss\circ Z$ (see the
proof of Lemma \ref{lmlll}). Its width is at most $n+2N$. Let $l(n)$
be the length of that computation and let $a(n)$ be its area. Then
by Remark \ref{rk90} and the description of the computation of
$\sss\circ Z$ (see Section \ref{compsa}), $$2\sum_{m=1}^{n-1}
(2^m+2N-3)+(2^n+2N-3)+2N-1<l(n)<2\sum_{m=1}^{n-1}(6\cdot
2^m+2N-3)+(6\cdot 2^n+2N-3)+2N-1.$$ So $3\cdot 2^n+C_1\le
l(n)<18\cdot 2^n+C_2$ for some constants $C_1, C_2$. In addition
$a(n)>C' l(n)\log_2 l(n)$ for some constant $C'$.

The sides of the corresponding trapezium $\Delta$ have the same
labels. Thus we can take $l(n)$ copies of $\Delta$ and glue them
side by side to obtain a new trapezium $\bar\Delta$. The perimeter
$d$ of $\bar\Delta$ is $O(l(n))$ and the area is $O(l(n)^2\log
l(n))$. By Lemmas \ref{NoAnnul} and \ref{computation}, there is
only one reduced \vk diagram with the same boundary label as
$\bar\Delta$.

Since $l(n)$ is between $3\cdot 2^{n}+C_1$ and $18\cdot 2^{n}+C_2$,
for every sufficiently large number $d$ there exists a number of the
form $l(n)$ between $d$ and $15d$. Indeed it is enough to find a
natural number $n$ satisfying $2^n>(d-C_1)/3$ and
$2^n<(15d-C_2)/18$. Such a number $n$ exists for all sufficiently
large $d$ since $15/18>2/3$.

Hence $n^2\log n$ is equivalent to a lower bound for the Dehn
function of $\sss\circ Z$.
\endproof

Lemmas \ref{lb} and \ref{main} show that the Dehn function of the
group $\sss\circ Z$ is equivalent to $n^2\log n$.

Finally the undecidability part of Theorem \ref{thmain} follows
from Lemmas \ref{conj} and \ref{lm1}. Indeed, by these lemmas if
$\LL$ is a non-recursive language then $\sss\circ Z$ has
undecidable conjugacy problem.

\begin{rk}\label{rklast} {\rm An easy example of a group with Dehn function
$n^2\log n$ constructed using the method described in this paper
is the group $\sss_0\circ Z$ where $\sss_0$ has only one
$q$-letter $\{k\}$, one $a$-letter $a$, and two rules $\tau_1$,
$\tau_2$, both equalling $[k\to ak]$ (and their inverses). One can
easily write down an explicit presentation of this group. It is an
HNN extension of a free group of rank 10 with 14 free letters.}
\end{rk}

\begin{rk}\label{rklast1} The method used in this paper allows one to construct groups
with other small Dehn functions. In fact, a similar argument to
the one used in the proof of Lemma \ref{main} proves the following
statement.

\bigskip

{\em  Let $f(n) > 1$ be any function. Let $S$ be an $S$-machine
satisfying the condition of Lemma \ref{width} with $\log_2 h$
replaced by $f(h)$. Then the Dehn function of the group $S$ does not
exceed $n^2f(n)$. }
\end{rk}

Using the technique from \cite{SBR} of converging Turing machines
into $S$-machines, and Remark \ref{rklast1}, one can construct
$S$-machines with many other Dehn functions between $n^2$ and $n^3$
(the fact that the Dehn function of any $S$-machine does not exceed
$n^3$ has been proved in \cite{SBR}).

\begin{minipage}[t]{3 in}
\noindent Alexander Yu. Ol'shanskii\\ Department of Mathematics\\
Vanderbilt University \\ alexander.olshanskiy@vanderbilt.edu\\
http://www.math.vanderbilt.edu/$\sim$olsh\\ and\\ Department of
Higher Algebra, MEHMAT\\
 Moscow State University\\
olshan@shabol.math.msu.su\\
\end{minipage}
\begin{minipage}[t]{3 in}
\noindent Mark V. Sapir\\ Department of Mathematics\\
Vanderbilt University\\
m.sapir@vanderbilt.edu\\
http://www.math.vanderbilt.edu/$\sim$msapir\\
\end{minipage}

\addtocontents{toc}{\contentsline {section}{\numberline {
}References \hbox {}}{\pageref{bibbb}}}

\end{document}